\documentclass[12pt,a4paper,leqno]{amsart}
\usepackage{amssymb}
\advance\textwidth by 2.5cm %%3.6
\advance\oddsidemargin by -1.6cm \advance\evensidemargin by -1.6cm

\input colordvi

\newcommand\new[1]{}
\newcommand\changed[2]{{#1}}
\newcommand\ok[1]{{#1}}

\usepackage{times}

\usepackage{mathrsfs}

\usepackage{verbatim}

\theoremstyle{plain}

\newtheorem{theorem}{Theorem}[section]

\theoremstyle{remark}

\newtheorem{remark}[theorem]{Remark}

\newtheorem{question}[theorem]{Question}

\theoremstyle{plain}

\newtheorem{corollary}[theorem]{Corollary}

\newtheorem{lemma}[theorem]{Lemma}

\newtheorem{proposition}[theorem]{Proposition}

\newtheorem{definition}[theorem]{Definition}

\newenvironment{del}[1]{}{\par\vspace{6pt}}

\numberwithin{equation}{section}

\def\todown{\searrow}

\newcommand{\lb}{\langle}

\newcommand{\rb}{\rangle}

\newcommand\diver{\operatorname{div}}

\newcommand\curl{\operatorname{curl}}

\newcommand\rot{\operatorname{rot}}

\newcommand\id{\operatorname{id}}

\newcommand{\eps}{\varepsilon}

\newcommand{\al}{\alpha}

\newcommand{\la}{\lambda}

\def\Nat{{\mathbb N}}

\def\Rnu{{\mathbb R}}

\def\Rn{{\Rnu^n~}}

\def\Rd{{\Rnu^d~}}

\def\liml{\lim\limits}

\def\suml{\sum\limits}

\def\supl{\sup\limits}

\def\intl{\int\limits}

\def\supl{\mathop{\sup}\limits}

\newcommand\tr{\mathop{tr}}

\begin{document}

\baselineskip 20pt

\title[Duality and the NSEs]{Duality, Vector advection and the Navier-Stokes equations}
\author{Z. Brze{\'z}niak}
\address{Department of Mathematics\\
The University of York\\
Heslington, York YO10 5DD, UK} \email{zb500@york.ac.uk}
\author{M. Neklyudov}
\address{School of Mathematics\\
The University of New South Wales\\
Sydney 2052, Australia} \email{misha.neklyudov@gmail.com}

\date{\today}

\thanks{The research of the second named authour was supported by an ORS award.}
\keywords{Navier-Stokes equations, Feynman Kac formula, Vector
advection} \subjclass{35Q30,60H30,76D05}
\begin{abstract}
In this article we show that three dimensional vector advection
equation is self dual in certain sense defined below. As a
consequence, we infer classical result of Serrin of existence of
strong solution of
Navier-Stokes equation. Also we deduce  %from generalization of Kelvin theorem
Feynman-Kac type formula for solution of the vector advection
equation and show that the formula is not unique i.e. there exist
flows which differ from standard flow along which vorticity is
conserved.
\end{abstract}
\maketitle

The purpose  of this paper is twofold. The first one is to
establish a certain self-duality formula for a vector advection
equation
 in the space $\mathbb{R}^3$. This formula can be understood
as generalization of the helicity invariance  for the Euler
equations , see Corollary \ref{cor:Duality-3} and Remark
\ref{Rem:gen_Kelvin_0}. As a byproduct,   see Corollary \ref{cor:Duality-3}, we give
 a new proof of  the classical result of Serrin \cite{Serrin} about the uniqueness of a
weak solution to the  Navier-Stokes equations (NSEs for short)
\begin{eqnarray}\label{eqn_NSES}
  \frac{\partial u}{\partial t}+(u\nabla)u &=& \nu\triangle u+\nabla p+f \\
  \diver u &=& 0\nonumber\\
  u(0) &=& u_0\nonumber
\end{eqnarray}
satisfying certain additional integrability condition.
\del{$L^r(0,T;\mathbb{L}^s(\Rnu^3))$,
$\frac{2}{r}+\frac{3}{s}=1$,}
 The second one, see Theorem \ref{thm:2DFlowDifferentLaw},  is to establish the
 existence of non-classical flows along with which the circulation of
 the solution of the vector advection equation is conserved in the mean.
 This problem seems to us important because it could potentially
 lead to the new a priori estimates of the solution of vector advection
 equation.

The importance of the vector advection equation stems from  the
fact that it appears in many different areas of hydrodynamics, e.g.
 the vorticity of a strong solution of the
3-dimensional \changed{NSEs}{Navier-Stokes equations} is its solution. Moreover, the major obstacle in
proving  the global existence of a strong solution to the
\changed{NSEs}{Navier-Stokes equations} is the appearance  of the
''vorticity stretching'' term in the vector advection equation. It
is necessary to underline that in a simpler case   of the scalar
advection equation, this conceptually important term is not
present and therefore  the  self-duality and other properties
described in this paper do not hold. Another application of the
vector advection equation is the equation for magnetic field in
MHD equations, see e.g. \cite{Moffat}.

Let us describe briefly  the main contributions of the paper.

In the first part of our paper   we study the following vector advection
equations
\begin{eqnarray}
\nonumber
\frac{\partial F(t,x)}{\partial t}&=&-\nu (A F)(t,x) \\
\label{eqn:DualEquationI}
&-&[\mathrm{P}((v(t,x)\nabla) F(t,x)-\nabla F(t,x) v(t,x))](t,x)+f(t,x),\; x\in \mathbb{R}^d,\\
F(0,x)&=&F_0(x),\; x\in \mathbb{R}^d,\nonumber
\end{eqnarray} where $v:[0,\infty)\times\mathbb{R}^d\to
\mathbb{R}^d$ is a given time-dependent vector field, $d=3$,
$\mathrm{P}$ is the Helmholtz projection onto the divergence free
vector fields and $A$ is the  Stokes operator. As usual by
$\mathbb{H}_{\rm sol}^{k,2}(\mathbb{R}^d)$, $d\in\mathbb{N}$,  we
denote  the space of  all divergence free vector fields that
belong to the Sobolev space $\mathbb{H}^{k,2}(\mathbb{R}^d)$. Let
us denote by $\mathcal{T}_T^{v}$  the  transport operator along
$v$, i.e. $\mathcal{T}_t^{v}F_0=F(t)$, for $t\geq 0$, where $f$ is
the unique solution to problem \eqref{eqn:DualEquationI}. The main
result here is Theorem \ref{thm:Duality} in which we formulate the
following self-duality formula.
\begin{equation}
(\curl F_0, \mathcal{T}_T^{S_{T}v}G_0)_H=(\curl
\mathcal{T}_T^{v}F_0,G_0)_H,\;F_0\in \curl^{-1}(H),\;G_0\in
H,\label{eqn:Duality-2I}
\end{equation}
where  $S_{T}$ is the time reversal operator, i.e.
$(S_{T}v)(t)=-v(T-t)$, $t\in[0,T]$. The self-duality formula \eqref{eqn:Duality-2I}
allows us to deduce certain properties of the operator
$\mathcal{T}_T^v$. In particular in Corollary
\ref{cor:TransportNorm} we show that the
$\mathcal{L}(\mathbb{H}_{\rm sol}^{k,2},\mathbb{H}_{\rm
sol}^{k,2})$-- norm of $\mathcal{T}_T^v$ is equal to its
$\mathcal{L}(\mathbb{H}_{\rm sol}^{1-k,2},\mathbb{H}_{\rm
sol}^{1-k,2})$-- norm. Moreover, in Corollary
\ref{cor:NaturalSpace}, we prove that the space
$\mathcal{L}(\mathbb{H}_{\rm sol}^{\frac{1}{2},2},\mathbb{H}_{\rm
sol}^{\frac{1}{2},2})$ is in a certain sense optimal for
$\mathcal{T}_T^v$.

  The main result in the
second part of the paper, Theorem  \ref{thm:2DFlowDifferentLaw},
 is about a  certain non-classical
Feynman-Kac type  formula for the solutions of the
vector advection equation \eqref{eqn:DualEquationI} in two
dimensions. We show that if the divergence free vector field $v$ is time-independent and sufficiently regular, then  the  stochastic flow of diffeomorphisms of $\mathbb{R}^2$ $X_{s}(t;\cdot),0\leq s\leq
t\leq T$, corresponding  to the following
SDE on $\mathbb{R}^2$, 
\begin{equation}
\label{eqn:flow-4I}
\left\{
\begin{array}{rcl}
d\changed{X_s(t;x)}{\changed{X_s(t;x)}{X_t^s}} &=&
\sqrt{2\nu}\sigma_1(\changed{X_s(t;x)}{X_t^s})\,dW(t), \;0\leq
s\leq t\leq
T,\\
X_s(s;x) &=& x.
\end{array}
\right.
\end{equation}
where, with a function $\phi:\mathbb{R}^2\to \mathbb{R} $   such that\footnote{Such
$\phi$ exists because $\diver v=0$. }  $v=\nabla^{\bot}\phi$,
\begin{displaymath}
\sigma_1 (x) = \left(\begin{array}{cc} \cos\frac{\phi(x)}{\nu} & -\sin\frac{\phi(x)}{\nu}\\
\sin\frac{\phi(x)}{\nu} & \cos\frac{\phi(x)}{\nu}
\end{array}\right), x\in\Rnu^2,
\end{displaymath}
 has  the following   properties: (i) its one-point
motion is a \textit{Brownian Motion} and (ii)  the circulation along it of
the solution of the two dimensional vector
advection equation \eqref{eqn:DualEquationI}, i.e.  with $d=2$,  
is a martingale. 
This
flow seems to be of interest on its own because the stream
function $\phi$ naturally arise in its construction. %This result can be viewed as a generalization of the Kelvin
%circulation Theorem, see \cite{???}. \Red{Precise reference needed.}

The question of the existence of an analogous flow in  the three
dimensional case  remains open, see Question
\eqref{quest:3DdifferentLawCase} for details. \del{This Question seems to us to be
of some importance because its positive solution  would mean that
there is a possibility of extension (in certain sense) of the
notion of stream function to the three dimensional case. But
maybe it is known that such extension is not possible? I am not
sure. Should we delete this paragraph?}

It should be noticed here  that a similar construction does not
work for the scalar advection equation because in this case the
Feynman-Kac type formula depends only upon the law of the flow
itself and not upon the law of the gradient of the flow. Also we
would like to point out that the main obstacle in getting a'priori
estimates for solutions of vector advection equation (in
particular, for vorticity of the solution to the 3-D NSEs) is lack
of an
 estimate for the gradient of the flow. Therefore, in connection
with this result,  a natural question is whether it is possible to
choose  the optimal flow for which gradient is bounded?

The main idea behind our approach to the Feynman-Kac type formula
for solutions of the vector advection equation is that in the case
with viscosity equal to $0$,  the  conservation law of
circulation, known also as Kelvin-Noether Theorem,  holds. In the
case of positive viscosity we are able to
 find an analog of this conservation law. \del{ for vector advection equation}
% We will show  that the circulation of a solution along  certain
%natural flow of diffeomorphisms is conserved in the average.
The Feynman-Kac formula is then an immediate consequence of  that
result. This idea has been used before in the papers \cite{Nekl2}
and \cite{Nekl} (though with quite sketchy proofs). \del{In
particular, in the paper \cite{Nekl2} this idea has been used to
prove that under certain assumptions there are no solutions for
the inverse Cauchy problem for Navier-Stokes equation.}%\Red{Any
%reason to write about this? It is quite disjoint from our
%paper.Ok.}
In the latter paper, see Theorem $5$ and Example $1$,
the Feynman-Kac formula for the solution of  vector advection
equation without incompressibility condition has been derived. A
somewhat similar idea has been also explored independently by
Constantin and Iyer in \cite{Const_Iyer_2008}, but see also
Flandoli et al. \cite{BusnelloFlandoliRomito} for a different
approach. Moreover, Flandoli et al. \cite{BusnelloFlandoliRomito}
proved Feynman-Kac formula for more general systems of parabolic
PDEs. However, we would like to point out that in all of the
articles mentioned
above only the "standard" stochastic flow %\del{
corresponding to the
following SDE
\begin{eqnarray}
d\changed{X_s(t;x)}{X_t^s}) &=&
v(t,\changed{X_s(t;x)}{X_t^s})\,dt+\sqrt{2\nu}\,dW(t),\; t\in
[s,T], \label{eqn:flow-1I}
\\X_s(s;x) &=& x.\nonumber
\end{eqnarray}
%}
\del{Can this equation be  deleted? I don't think so: In this
paragraph we explain the difference between our result and
different other Feynman-Kac formulas.} has been used and,
correspondingly, the problems discussed here does not appear in
their framework.

%\added{
One possible application of Theorem \ref{thm:2DFlowDifferentLaw}
is the extension of Le Jan and  Raimond's
theory of statistical solutions  of the scalar advection equations, see \cite{[LeJanRaimond2002]},  to the
2D vector advection case. Indeed, Le Jan, Raimond theory defines
statistical solution $X_s(t;x)$ of SDE \eqref{eqn:flow-1I}
(corresponding to a solution of scalar advection equation in a
natural way) with velocity $v$ given by
\begin{equation}
dv^i(t,x)=\suml_{k=1}^{\infty}\sigma_k^i(x)\,dW(t)^k,\;x\in\Rnu^n,\; t\geq
0,\; i=1,\ldots,n,\label{eqn:GaussianVelocity}
\end{equation}
where $\sigma_k^i(\cdot)$ are H\"{o}lder continuous and
$\{W(t)^k\}_{k=1}^{\infty}$ is a family of i.i.d. Wiener
processes. In the case of the 2D vector advection,   Theorem
\ref{thm:2DFlowDifferentLaw} implies that we don't need to define
process $X_s(t;x)$ (It is just Brownian motion!). We only need to
show that the \textbf{linear} equation \eqref{eqn:2DGradientFormula}
 for  the gradient of the flow $\nabla X_s(t;x)$ has a
strong solution. At this moment, there appears certain difficulty
with the definition of the right hand side of equation
\eqref{eqn:2DGradientFormula} for irregular vector field $v$  of
the form \eqref{eqn:GaussianVelocity}. We  are of the impression
that the white noise calculus could be of some help here.
%}

Finally, the idea of generalization of the conservation laws has
been extensively studied in physical literature, where it is
called statistical integral of motion or zero mode, see e.g. the
survey \cite[part II.E, p.932]{FalkovichGawedzkiVergasola}, and
references therein.

\textbf{Note:} After we had proved Corollary
\ref{cor:LocMartingaleProperty} we became aware that independently
of us a similar result was proved  recently by Constantin and Iyer
in \cite{Const_Iyer_2008}.

\subsection*{Acknowledgments} We would like to thank T. Komorowski
and B. Go{\l}dys for their useful remarks, in particular to the
former one for informing us about the work by Constantin and Iyer
\cite{Const_Iyer_2008}. The present article derives from work done as
part of the Ph. D. thesis of the second named authour  at the University of
York, while supported by the ORS award, University of
York scholarship and, later, by an ARC Discovery project
DP0558539. The research of the first named author was supported by
an ARC Discovery grant DP0663153.

\section{Notations and hypotheses}\label{sec:notat}
Let $D$ be either $\Rd$ or an open, bounded and connected set in
$\Rd$. In the latter case, we assume that the
 boundary $\Gamma=\partial D$ of $D$ is of $C^3$ class and we denote by
$\overrightarrow{n}$  the outer normal vector field to the
boundary $\Gamma$. We denote by $C^{\infty}(D,\Rn)$ the space of
infinitely differentiable functions from $D$ to $\Rn$ and by
$C_0^{\infty}(D,\Rn)$ the subspace of those functions belonging to
$C^{\infty}(D,\Rn)$ which have a compact support. Finally, let us
denote

$$\mathcal{D}(D)=\{f\in C_0^{\infty}(D,\Rd):\diver f=0\}.$$

For $k\in\mathbb{N}$ and $p\in[1,\infty)$, let $H_0^{k,p}(D,\Rn)$,
respectively $H^{k,p}(D,\Rn)$, be the completion of
$C_0^{\infty}(D,\Rn)$, respectively $C^{\infty}(D,\Rn)$, with
respect to norm
$$
\Vert f\Vert_{k,p}=(\suml_{l=0}^k\suml_{|\al|\leq
l}\intl_{D}|D^{\al}f(x)|_{\Rn}^p\,dx)^{1/p}.
$$
We will use the following notation
\del{\begin{eqnarray*}
H^{k,p}(D)&=&H^{k,p}(D,\Rnu),\\
H_0^{k,p}(D)&=&H_0^{k,p}(D,\Rnu),\\
\mathbb{H}^{k,p}(D)&=&H^{k,p}(D,\Rd),\\
\mathbb{H}_0^{k,p}(D)&=&H_0^{k,p}(D,\Rd),\\
\mathbb{H}^{k}(D)&=&\mathbb{H}^{k,2}(D),\\
L_0^{p}(D)&=&H_0^{0,p}(D,\Rnu),\\
\mathbb{L}^{p}(D)&=&H^{0,p}(D,\Rd),\\
\mathbb{L}_0^{p}(D)&=&H_0^{0,p}(D,\Rd).
\end{eqnarray*}
}

$$
\begin{array}{rlcrlcrlc}
H^{k,p}(D)&=&H^{k,p}(D,\Rnu),&&\quad&
H_0^{k,p}(D)&=&H_0^{k,p}(D,\Rnu),\\
\mathbb{H}^{k,p}(D)&=&H^{k,p}(D,\Rd),&&&
\mathbb{H}_0^{k,p}(D)&=&H_0^{k,p}(D,\Rd),\\
\mathbb{H}^{k}(D)&=&\mathbb{H}^{k,2}(D),&&&
%L_0^{p}(D)&=&H_0^{0,p}(D,\Rnu),\\
%\mathbb{L}_0^{p}(D)&=&H_0^{0,p}(D,\Rd),&&&
\mathbb{L}^{p}(D)&=&H^{0,p}(D,\Rd).
\end{array}
$$

Finally, let us denote
\begin{eqnarray*}
H&=&\{f\in\mathbb{L}^2(D):\diver
f=0,(f\cdot\overrightarrow{n})|_{\Gamma}=0\},\\
V&=&\mathbb{H}_0^{1,2}(D)\cap H.
\end{eqnarray*}

Equipped  with the norm $\Vert \cdot\Vert _{0,2}$, $H$ is a
Hilbert space. Similarly, $V$ is a Hilbert space when equipped
with the norm $\Vert \cdot\Vert_{1,2}$.
 The norms in $H$ and $V$ will be  denoted by $|\cdot|$ and $\Vert \cdot\Vert$. See also \cite[pp. 9-15]{Temam_2001}  for the definition and different
characterizations of the spaces $H$ and $V$.

By $\mathbb{H}_{\rm sol}^{k,p}(D)$ we will denote the  completion
of $\mathcal{D}(D)$ w.r.t. the norm $\Vert \cdot\Vert_{k,p}$.
We will often write $\mathbb{H}_{\rm sol}^{k,p}$ instead of
$\mathbb{H}_{\rm sol}^{k,p}(\Rnu^3)$. We also denote by
$\mathbb{H}_{h,\rm sol}^{k,p}$ the completion of
$\mathcal{D}(\Rnu^3)$ w.r.t. the homogeneous norm 
$$
\Vert
f\Vert_{k,p}^h=(\intl_{\Rnu^3}|\curl^{k}f|_{\Rnu^3}^p\,dx)^{1/p},k\in\Nat,\in[1,\infty).
$$
Let us also denote $\mathbb{H}_{h,\rm
sol}^{-k,2}=(\mathbb{H}_{h,\rm sol}^{k,2})^*,k\in\Nat$ %\Red{I have
%changed $p$ to $2$ since I think we only need these spaces with
%$p=2$. Correct}
and define the spaces with fractional order via
the complex interpolation, i.e.
$$\mathbb{H}_{h,\rm
sol}^{\al,p}=\big[\mathbb{H}_{h,\rm sol}^{[\al],p},\mathbb{H}_{h,\rm
sol}^{[\al]+1,p}\big]_{\al-[\al]},\al\in\Rnu,
$$
where $[\cdot,\cdot]_{\beta}$ is a complex interpolation space of
order $\beta$.

Let $\mathrm{P}:\mathbb{H}^{k,p}(D)\to \mathbb{H}_{\rm
sol}^{k,p}(D)$ be the Helmholtz projection onto the
divergence-free vector fields, see \cite{Fuj_Mor_1977} or
\cite{Temam_2001}. \del{I think there is an older reference!}

%$$D(A_2)=V_2\cap\mathbb{H}^{2,2}(D),A_2=-\mathrm{P}\triangle:D(A_2)\to H,\mbox{$A_2$-Stokes operator.}$$
From now on we consider the case $d=3$. By $\times$ we will denote
the vector product in $\Rnu^3$. We will often use the following
properties of the vector product.
\begin{eqnarray}
(a\times b,c)_{\Rnu^3}=(a,b\times c)_{\Rnu^3}\label{eqn:Wedge-1}\\
|a\times b|_{\Rnu^3}\leq
|a|_{\Rnu^3}|b|_{\Rnu^3}.\label{eqn:Wedge-2}
\end{eqnarray}

We will identify the dual $H^\prime$ with $H$ and so we can assume
that $H\subset V^\prime$. In particular,  $$V\subset H\cong
H^\prime \subset V^\prime$$ is Gelfand triple. We will need the
following results borrowed from the monograph  \cite{Lions} by
Lions and Magenes, see Theorem 3.1, p. 19 and Proposition 2.1, p.
18.
\begin{lemma}\label{lem:ContinuityLem}
Suppose that $\mathcal{V}\subset\mathcal{H}\subset\mathcal{V}'$ is
a Gelfand triple with the duality relation $\langle
\cdot,\cdot\rangle_{\mathcal{V}',\mathcal{V}}$. If $u\in
L^2(0,T;\mathcal{V})$, $u'\in L^2(0,T;\mathcal{V}')$, then $u$ is
almost everywhere equal to a continuous function from $[0,T]$ into
$\mathcal{H}$ and we have the following equality, which holds in
the scalar distribution sense on $(0,T)$:
\begin{equation}
\frac{d}{dt}|u|^2=2\langle u',u\rangle .
\end{equation}
\end{lemma}
As a consequence we have the following result.
\begin{corollary}\label{cor:DiffRule}
If $f,g\in L^2(0,T;\mathcal{V})$ with $f^\prime,g'\in
L^2(0,T;\mathcal{V}')$ then $(f,g)_{\mathcal{H}}$ is almost
everywhere equal to weakly differentiable function and
\begin{equation}
\frac{d}{dt}(f,g)_{\mathcal{H}}=\langle
f^\prime,g\rangle_{\mathcal{V}',\mathcal{V}}+\langle
f,g'\rangle_{\mathcal{V}',\mathcal{V}}.\label{eqn:DiffRule}
\end{equation}
\end{corollary}

We also recall the following result from \cite{Lions}, see Theorem
4.1, p. 238 and Remark 4.3, p. 239
\begin{theorem}\label{thm:LionsPDEexun}
Assume that
\begin{equation}
A\in
L^{\infty}([0,T],\mathcal{L}(\mathcal{V},\mathcal{V}'))\label{eqn:BoundednCond}
\end{equation}
satisfies the following coercivity condition. There exist $\al>0$
and $\lambda\in\Rnu$ such that
\begin{equation}
\langle A(t)u,u\rangle_{\mathcal{V}',\mathcal{V}}\geq\al
|u|_{\mathcal{V}}^2+\lambda
|u|_{\mathcal{H}}^2,\,\,\,u\in\mathcal{V}.\label{eqn:CoercivCond}
\end{equation}
Then for all $u_0\in\mathcal{H}$ and $f\in L^2(0,T;\mathcal{V}')$
the problem
$$
\left \{\begin{array}{lll}
\frac{du}{dt}+Au = f,\\
u(0) = u_0
\end{array}\right.
$$
has unique solution $u\in L^2(0,T;\mathcal{V})$ such that $u'\in
L^2(0,T;\mathcal{V}')$. Moreover, this unique solution $u$
satisfies the following inequality
\begin{equation}
|u|_H^2(t)+\al\intl_0^t|u(s)|_V^2\,ds\leq (1+2\lambda
t)e^{2\lambda
t}(|u_0|_H^2+\frac{1}{4\al}\intl_0^t|f|_{V^\prime}^2\,ds),t\in
[0,T].\label{eqn:AprioriEstim-1}
\end{equation}
\end{theorem}
We will also need the following result.
\begin{proposition}\label{prop:CoercivSemigroup}
Assume that an operator
$A\in\mathcal{L}(\mathcal{V},\mathcal{V}')$ satisfies the
coercivity condition \eqref{eqn:CoercivCond}. {Let us denote
$D(A)=\{x\in \mathcal{H}|Ax\in \mathcal{H}\}$.} Then for all $f\in
L^2(0,T;\mathcal{H})$ and $u_0\in \mathcal{V}$ there exists a
unique solution $u\in L^2(0,T;D(A))\cap C([0,T];\mathcal{V})$ of
the problem:
\begin{eqnarray}
\frac{du}{dt}+\nu Au &=& f,\label{eqn:SParabolEq-1}\\
u(0) &=& u_0\nonumber
\end{eqnarray}
and it satisfies $u'\in L^2(0,T;\mathcal{H})$. Moreover, for a
constant $C=C(\lambda,T,\nu)$ independent of $u_0$ and $f$, such
that
\begin{eqnarray}
|u'|_{L^2(0,T;\mathcal{H})}^2+\nu^2|u|_{L^2(0,T;D(A))}^2\leq
C(|f|_{L^2(0,T;\mathcal{H})}^2+|u_0|_{\mathcal{V}}^2).\label{eqn:AprioriEstim-2}
\end{eqnarray}
\end{proposition}
\begin{proof}[Proof of Proposition \ref{prop:CoercivSemigroup}]
It follows from Theorem 3.6.1 p.76 of \cite{Tanabe} that $-A$
generates an analytic semigroup in $\mathcal{H}$. Therefore, the
existence and the uniqueness of solution $u$ follows from Theorem
3.2 p.22 of \cite{LionsMagenes-2}. It remains to show the
inequality \eqref{eqn:AprioriEstim-2}. let us define a Banach
space $X=\{u\in L^2(0,T;D(A)):u'\in L^2(0,T;H)\}$ and a bounded
linear operator $\mathcal{Q}:X\ni u\mapsto (u(0),u'+Au)\in V\times
L^2(0,T;H))$. Since  $\mathcal{Q}$ is a bijection, according to
the Open Mapping Theorem, there exists the inverse continuous
operator $\mathcal{Q}^{-1}$, i.e.
$\mathcal{Q}^{-1}\in\mathcal{L}(V\times L^2(0,T;H),X)$. Hence the
inequality \eqref{eqn:AprioriEstim-2} follows.

%Multiplying \eqref{eqn:SParabolEq-1} on
%$Au$, integrating it w.r.t. time and using coercivity estimate
%\eqref{eqn:CoercivCond} we get
%\begin{equation}
%\al |u|_{\mathcal{V}}^2(t)+\nu|u|_{L^2(0,T;D(A))}^2\leq
%\frac{1}{4\nu}|f|_{L^2(0,T;\mathcal{H})}^2+\langle Au_0,u_0\rangle +\lambda|u|_{\mathcal{H}}^2(t)
%\end{equation}
%Now the result follows from \eqref{eqn:AprioriEstim-1}.
\end{proof}

\begin{definition}\label{def:CoerciveForm-1}
Let  us define a bilinear form $\tilde{a}:V\times V\to \Rnu$ by
\begin{eqnarray*}
\tilde{a}(u,v)=\suml_{i,j=1}^3\intl_D\nabla_i u^j\nabla_i
v^jdx,u,v\in V.
\end{eqnarray*}
\end{definition}
\begin{lemma}\label{lem:coercivity}
The form $\tilde{a}:V\times V\to \Rnu^1$ is positive, bilinear,
continuous and symmetric.
\end{lemma}
\begin{proof}
Proof is omitted.
\end{proof}
It follows from Lemma \ref{lem:coercivity} and the Lax-Milgram
Theorem that for any $f\in V^\prime$ there exists unique $u\in V$
such that
\begin{equation}
\tilde{a}(u,v)+\lambda (u,v)=\langle f,v\rangle_{V^\prime,V},v\in
V.\label{eqn:Eq-VarForm1}
\end{equation}
\begin{definition}\label{def_Stokes}
Define $A\in\mathcal{L}(V,V^\prime)$ by an identity
$$
\tilde{a}(u,v)=\langle Au,v\rangle_{V,V^\prime},u,v\in V.
$$
\end{definition}
\begin{remark}
The operator $A$ defined above is often called the Stokes
operator.
\end{remark}
\begin{corollary}
The operator $A$ defined in Definition \ref{def_Stokes} is
self-adjoint and positive definite.
\end{corollary}
\begin{proof}
Follows from the symmetry of the form $\tilde{a}$, Theorem 2.2.3,
Remark 2.2.1, p.29 of \cite{Tanabe}.
\end{proof}
\begin{definition}
Let us define trilinear form
$\tilde{b}:C_0^{\infty}(D)\times\mathcal{D}\times\mathcal{D}\to\Rnu$
by
\begin{equation}
\tilde{b}(v,f,\phi)=\langle \mathrm{P}(v\times\curl
f),\phi\rangle_{V^\prime,V},\; (v,f,\phi)\in
C_0^{\infty}(D)\times\mathcal{D}\times\mathcal{D}.
\end{equation}
\end{definition}
\begin{lemma}\label{lem:NLTermEstimates}
For any $\delta$ there exists $C_{\delta}>0$ such that for all
$\eps>0$ and all $(v,f,\phi)\in
C_0^{\infty}(D)\times\mathcal{D}\times\mathcal{D}$,
\begin{eqnarray}
\label{eqn:NlinEstimate1} |\tilde{b}(v,f,\phi)|^2& \leq
&|f|_V^2|\phi|_V^2(\eps^{1+\delta/3}+
\frac{C_{\delta}}{\eps^{1+3/\delta}}|v(t)|_{\mathbb{L}^{3+\delta}(D)}^{2+\frac{6}{\delta}}),\\
\label{eqn:NlinEstimate2} |\tilde{b}(v,f,\phi)|&\leq&
\frac{1}{2}\Vert f\Vert_{V}^2+\frac{1}{2}(\eps^{1+\delta/3}\Vert
\phi\Vert _{V}^2+
\frac{C_{\delta}}{\eps^{1+3/\delta}}|v(t)|_{\mathbb{L}^{3+\delta}(D)}^{2+\frac{6}{\delta}}|\phi|_{H}^2).
\end{eqnarray}
Moreover, if we assume that $f\in D(A)$, then for any $\phi\in V$
the following inequality holds
\begin{equation}
|\tilde{b}(v,f,\phi)|^2\leq |\phi|_H^2(\eps^{1+\delta/3} \Vert
f\Vert
_{D(A)}^2+\frac{C_{\delta}}{\eps^{1+3/\delta}}|v|_{\mathbb{L}^{3+\delta}(D)}^{2+\frac{6}{\delta}}|f|_{V}^2)
\end{equation}
\end{lemma}
To prove Lemma \ref{lem:NLTermEstimates} we will need the
following auxiliary result.
\begin{lemma}\label{lem:EstimGag-NirType}
For any $\delta$ there exists $C_{\delta}>0$ such that for all
$\eps>0$
\begin{equation}
\Vert f\times g\Vert _{\mathbb{L}^2(D)}^2\leq \eps^{1+\delta/3}
\Vert f\Vert
_{V}^2+\frac{C_{\delta}}{\eps^{1+3/\delta}}|g|_{\mathbb{L}^{3+\delta}(D)}^{2+\frac{6}{\delta}}|f|_{H}^2,\,f\in
V,g\in H.\label{eqn:EstimGag-NirType}
\end{equation}
\end{lemma}
\begin{proof}[Proof of Lemma \ref{lem:EstimGag-NirType}]
Let us $p=3-\frac{2\delta}{1+\delta}$, $q=\frac{3+\delta}{2}$,
$\theta=\frac{3}{3+\delta}$. Then $\frac{1}{p}+\frac{1}{q}=1$ and
therefore by the inequality \eqref{eqn:Wedge-2}, \del{second
follows from} the   H\"{o}lder inequality, the \del{third follows
from} Gagliardo-Nirenberg inequality (see Theorem 9.3, p.24 in
\cite{Friedman}) and \del{fourth inequality follows} from the
Young inequality we infer the following train of inequalities
\begin{eqnarray*}
\Vert f\times g\Vert _{\mathbb{L}^2(D)}^2 & \leq
&\intl_{D}|f|^2|g|^2\,dx\leq
|f|_{\mathbb{L}^{2p}(D)}^2|g|_{\mathbb{L}^{2q}(D)}^2\\
\leq (\Vert f\Vert
_{V}^{\theta}|f|_H^{1-\theta})^2|g|_{\mathbb{L}^{2q}(D)}^2 &\leq&
\eps^{1+\delta/3} \Vert f\Vert
_{V}^2+\frac{C_{\delta}}{\eps^{1+3/\delta}}|g|_{\mathbb{L}^{3+\delta}(D)}^{2+\frac{6}{\delta}}|f|_{H}^2.
\end{eqnarray*}
\del{where the first inequality follows from }
\end{proof}
\begin{proof}[Proof of Lemma \ref{lem:NLTermEstimates}]
Let us fix $(v,f,\phi)\in
C_0^{\infty}(D)\times\mathcal{D}\times\mathcal{D}$. Then by
equality \eqref{eqn:Wedge-1}, \del{second inequality follows from}
and  Lemma \ref{lem:EstimGag-NirType} we have
\begin{eqnarray}
\label{eqn:Est-1} |\tilde{b}(v,f,\phi)|^2 &= &|\langle v(t)\times
\phi,\curl f\rangle_{V^\prime,V}|^2
\\
\leq |\curl f|_H^2|v(t)\times \phi|_H^2 &\leq & \Vert f\Vert
_{V}^2(\eps^{1+\delta_0/3}\Vert \phi\Vert _{V}^2+
\frac{C_{\delta_0}}{\eps^{1+3/\delta_0}}|v(t)|_{\mathbb{L}^{3+\delta_0}(D)}^{2+\frac{6}{\delta_0}}|\phi|_{H}^2)\nonumber\\
 &\leq& |f|_V^2|\phi|_V^2(\eps^{1+\delta_0/3}+
\frac{C_{\delta_0}}{\eps^{1+3/\delta_0}}|v(t)|_{\mathbb{L}^{3+\delta_0}(D)}^{2+\frac{6}{\delta_0}}).\nonumber
\end{eqnarray}
\del{where equality \eqref{eqn:Est-1} follows from }

Similarly,
\begin{eqnarray}
|\tilde{b}(v,f,\phi)| &=& |\langle v(t)\times \phi,\curl
f\rangle_{V^\prime,V}|\leq
|\curl f|_H|v(t)\times \phi|_H\\
\nonumber \leq \frac{1}{2}\Vert f\Vert
_{V}^2+\frac{1}{2}|v(t)\times \phi|_H^2 &\leq& \frac{1}{2}\Vert
f\Vert _{V}^2+\frac{1}{2}(\eps^{1+\delta_0/3}\Vert \phi\Vert
_{V}^2+
\frac{C_{\delta_0}}{\eps^{1+3/\delta_0}}|v(t)|_{\mathbb{L}^{3+\delta_0}(D)}^{2+\frac{6}{\delta_0}}|\phi|_{H}^2),
\nonumber
\end{eqnarray}
and
\begin{eqnarray}
|\tilde{b}(v,f,\phi)|^2&=&|\langle v(t)\times\curl
f,\phi\rangle_{V^\prime,V}|^2\leq
|\phi|_H^2|v(t)\times\curl f|_H^2\\
\nonumber &\leq & |\phi|_H^2(\eps^{1+\delta/3} |f|_{D(A)}^2
+\frac{C_{\delta}}{\eps^{1+3/\delta}}|v|_{\mathbb{L}^{3+\delta}(D)}^{2+\frac{6}{\delta}}|f|_{V}^2).
\end{eqnarray}
\end{proof}
Fix $\delta_0>0$. It follows from inequality
\eqref{eqn:NlinEstimate1} that \del{for any $\delta_0>0$} the
trilinear form $\tilde{b}$ is continuous with respect to the
$L^{3+\delta_0}(D)\times V\times V$ topology. Therefore, \del{for
any $\delta_0>0$} there exist continuous trilinear form
$b:L^{3+\delta_0}(D)\times V\times V\to\Rnu$ such that
$$
b(\cdot,\cdot,\cdot)|_{C_0^{\infty}(D)\times\mathcal{D}\times\mathcal{D}}=\tilde{b}.
$$
Moreover,
\begin{equation}
b(v,f,\phi)=-(v\times \phi,\curl f)_H,\; (v,f,\phi)\in
\ok{L^{3+\delta_0}(D)\times V\times V}.\label{eqn:bFormFormula-1}
\end{equation}
Indeed, the form on the left hand side of equality
\eqref{eqn:bFormFormula-1} is equal to the form on the right hand
side of equality \eqref{eqn:bFormFormula-1} for $(v,f,\phi)\in
C_0^{\infty}(D)\times\mathcal{D}\times\mathcal{D}$ and both forms
are continuous in $L^{3+\delta_0}(D)\times V\times V$.

%Indeed, it follows from inequality \eqref{eqn:NlinEstimate1} that
%$b$ is continuous in topology of $L^{3+\delta_0}(D)\times V\times V$
\begin{definition}
Let us define a bilinear operator $B:L^{3+\delta_0}(D)\times V\to
V^\prime$ by
$$
\langle B(v,f),\phi\rangle_{V^\prime,V}=b(v,f,\phi),v\in
L^{3+\delta_0}(D),f\in V,\phi\in V.
$$
\end{definition}
\begin{corollary}\label{cor:NonlinEstimates}
Assume that $\delta_0>0$.  Then there exists a constant
$C_{\delta_0}>0$ independent  such that
\begin{equation}
|B(v,f)|_{V^\prime}^2\leq \Vert f\Vert _V^2(\eps^{1+\delta_0/3}+
\frac{C_{\delta_0}}{\eps^{1+3/\delta_0}}|v(t)|_{\mathbb{L}^{3+\delta_0}(D)}^{2+\frac{6}{\delta_0}}),\;
(v,f)\in L^{3+\delta_0}(D)\times V.
\end{equation}
Moreover, if $(v,f)\in L^{3+\delta_0}(D)\times D(A)$ then
$B(v,f)\in H$ and
\begin{equation}
|B(v,f)|_{H}^2\leq (\eps^{1+\delta_0/3} \Vert f\Vert
_{D(A)}^2+\frac{C_{\delta_0}}{\eps^{1+3/\delta_0}}|v|_{\mathbb{L}^{3+\delta_0}(D)}^{2+\frac{6}{\delta}}|f|_{V}^2)\label{eqn:NlinEstimate3}
\end{equation}
\end{corollary}
\begin{proof}[Proof of Corollary \ref{cor:NonlinEstimates}]
Proof immediately follows from Lemma \ref{lem:NLTermEstimates}.
\end{proof}

\section{Duality}\label{sec:Exist&Uniqueness}
Assume that $F_0\in H, f\in L^2(0,T;V^\prime)$. We consider the
following two problems:
\begin{eqnarray}
\frac{\partial F}{\partial t}&=&-\nu A F-B(v(t),F)+f,\;\label{eqn:DualEquation}\\
F(0)&=&F_0,\label{in_v}\\
\frac{\partial G}{\partial t}&=&-\nu A G-\curl{(v(t)\times
G)}+f,\;\label{eqn:DualEquation-2}\\
G(0)&=G_0.\label{in_v-2}
\end{eqnarray}

\del{Changed.}
\begin{definition}
We will call an element $F$, respectively $G$,  of $L^2(0,T;V)\cap
L^{\infty}(0,T;H)\cap C([0,T];H_{\rm w})$, where $H_{\rm w}$ is
equal to  $H$ endowed with the weak topology,  a solution of
problem (\ref{eqn:DualEquation}-\ref{in_v}), resp.
(\ref{eqn:DualEquation-2}-\ref{in_v-2}),  iff $F$, resp. $G$,
satisfies equation (\ref{eqn:DualEquation}), resp. problem
\eqref{eqn:DualEquation-2},  in the distribution sense and $F$,
resp. $G$,  satisfies \eqref{in_v}, resp. \eqref{in_v-2} as
elements of $C([0,T];H_{\rm w})$.
\end{definition}
In the next two Propositions  we will deal with the  existence and
regularity results for solutions of
(\ref{eqn:DualEquation}-\ref{in_v}) and
(\ref{eqn:DualEquation-2}--\ref{in_v-2}). These results are
probably known, but since we have been unable to find  them (the
form we need)  in the literature, we have decided to present them
for the sake of the completeness of the paper.

\begin{proposition}\label{prop:ExistUniqSolution}
Suppose that $(F_0,f)\in H\times L^2(0,T;V^\prime)$ and
\begin{equation}
v\in\bigcup_{\delta_0>0}
L^{2+\frac{6}{\delta_0}}(0,T;\mathbb{L}^{3+\delta_0}(D)).\label{eqn:IntegrCond-1}
\end{equation}
Then
\begin{trivlist}\item[(i)] there exists the unique solution $F$ of
problem (\ref{eqn:DualEquation}-\ref{in_v}) and for each
$\delta_0>0$ there exists a constant
$K_1=K_1(|v|_{L^{2+\frac{6}{\delta_0}}(0,T;\mathbb{L}^{3+\delta_0}(D))},\nu)$
independent of $(F_0,f)$  such that $F$ satisfies the following
inequality
\begin{eqnarray}
\nonumber |F(t)|_H^2&+&\nu\intl_0^t\Vert F(s)\Vert _{V}^2\,ds \\
&\leq&
K_1\big(\intl_0^t|v(s)|_{\mathbb{L}^{3+\delta_0}}^{2+6/\delta_0}\,ds,\nu\big)
\Big(|F_0|_H^2+\frac{C}{\nu}\intl_0^t|f(s)|_{V^\prime}^2\,ds\Big),\;
t\in [0,T].\label{eqn:EnergyInequality-3}
\end{eqnarray}
Furthermore, $F\in C([0,T],V^\prime)$ and  $F'\in
L^{1+\frac{3}{2\delta_0+3}}(0,T;V^\prime)$. Moreover, if $v$
satisfies the following, stronger  than \eqref{eqn:IntegrCond-1},
condition
\begin{equation}
v\in L^{\infty}(0,T;\mathbb{L}^{3+\delta_0}(D))\mbox{ for some
}\delta_0>0,\label{eqn:IntegrCond-2}
\end{equation}
then $F^\prime\in L^2(0,T;V^\prime)$.

\item[(ii)] If in addition $(F_0,f)\in V\times L^2(0,T;H)$ and the
condition \eqref{eqn:IntegrCond-2} is satisfied, then $F\in
C([0,T],V)\cap L^2(0,T;D(A))$.
\item[(iii)]%\Blue{I still don't like this. One needs to spend a
%bit more time and write down those formulas for $F^{(k)}(0)$ and
%instead of making a meaningless assumption that $F^{(k)}(0)\in V$.
%}
Assume that $n\in\mathbb{N}$. Suppose $f^{(n)}\in L^2(0,T;H)$,
there exists $\delta_0>0$ such that $v\in
C^{n-1}(0,T;\mathbb{L}^{3+\delta_0}(D))$, $v^{(n)}\in
L^{\infty}(0,T;\mathbb{L}^{3+\delta_0}(D))$  and $g_k\in V$, for
$k=0,1,\ldots,n$, where sequence $\{g_k\}_{k=0}^{\infty}$ is
defined by formula\footnote{It is easy to see that formally system
(\ref{eqn:DualEquation}-\ref{in_v}) uniquely defines $F^{(k)}(0)$.
Indeed, if we formally put $t=0$ in the system we immediately get
expression for $F'(0)$ through known parameters. Similarly,
differentiating equation \ref{eqn:DualEquation} w.r.t. time we get
recurrent formula \eqref{eqn:RecurrenceFormula-1} for
$F^{(k)}(0),k\in\Nat$. So, the condition $g_k\in V$ is
compatibility condition.}
\begin{eqnarray}
g_m &=& -\nu A g_{m-1}-\suml_{k=0}^{m-1}B(v^{(m-k-1)})(0),g_k)+f^{m-1}(0), m=1,\ldots,n\label{eqn:RecurrenceFormula-1}\\
g_0 &=& F_0.\nonumber
\end{eqnarray}
Then $F\in C^n([0,T],V)$.
\end{trivlist}
\end{proposition}

\begin{remark}
We should notice that on the one hand, our class
$\bigcup_{\delta>0}L^{2+\frac{6}{\delta_0}}(0,T;\mathbb{L}^{3+\delta_0}(D))$
is the Serrin regularity class. Indeed, if
$r=2+\frac{6}{\delta_0}$, $s=3+\delta_0$ then
$\frac{2}{r}+\frac{3}{s}=1$.  Therefore, any weak solution of the
\changed{NSEs}{Navier-Stokes equations} belonging to this class is
a strong solution. On the other hand, we have been unable to prove
that under the assumption \eqref{eqn:IntegrCond-1} a solution $F$
of problem (\ref{eqn:DualEquation}-\ref{in_v}) is such that
$F^\prime\in L^2(0,T,V^\prime)$. A problem that arises  here is
similar to the problem wether a weak solution $u$ of the
\changed{NSEs}{Navier-Stokes equations}, see \cite{Temam_2001}, p.
191 Problem 3.2 and Theorem 3.1, satisfies  $u'\in
L^2(0,T;V^\prime)$.
\end{remark}

For the second equation we have:
\begin{proposition}\label{prop:ExistUniqSolution-2}
Suppose that a time dependent vector field  $v$ satisfies the
assumption \eqref{eqn:IntegrCond-1}. Then
\begin{trivlist}\item[(i)] for every $(F_0,f)\in H\times L^2(0,T;V^\prime)$  there exists unique solution $G$ of the problem
(\ref{eqn:DualEquation-2}-\ref{in_v-2})  such that $G'\in
L^2(0,T;V^\prime)$ and for each $\delta_0>0$ there exits a
constant
$K_1\left(\intl_0^t|v(s)|_{\mathbb{L}^{3+\delta_0}}^{2+6/\delta_0}\,ds,\nu\right)$
such that  $G$ satisfies  the following inequality for $t\in
[0,T]$,
\begin{eqnarray}\nonumber
|G(t)|_H^2&+&\nu\intl_0^t\Vert G(s)\Vert _{V}^2\,ds\\
&\leq&
K_1\Big(\intl_0^t|v(s)|_{\mathbb{L}^{3+\delta_0}}^{2+6/\delta_0}\,ds,\nu\Big)
\Big(|G_0|_H^2+\frac{C}{\nu}\intl_0^t|f(s)|_{V^\prime}^2\,ds\Big),\;
t\in [0,T].\label{eqn:EnergyInequality-2}
\end{eqnarray}
\item[(ii)] If in addition  $v\in L^2(0,T,V)$ and $(F_0,f)\in
V\times L^2(0,T;H)$,  then the solution $G$ from part (i)
satisfies $G\in
C([0,T],V)\cap L^2(0,T;D(A))$. \item[(iii)] Fix $n\in \Nat$. If %\Blue{You need to make changes
%similar to those I did in Proposition 2.2. Done.}
$f^{(n)}\in L^2(0,T;H)$, there exists $\delta_0>0$ such that $v\in
C^{n-1}(0,T;\mathbb{L}^{3+\delta_0}(D))$, $v^{(n)}\in
L^{\infty}(0,T;\mathbb{L}^{3+\delta_0}(D))\cap L^2(0,T,V)$ and
$l_k\in V$, for $k=0,1,\ldots,n$, where sequence
$\{l_k\}_{k=0}^{\infty}$ is defined by formula
\begin{eqnarray}
l_m &=& -\nu A l_{m-1}-\suml_{k=0}^{m-1}\curl(v^{(m-k-1)}(0)\times l_k)+f^{m-1}(0), m=1,\ldots,n\label{eqn:RecurrenceFormula-2}\\
l_0 &=& G_0.\nonumber
\end{eqnarray}
Then $G\in C^n([0,T],V)$.
%\footnote{The value of $G^{(k)}(0)$ in
%condition $G^{(k)}(0)\in V,k\in\Nat$ is defined here in the same
%way as $F^{(k)}(0),k\in\Nat$ was defined in the Proposition
%\ref{prop:ExistUniqSolution}. See footnote on the previous page.}
\end{trivlist}
\end{proposition}
%\begin{remark}
%Notice that in the case of (\ref{eqn:DualEquation-2}-\ref{in_v-2}) we have
%$\diver{G}=0,t\in [0,T]$ in a weak sense and therefore, Proposition
%\ref{prop:ExistUniqSolution-2} is correct. Indeed, applying $\diver$
%operator to equation (\ref{eqn:DualEquation-2}-\ref{in_v-2}) and denoting
%$M=\diver{G}$ we get
%$$
%M'=\diver f, M(0)=0,
%$$
%and therefore $M=0,t\in [0,T]$.
%\end{remark}
\begin{corollary}\label{cor:InfReg-1}%\Blue{Comparison has to be done outside the proof. Make one more remark, %please.Done.}
Assume that $F_0\in H$,  $f,v\in C^\infty([0,T];H)$. If for each
$k\in\mathbb{N}$,  $v^{(k)}$ satisfies the condition
\eqref{eqn:IntegrCond-1},  then the solution of the problem
 (\ref{eqn:DualEquation}-\ref{in_v}) is in $C^{\infty}((0,T]\times
D)$.
\end{corollary}
\begin{proof}[Proof of Corollary \ref{cor:InfReg-1}]
It follows from Remark 3.2, p.90 in \cite{Temam_1981}.
\end{proof}
Similarly for the problem  (\ref{eqn:DualEquation-2}-\ref{in_v-2})
we have
\begin{corollary}\label{cor:InfReg-2}
Assume that $F_0\in H$, $f,v\in C^\infty([0,T];H)$.  If for each $
k\in\Nat$ $v^{(k)}$ satisfies the condition
\eqref{eqn:IntegrCond-2}, then the solution of the problem
(\ref{eqn:DualEquation-2}-\ref{in_v-2}) is in
$C^{\infty}((0,T]\times D)$.
\end{corollary}
\begin{proof}[Proof of Corollary \ref{cor:InfReg-2}]
It follows from Remark 3.2, p.90 in \cite{Temam_1981}.
\end{proof}
The main result in this section is the following.
\begin{theorem}\label{thm:Duality}
Suppose that $F_0\in H$, $G_0\in H$ and  $v\in
\bigcup_{\delta_0>0}
L^{2+\frac{6}{\delta_0}}(0,T;\mathbb{L}^{3+\delta_0}(D))$. Let $F$
and $G$ be solutions of respectively problems \eqref{eqn:DualEq-1}
and \eqref{eqn:DualEq-2} below.
\begin{eqnarray}
\frac{\partial F}{\partial t} &=& -\nu A F-B(v(t),F), \;t\in(0,T),\label{eqn:DualEq-1}\\
F(0,\cdot) &=& F_0,\nonumber\\
\frac{\partial G}{\partial t} &=& -\nu A G+\curl{(v(T-t)\times G)},\; t\in(0,T),\label{eqn:DualEq-2}\\
G(0,\cdot) &=& G_0.\nonumber
\end{eqnarray}
Then,  the following identity holds
\begin{equation}
(F(t),G(T-t))_{H}=(F(0),G(T))_{H},\;t\in [0,T].\label{eqn:Duality}
\end{equation}
\end{theorem}
From now on we will only consider the case $D=\Rnu^3$. We notice
that now  if $F$ is a solution of the problem
(\ref{eqn:DualEquation}-\ref{in_v}) with data  $(F_0,f,v)$, then
$\curl F$ is a solution of the problem
(\ref{eqn:DualEquation-2}-\ref{in_v-2}) with data  $(\curl
F_0,\curl f,v)$.
\begin{definition}\label{def:VectorTransportOperator}
Let $\mathcal{T}_T^v:H\to H$ be the vector transport operator
defined by $\mathcal{T}_T^v (F_0)=F(T)$, where $F$ is the unique
solution of the problem
\eqref{eqn:DualEq-1} with data $(F_0,v)$. %and $v\in
%L^{2+\frac{6}{\delta_0}}(0,T;\mathbb{L}^{3+\delta_0}(\Rnu^3))$.
\end{definition}
Define also the  time reversal operator
$$
S_{T}:\bigcup_{\delta_0>0}
L^{2+\frac{6}{\delta_0}}(0,T;\mathbb{L}^{3+\delta_0}(D))\to
\bigcup_{\delta_0>0}
L^{2+\frac{6}{\delta_0}}(0,T;\mathbb{L}^{3+\delta_0}(D))
$$
by $(S_{T}v)(t)=-v(T-t)$. \del{Changed.} Then from Theorem
\ref{thm:Duality} we infer that
\begin{corollary}\label{cor:Duality-2}
Assume that $F_0\in V$, $G_0\in H$ and  $v\in \bigcup_{\delta_0>0}
L^{2+\frac{6}{\delta_0}}(0,T;\mathbb{L}^{3+\delta_0}(\Rnu^3))$.
Then the following duality relation holds,
\begin{equation}
(\curl F_0, \mathcal{T}_T^{S_{T}v}G_0)_H=(\curl
\mathcal{T}_T^{v}F_0,G_0)_H.\label{eqn:Duality-2}
\end{equation}
\end{corollary}
%\begin{remark}
%This Corollary can be used to define operator $\mathcal{T}_T^v$ on
%the functions from Sobolev spaces with negative index.
%\Red{Unclear. Do we need this anyhow? If yes, we need to be more
%precise here.Changed. I decided to get rid of this remark.}
%\end{remark}
\begin{corollary}\label{cor:TransportNorm}
Assume that $v$ satisfies the assumption
\eqref{eqn:IntegrCond-1} %or that there exists the unique solution $F\in
%L^{\infty}(0,T;H)$ of the problem  \eqref{eqn:DualEq-1} with the initial
%data $F_0\in \mathcal{D}(\Rnu^3)$ \Red{One fixed $F_0$ or any? What about $G_0$ in \eqref{eqn:Duality-2}?}
such that duality relation \eqref{eqn:Duality-2} holds. Then
\begin{equation}
\Vert \mathcal{T}_T^v\Vert _{\mathcal{L}(\mathbb{H}_{h,\rm
sol}^{\al,2},\mathbb{H}_{h,\rm sol}^{\al,2})}= \Vert
\mathcal{T}_T^{S_{T}v}\Vert _{\mathcal{L}(\mathbb{H}_{h,\rm
sol}^{1-\al,2},\mathbb{H}_{h,\rm sol}^{1-\al,2})},\;\al\in [0,1].
\end{equation}
\end{corollary}
\begin{proof}[Proof of Corollary \ref{cor:TransportNorm}]
Because $\mathbb{H}_{h,\rm sol}^{\al,2}$ is the complex
interpolation space between $\mathbb{H}_{h,\rm
sol}^{0,2}=\mathbb{L}_{sol}^2$ and $\mathbb{H}_{h,\rm sol}^{1,2}$
of order $\al$, it is enough to consider the cases
$\al\in\{0,1\}$. Furthermore, we can restrict ourselves to the
case of $\al=0$ because $S_{T}\circ S_{T}=\id$.

 From equality \eqref{eqn:Duality-2} it follows that
\begin{eqnarray*}
\Vert \mathcal{T}_T^v\Vert
_{\mathcal{L}(\mathbb{L}_{sol}^2,\mathbb{L}_{sol}^2)}=\supl_{\phi,\psi\in\mathcal{D}(\Rnu^3)
}\frac{|\langle \mathcal{T}_T^v\phi,\psi\rangle |}{\Vert \phi\Vert _{\mathbb{L}_{sol}^2}\Vert \psi\Vert _{\mathbb{L}_{sol}^2}}=\\
\supl_{\phi,\psi\in\mathcal{D}(\Rnu^3) }\frac{|\langle
\curl
\mathcal{T}_T^v\phi,\curl^{-1}\psi\rangle |}{\Vert \phi\Vert _{\mathbb{L}_{sol}^2}\Vert \psi\Vert _{\mathbb{L}_{sol}^2}}=\\
\supl_{\phi,\psi\in\mathcal{D}(\Rnu^3)
}\frac{|\langle \curl\phi,\mathcal{T}_T^{S_{T}v}\curl^{-1}\psi\rangle |}{\Vert \phi\Vert _{\mathbb{L}_{sol}^2}\Vert \psi\Vert _{\mathbb{L}_{sol}^2}}=\\
\supl_{\phi,\psi\in\mathcal{D}(\Rnu^3)
}\frac{|\langle \phi,\mathcal{T}_T^{S_{T}v}\psi\rangle |}{\Vert \phi\Vert _{\mathbb{H}_{h,\rm sol}^{-1,2}}\Vert \psi\Vert _{\mathbb{H}_{h,\rm sol}^{1,2}}}=\\
\Vert \mathcal{T}_T^{S_{T}v}\Vert _{\mathcal{L}(\mathbb{H}_{h,\rm
sol}^{1,2},\mathbb{H}_{h,\rm sol}^{1,2})}
\end{eqnarray*}
\end{proof}
\begin{definition}
By $X_{\alpha}$ we denote the class of all functions
$u:[0,\infty)\times\Rnu^3\to\Rnu^3$ satisfying the following three
conditions. %\Blue{Changed. I put $T=\infty$ here.}
\begin{trivlist}
\item[(i)] $u\in L_{loc}^{\infty}([0,\infty);H)$. \item[(ii)]  For
all $t\in[0,\infty)$ there exists a unique solution of equation
\eqref{eqn:DualEq-1} with parameters $u^t=u|_{[0,t]}$ and
$v^t=S^t(u|_{[0,t]})$. Furthermore, the duality relation
\eqref{eqn:Duality-2} with the vector field $v$ replaced by the
vector field  $u^t$ holds. \del{Changed.} \item[(iii)] For every
$t\in [0,\infty)$, $\mathcal{T}_t^{u^t}
\in\mathcal{L}(\mathbb{H}_{h,\rm sol}^{\alpha,2},\mathbb{H}_{h,\rm
sol}^{\alpha,2})$. %\Blue{I replaced $\mathcal{T}_t^{u}$ by
%$\mathcal{T}_t^{u^t}$. }
\end{trivlist}
\end{definition}
Then the following result follows from Corollary
\ref{cor:TransportNorm}
\begin{corollary}\label{cor:NaturalSpace}
Assume that $\al\in [0,1]$. Then $X_{\alpha}=X_{1-\alpha}\subset
X_{\frac{1}{2}}$ and  the space $X_{\frac{1}{2}}$ is invariant
with respect to scalings $\Psi_{\la}$, $\la\in (0, 1]$, where
$(\Psi_{\la}u)(t,x)=\la u(\la^2 t,\la x)$,
 $t\in [0,\infty)$, $x\in\Rnu^3$.\\
 \del{Changed (Put $T=\infty$).}
\end{corollary}
\begin{proof}[Proof of Corollary \ref{cor:NaturalSpace}]
Property $X_{\alpha}=X_{1-\alpha}$ is a direct consequence of
Corollary \ref{cor:TransportNorm} and the definition of\,\,
$X_{\alpha}$. We will show that $X_{\alpha}\subset
X_{\frac{1}{2}}$. Let $u\in X_{\alpha}$. Then for all $t\geq 0$,
$$
\mathcal{T}_t^u\in \mathcal{L}(\mathbb{H}_{h,\rm
sol}^{\al,2},\mathbb{H}_{h,\rm sol}^{\al,2}), \mathcal{T}_t^u\in
\mathcal{L}(\mathbb{H}_{h,\rm sol}^{1-\al,2},\mathbb{H}_{h,\rm
sol}^{1-\al,2}).
$$
Indeed, it follows by definition of $X_{\al}$ that
$$
|\mathcal{T}_t^u|_{ \mathcal{L}(\mathbb{H}_{h,\rm
sol}^{1-\al,2},\mathbb{H}_{h,\rm
sol}^{1-\al,2})}=|\mathcal{T}_t^{S^t(u|_{[0,t]})}|_{
\mathcal{L}(\mathbb{H}_{h,\rm sol}^{\al,2},\mathbb{H}_{h,\rm
sol}^{\al,2})},t\in [0,\infty).
$$
Therefore, by the Interpolation Theorem, see  \cite[Theorems
1.9.4, p. 59 and 1.15.3, p. 103]{Triebel},  we have that
$$
\mathcal{T}_t^u\in \mathcal{L}([\mathbb{H}_{h,\rm
sol}^{\al,2},\mathbb{H}_{h,\rm
sol}^{1-\al,2}]_{1/2},[\mathbb{H}_{h,\rm
sol}^{\al,2},\mathbb{H}_{h,\rm sol}^{1-\al,2}]_{1/2}),t\in
[0,\infty),
$$
i.e.
$$
\mathcal{T}_t^u\in \mathcal{L}(\mathbb{H}_{h,\rm
sol}^{\frac{1}{2},2},\mathbb{H}_{h,\rm sol}^{\frac{1}{2},2}),t\in
[0,\infty).
$$
Third property follows from identity
$$
\mathcal{T}_t^{\Psi_\la(u)}\Psi_{\la}(F_0)=\Psi_{\la}(\mathcal{T}_t^{u}F_0),t\in
[0,\infty)
$$
and boundedness of scaling operators $\Psi_{\la}$ and
$\Psi_{\la}^{-1}=\Psi_{\frac{1}{\la}}$ in $\mathbb{H}_{h,\rm
sol}^{\frac{1}{2},2}$.
\end{proof}
%Consider now a family of pairs of  spaces $(\mathbb{H}_{h,\rm
%sol}^{\al,2},X_{\al}),\al\in [0,1]$. Then Corollary
%\ref{cor:NaturalSpace} means  that the pair $(\mathbb{H}_{h,\rm
%sol}^{\frac{1}{2},2},X_{\frac{1}{2}})$ is optimal with respect to
%the  vector transport operator $\mathcal{T}_t^u$ in the following
%sense.  For $\al\in[0,1]$ let ${\tilde X}_{\al}$ be the set of all
%$u\in X_{\al}$ such that $
%\mathcal{T}_t^u\in\mathcal{L}(\mathbb{H}_{h,\rm
%sol}^{\al,2},\mathbb{H}_{h,\rm sol}^{\al,2}) $ \Blue{What about a quantifier for $t$. It seems to me that according to the definition ${\tilde X}_{\al}=X_{\al}$ }. Then ${\tilde
%X}_{\frac12}$ is the biggest set from all of this family. \Red{I
%have made some changes here.}
%\Blue{My suggestion is to remove the whole paragraph here. I did
%it.}

%We would like to point out that
%\begin{corollary}If $u$ is a weak solution of Navier-Stokes
%and $u\in X_{\frac{1}{2}}$ then $u$ is a strong solution of
%Navier-Stokes equation.
%\end{corollary}

The first part of our next result is the classical  result of
Serrin-Prodi- Ladyzhenskaya (\cite{Serrin,Prodi,Ladyzhenskaya}).
But the second part, i.e. inequalities \eqref{eqn:Duality-3} and
\eqref{eqn:CurlEstimate-1} are new.
%which in our context means reads
%$\cupl_{\delta_0>0}L^{2+\frac{6}{\delta_0}}(0,T;\mathbb{L}^{3+\delta_0}(\Rnu^3))\subset
%X_1=X_0$. \Red{I do not understand the logic here. We have proved before that $\cupl_{\delta_0>0}L^{2+\frac{6}{\delta_0}}(0,T;\mathbb{L}^{3+\delta_0}(\Rnu^3))\subset
%X_1=X_0$ haven't we? }

\begin{corollary}\label{cor:Duality-3}
Assume that $u$ is a weak solution of the
\changed{NSEs}{Navier-Stokes equations} \eqref{eqn_NSES} with the
external force  $0$. Assume that $u$ satisfies the Serrin
condition, i..e  $u\in \bigcup_{\delta_0>0}
L^{2+\frac{6}{\delta_0}}(0,T;\mathbb{L}^{3+\delta_0}(\Rnu^3))$ and
$u(0)\in V$. Then   $u\in L^{\infty}(0,T;V)$, i.e. $u$ is a strong
solution of \eqref{eqn_NSES}. Moreover, if $G_0\in H$,  then
\begin{eqnarray}
(\curl u(0), \mathcal{T}_T^{S_{T}(u)}G_0)_H=(\curl
u(T),G_0)_H,\label{eqn:Duality-3}\\
\Vert \curl u(T)\Vert _H\leq \Vert \mathcal{T}_T^{S_{T}(u)}\Vert
_{\mathcal{L}(H,H)}\Vert \curl u(0)\Vert.
_H\label{eqn:CurlEstimate-1}
\end{eqnarray}

\end{corollary}
\begin{remark}\label{Rem:gen_Kelvin_0}
Let us observe  that the equality \eqref{eqn:Duality-3} is a
generalization of the helicity invariance \del{What is $h$ in the
integral?Changed. Just Misprint.} $\intl_{\Rnu^3}(u,\curl
u)_{\Rnu^3}\,dx$, see e.g. p. 120 -- 121 in  \cite{Moffat} for the
solutions of the Euler equations\del{ I made some changes here}.
Indeed, if we consider the transport operator
$\mathcal{T}_T^{\cdot}$ for $\nu=0$ and take $G_0=u(T)$ on  the
right hand side of equality \eqref{eqn:Duality-3} then, under the
assumption that the Euler equation has a unique solution, we infer
that $\mathcal{T}_T^{S_{T}(u)}u(T)=u(0)$.
\end{remark}
\begin{proof}[Proof of Corollary \ref{cor:Duality-3}]
By Proposition \ref{prop:ExistUniqSolution} there exist unique
solution $F\in L^2(0,T;V)\cap L^{\infty}(0,T;H)$ of equation
(\ref{eqn:DualEquation}-\ref{in_v}) with initial condition
$F_0=u(0)$ and $v=u$. We can notice that $u$ is also solution of
(\ref{eqn:DualEquation}-\ref{in_v}) by Navier-Stokes equation.
Thus, $F=u$ and we have \eqref{eqn:Duality-3} by Theorem
\ref{thm:Duality}. Therefore, we have $$\Vert \curl u(t)\Vert
_H\leq \Vert \mathcal{T}_T^{S_{T}(u)}\Vert
_{\mathcal{L}(H,H)}\Vert \curl u(0)\Vert _H$$ and by boundedness
of operator $\mathcal{T}_T^{S_{T}(u)}$ (Proposition
\ref{prop:ExistUniqSolution}) we get the result.
\end{proof}
\del{Hence, in order to get an estimate on $\Vert \curl u(t)\Vert
_{\mathbb{L}^2}$ (where $u$-is a strong solution of Navier-Stokes
equation) from above we need to study operator of vector transport
$\mathcal{T}_T^{u}$ for very irregular $u$. In next section we
will get a Feynman-Kac type formula for $\mathcal{T}_T^{u}$.}
%From theorem \eqref{prop:ExistUniqSolution} we immediately get
%classical result
%\begin{corollary}
%If $u$ satisfies Serrin condition \eqref{eqn:IntegrCond-1}
%\end{corollary}
\del{END. I finished here for now. Not much, but some progress.}

\section{Formulae of Feynman-Kac
Type.}\label{sec:FeynmanKacFormulae}

 The aim of this section is twofold. Firstly,  we will discuss the
physical meaning  of the operator $\mathcal{T}_T^{S_{T}(\cdot)}$.
Secondly,  \del{and in the same time} we will deduce a formula of
Feynman-Kac type. \del{From now on} In the whole section  we
suppose
that  $D=\Rn$. %and $v\in C_0^{\infty}([0,T]\times\Rn)$.\Red{Why so
%very regular?}
We also assume that $(\Omega,\mathcal{F},\{\mathcal{F}_t\}_{t\geq
0},\mathbb{P})$ is a complete filtered probability space and that
$(W(t))_{t\geq 0}$ is an $\mathbb{R}^m$-valued  Wiener process on
this space. We have
the following Proposition.%, see \cite{Nekl}.

\begin{proposition}\label{prop:ContourConservation}
Assume that $\alpha \in (0,1)$, $\sigma(\cdot,\cdot)\in
L^1(0,T;C_{b}^{2,\al}(\Rn ,\Rn\otimes\Rnu^{m}))$,
$a(\cdot,\cdot)\in L^1(0,T;C_{b}^{1,\al}(\Rn ,\Rn))$. Let us
assume that  a continuous and adapted process
$X=[0,T]\times\Rn\times\Omega\to\Rn$ is a unique %\Blue{Do we really use uniqueness here? In the proof you mention %that $X$ generates a flow of diffeomorphisms, but in the formula \eqref{eqn:ContourItoFormula} below, this flow %doesn't appear. If I am right, we can substantially weaken the assumptions. Probably, I am wrong. We need flow to %define $X_s(\Gamma)$. Maybe make a remark of this sort?}
solution to the problem
$$
\begin{array}{l}
dX_t(x) = a(t,X_t(x))\,dt+\sigma(t,X_t(x))\,dW(t),\\
X_0(x) = x.
\end{array}
$$
Then for any     $C^1$ class  closed loop $\Gamma$ in $\Rn$,  any
$F\in C^{1,2}([0,T]\times\Rn,\Rn)$ and any $t\in [0,T]$, we have
$\mathbb{P}$-a.s.,
\begin{eqnarray}\label{eqn:ContourItoFormula}
&&\intl_{X_t(\Gamma)}\suml_{k=1}^{n}F^k(t,x)\,dx_k=\intl_{\Gamma}\suml_{k=1}^{n}F^k(0,x)\,dx_k\\
&+&
\intl_0^t\intl_{X_s(\Gamma)}\suml_{k=1}^{n}\left(\frac{\partial
F^k}{\partial t}+\suml_{j=1}^{n}a^{j}(\frac{\partial F^k}{\partial
x_j}-\frac{\partial F^j}{\partial
x_k})+\frac{1}{2}\suml_{i,j=1}^n\frac{\partial^2 F^k}{\partial
x_i\partial
x_j}\suml_{m=1}^{n}\sigma^{im}\sigma^{jm}\right)\,dx_kds\nonumber\\
&&+\frac{1}{2}\intl_0^t\intl_{X_s(\Gamma)}\suml_{k=1}^{n}
\left(\suml_{j,l}\frac{\partial F^j}{\partial
x_l}\suml_m\sigma^{lm}\frac{\partial \sigma^{jm}}{\partial
x_k}\right)\,dx_kds\label{eqn:correction_term}\\
&+&\intl_0^t\intl_{X_s(\Gamma)}
\suml_{k,j=1}^{n}F^j(s,x)\frac{\partial \sigma^{jl}}{\partial
x_k}\,dx_kdw_s^l+
\intl_0^t\intl_{X_s(\Gamma)}\suml_{k=1}^n\left(\suml_{i,l=1}\frac{\partial
F^k}{\partial x_i}\sigma^{il}\right)\,dx_k\,dW^{l}(s).\nonumber
\end{eqnarray}
% \Red{Can you think of a compact version of the above formula. For
% example (this is easy) instead of
% $\intl_{X_t(\Gamma)}\suml_{k=1}^{n}F^k(t,x)\,dx_k$ we could write
% $\intl_{X_t(\Gamma)}F(t,x)\cdot dx$.}
\end{proposition}
%\Green{
\begin{remark}
The term \eqref{eqn:correction_term} is of major interest for us.
Its appearance allows us  to " emulate" drift in two dimensional
case i.e. to consider flow without drift such that this term
"creates" necessary drift (see subsections \ref{example-2D},
\ref{example-3D} and Theorem \ref{thm:2DFlowDifferentLaw}for
detailed explanation).
\end{remark}
%}
\begin{proof}[Proof of Proposition \ref{prop:ContourConservation}]
It follows from Theorems 3.3.3, p.94 and  4.6.5, p.173 of
\cite{Kunita} that $X_t(\cdot),t\in [0,T]$ is a flow of
$C^1$--diffeomorphisms and $\nabla X_t(\cdot)$ satisfies
corresponding equation for gradient of the flow. Then formula
\eqref{eqn:ContourItoFormula} immediately follows from the It\^o
formula, see \cite{Nekl} for calculations.
%\Blue{We need to write what is the different between our result here and your result in \cite{Nekl}.}
\end{proof}
\begin{corollary}\label{cor:LocMartingaleProperty}
 Assume that  $\nu>0$ and  $v\in L^1(0,T;C_{b}^{1,\al}(\Rn ,\Rn))$ for some
$\al\in (0,1)$. Let   $(\changed{X_s(t;x)}{X_t^s})_{0\leq s\leq
t\leq T}$, be a stochastic flow corresponding to the following SDE
\begin{eqnarray}
d\changed{X_s(t;x)}{X_t^s} &=&
v(t,\changed{X_s(t;x)}{X_t^s})\,dt+\sqrt{2\nu}\,dW(t),\; t\in
[s,T], \label{eqn:flow-1}
\\X_s(s;x) &=& x.\nonumber
\end{eqnarray}
Assume that  $F_0\in C^2(\Rn)$ and  let $F\in
C^{1,2}([0,T]\times\Rn,\Rn)$ be a solution of the following linear
equation \footnote{which coincides with Problem
\eqref{eqn:DualEq-1} in the case $n=3$}
\begin{eqnarray}
\frac{\partial F(t)}{\partial t} &=& -\nu AF+\mathrm{P}((v(T-t)\nabla)F-\nabla F v(T-t)),t\in (0,T),\label{eqn:eq-3}\\
F(0) &=& F_0, \label{eqn:eq-3b}
\end{eqnarray}
 Then for any $s \in [0,T]$ a
process $(M_s(t))_{t\in [T-s,T]}$ defined by the following formula
$$M_s(t)=\intl_{X_{T-s}(t;\Gamma)}\suml_{k=1}^nF^k(T-t)\,dx_k,\; t\in[T-s,T]$$
is a local martingale.
\end{corollary}
\begin{proof}[Proof of Corollary \ref{cor:LocMartingaleProperty}]
This follows immediately from Proposition
\ref{prop:ContourConservation}.
\end{proof}
\begin{remark}\label{Rem:gen_Kelvin}
Corollary \ref{cor:LocMartingaleProperty}, whose  idea is taken
from \cite{Nekl}, can be seen as  a generalization of the Kelvin
circulation Theorem,  see e.g. \cite[ p. 26]{MarchioroPulvirenti}.
Indeed, if $\nu=0$, then $\changed{X_s(t;x)}{X_t^s}$ is a position
of a particle at time $t$ starting from point $x$ at time $s$,
moving in the deterministic velocity field $v$. Moreover, $M_s$ is
the circulation along a curve $\gamma$  moved by the flow
generated by $v$. Hence, by Proposition
\ref{prop:ContourConservation} it follows  that the local
martingale $M_s$ is constant in time. A similar  result has
recently been independently derived  by Constantin and Iyer, see
\cite[Proposition 2.9]{Const_Iyer_2008}.
\end{remark}

Next  we deduce from the corollary \ref{cor:LocMartingaleProperty}
the following formula of the Feynman-Kac type for the solution of
equation \eqref{eqn:eq-3}.

\begin{proposition}\label{prop:FeynKacformula}
Assume that $v\in L^1(0,T;C_{b}^{2,\al}(\Rn ,\Rn))$ for some
$\al\in (0,1)$ and
\begin{equation}
\intl_0^T(|v|_{L^{\infty}}(s)+|\nabla
v|_{L^{\infty}}(s))\,ds<\infty.\label{eqn:UniIntegrCond-2}
\end{equation}
Assume that   $F:[0,T]\times\Rn\to\Rn$ is a solution of the
problem \eqref{eqn:eq-3}-\eqref{eqn:eq-3b} with $F_0\in
C^{2}(\Rnu^n)\cap L^4(\Rnu^n)$ and
$(\changed{X_s(t;x)}{X_t^s})_{0\leq s\leq t\leq T}$ is a
stochastic flow corresponding to SDE \eqref{eqn:flow-1}. Assume
also that there exists $\beta>0$ such that for  any $\Gamma\in
C^1(\mathbb{S}^1,\Rn)$, where $\mathbb{S}^1$ is the unit circle,
for all $s,t\in[0,T]$ such that  $ T-s\leq t$,
\begin{equation}
\mathbb{E}|\intl_{X_{T-s}(t;\Gamma)}F^k(T-t,x)\,dx_k|^{1+\beta}<\infty,\;
k=1,\cdots,n.\label{eqn:UniIntegrCond}
\end{equation}
Fix $s\in [0,T]$ and define a functions $Q_s:\mathbb{R}^n\to
\mathbb{R}^n$ by
$$Q_s(x):=\mathbb{E}(F_0(\changed{X_{T-s}(T;x)}{X_T^{T-s}(x)})\nabla \changed{X_{T-s}(T;x)}{X_T^{T-s}(x)})),\; \; x\in\Rn.$$

Then, $Q_s\in L^2(\mathbb{R}^n,\mathbb{R}^n)\cap
C^{1+\eps}(\mathbb{R}^n,\mathbb{R}^n)$,  $0<\eps<\alpha$ and
\begin{equation}\label{eqn:FeynKacFormula}
F(s,x)=[\mathrm{P}(Q_s)](x),\; x\in \mathbb{R}^n,\; s\in [0,T].\;
\end{equation}
% \Red{Changed. The above formula (and other related) have to be rewritten.
% Firstly we fix $s\in [0,T]$ and define a functions
% $Q_s:\mathbb{R}^n\to \mathbb{R}^n$ by
% $$Q_s(x):=\mathbb{E}(F_0(\changed{X_{T-s}(T;x)}{X_T^{T-s}(x)})\nabla \changed{X_{T-s}(T;x)}{X_T^{T-s}(x)}),\; \; x\in\Rn.$$
% Secondly, we prove that $Q_s$ belongs to such a space that $\mathrm{P}(Q_s)$
% makes sense. Finally, we write $$ F(s,x)=[\mathrm{P}(Q_s)](x),\; x\in
% \mathbb{R}^n.$$ We need also to check that $\mathrm{P}(Q_s)\in C(\Rn,\Rn)$ so
% that the above formula holds true for all $x\in\Rn$. }
\end{proposition}
%\begin{proof}\Red{what about the proof? Is it true that the mean
%of a local martingale is constant in time?}
%\end{proof}
\begin{remark}
In connection with the formula \eqref{eqn:FeynKacFormula} we can
ask whether the flow $(\changed{X_s(t;x)}{X_t^s})_{0\leq s\leq
t\leq T}$ associated to the SDE \eqref{eqn:flow-1} is the only
flow such that the function $F$ defined by  the formula
\eqref{eqn:FeynKacFormula} is  a solution to problem
(\ref{eqn:eq-3}-\ref{eqn:eq-3b})? It turns out that the answer to
this question is negative. In the subsections \ref{example-2D} and
\ref{example-3D} we will consider separately two and three
dimensional examples.
\end{remark}
\begin{remark}
Condition \eqref{eqn:UniIntegrCond} is satisfied if, for instance,
$F\in L^{\infty}([0,T]\times\Rn)$ and
$$
\intl_0^T|\nabla v|_{L^{\infty}}(s)\,ds<\infty.
$$
Indeed, in this case we have  the following inequality
$$
|\nabla \changed{X_s(t;\cdot)}{X_t^s}|_{L^{\infty}}\leq
exp(\intl_s^t|\nabla v|_{L^{\infty}}(r)\,dr),\;s\leq t\leq T,
$$
and hence the result follows.
\end{remark}
\begin{proof}[Proof of Proposition \ref{prop:FeynKacformula}]
For fixed $s\in [0,T)$ let us denote
\begin{equation}
M_s(t)=\intl_{X_{T-s}(t;\Gamma)}\suml_{k=1}^3F^k(T-t)\,dx_k,\,
t\in[T-s,T].\label{eqn:LocMartingale-1}
\end{equation}
Then by Corollary \ref{cor:LocMartingaleProperty} the process
$(M_s(t))$, $t\in[T-s,T]$ is a local martingale. Hence, by the
uniform integrability condition \eqref{eqn:UniIntegrCond} we infer
that $M_s$ is martingale and so
$\mathbb{E}M_s(T-s)=\mathbb{E}M_s(T)$. In particular,
\begin{equation}
\intl_{\Gamma}F^k(s,x)\,dx_k=\intl_{\Gamma}Q_s^k(x)\,dx_k,\Gamma\in
C^1(S^1,\Rn).
\end{equation}
It immediately follows from Theorems 3.3.3, p.94 and  4.6.5, p.173
of \cite{Kunita} that $Q_s\in C^{1+\eps}(\Rn,\Rn)$,
$0<\eps<\alpha$. Furthermore, $Q_s\in L^2(\Rn,\Rn)$. Indeed, by
definition of the flow \eqref{eqn:flow-1} we have
\begin{equation}
\supl_x|\nabla X_{T-s}(T;x)|\leq e^{\intl_0^T|\nabla
v|_{L^{\infty}}(r)\,dr}.\nonumber
\end{equation}
Hence
\begin{eqnarray}
\intl_{\Rn}|Q_s(x)|^2\,dx &\leq& \intl_{\Rn}\mathbb{E}|F_0(\changed{X_{T-s}(T;x)}{X_T^{T-s}(x)})\nabla \changed{X_{T-s}(T;x)}{X_T^{T-s}(x)}|^2\,dx\nonumber\\
&\leq& \mathbb{E}\Big(\supl_x|\nabla
\changed{X_{T-s}(T;x)}{X_T^{T-s}(x)}|^2\intl_{\Rn}|F_0(\changed{X_{T-s}(T;x)}{X_T^{T-s}(x)})|^2\,dx\Big)\nonumber\\
&\leq&
e^{\intl_0^T|\nabla v|_{L^{\infty}}(r)\,dr}\mathbb{E}\intl_{\Rn}|F_0(\changed{X_{T-s}(T;x)}{X_T^{T-s}(x)})|^2\,dx\label{eqn:auxEst-1}\\
&\leq& e^{\intl_0^T|\nabla
v|_{L^{\infty}}(r)\,dr}\intl_{\Rn}\mathbb{\tilde{E}}(|F_0(x+\sqrt{2\nu}(W_T-W_{T-s}))|^2\mathcal{E}_{T-s}^T)\,dx,\nonumber
\end{eqnarray}
where
$\mathcal{E}_{T-s}^T=e^{\intl_{T-s}^Tv(r,X_{T-s}(r;x))\,dW_r-1/2\intl_{T-s}^T|v(r,X_{T-s}(r;x))|^2\,dr}$
is a stochastic exponent. We can notice that
\begin{equation}
\mathbb{\tilde{E}}|\mathcal{E}_{T-s}^T|^2\leq
e^{2\intl_0^T|v(r)|_{L^{\infty}}(r)\,dr}\label{eqn:auxEst-2}
\end{equation}
and, therefore, combining \eqref{eqn:auxEst-1} and
\eqref{eqn:auxEst-2} we get
\begin{equation}
\intl_{\Rn}|Q_s(x)|^2\,dx\leq
e^{\intl_0^T(|v|_{L^{\infty}}(r)+|\nabla
v|_{L^{\infty}}(r))\,dr}\intl_{\Rn} |F_0|^4dx<\infty.
\end{equation}
It remains to notice that operator
$\mathrm{P}:C^{\beta}(\Rn,\Rn)\to C^{\beta}(\Rn,\Rn),\beta\in
(0,1)$ is bounded. Indeed, it follows from representation of
$\mathrm{P}$ as pseudodifferential operator
(\cite{Fuj_Mor_1977},\cite{Temam_2001}) and Theorem 7.9.6 in
\cite{Hormander-1990}.
\end{proof}
\begin{remark}
Another method of proving the formula \eqref{eqn:FeynKacFormula}
is presented in the article  \cite{BusnelloFlandoliRomito} by
Busnello et al., see also literature therein. The approach used
there is based upon an extension of the standard Feynman-Kac
formula for parabolic equations to more general system of linear
parabolic equations with a potential term (see the system (3.2)
in \cite[p.306]{BusnelloFlandoliRomito}). This extension is
carried out by using  the new variables method   introduced
earlier by Krylov \cite{Krylov}. One should mention here that the
formula \eqref{eqn:FeynKacFormula} is used in
\cite{BusnelloFlandoliRomito} to prove the local existence and
uniqueness of strong solutions to the \changed{NSEs}{Navier-Stokes
equations}.
\end{remark}

\subsection{Examples of nontrivial flows in $\mathbb{R}^2$}\label{example-2D}

In this subsection we provide nontrivial examples of the flows
which can be used in the Feynman-Kac type formula
\eqref{eqn:FeynKacFormula} in two dimensional case.

\begin{proposition}\label{prop:2DFlowTheSameLaw}
Suppose that $v\in C_0^{\infty}([0,T]\times\Rnu^2,\Rnu^2)$, $\psi:
\Rnu \to \Rnu$ is a $C^1$-class diffeomorphism, $\phi=\psi\circ
\rot v$ and  $F_0\in C_0^{\infty}(\Rn)$. Let
$(\changed{X_s(t;x)}{X_t^s})$, $0\leq s\leq t\leq T$ be the
stochastic flow corresponding to the following SDE
\begin{eqnarray}
d\changed{X_s(t;x)}{X_t^s} &=& v(t,\changed{X_s(t;x)}{X_t^s})\,dt+\sqrt{2\nu}\sigma_1(\changed{X_s(t;x)}{X_t^s})\,dW(t), \label{eqn:flow-2}\\
X_s(s;x) &=& x,\nonumber
\end{eqnarray}
where
\begin{displaymath}
\sigma_1 (x) = \left(\begin{array}{cc} \cos\phi(x) & -\sin\phi(x)\\
\sin\phi(x) & \cos\phi(x)
\end{array}\right), x\in\Rnu^2.
\end{displaymath}
Assume that $F:[0,T]\times\Rn\to\Rn$ is a solution to problem
(\ref{eqn:eq-3}-\ref{eqn:eq-3}) such that for some $\beta>0$ and
any
 $\Gamma\in C^1(\mathbb{S}^1,\mathbb{R}^2)$ the condition \eqref{eqn:UniIntegrCond} is
satisfied. Then, the  formula \eqref{eqn:FeynKacFormula} holds
true.
\end{proposition}
\begin{proof}[Proof of Proposition \ref{prop:2DFlowTheSameLaw}]
Suppose that the condition \eqref{eqn:UniIntegrCond} is fulfilled.
Then, it is enough to show that process $(M_s(t)),t\in[T-s,T]$
defined by formula \eqref{eqn:LocMartingale-1} above (where flow
$(\changed{X_s(t;x)}{X_t^s}),0\leq s\leq t\leq T$ is given by
\eqref{eqn:flow-2}) is a local martingale. We have
\begin{eqnarray}
\intl_{X_{T-s}(t;\Gamma)}\suml_{k=1}^{n}F^k(T-t,x)\,dx_k=\intl_{\Gamma}\suml_{k=1}^{n}F^k(s,x)\,dx_k\nonumber
\end{eqnarray}
\begin{eqnarray}
+\intl_{T-s}^t\intl_{X_{T-s}(\tau;\Gamma)}\suml_{k=1}^{n}\left(\frac{\partial
F^k}{\partial t}+\suml_{j=1}^{n}v^{j}(\frac{\partial F^k}{\partial
x_j}-\frac{\partial F^j}{\partial
x_k})+\nu\suml_{i,j=1}^n\frac{\partial^2 F^k}{\partial x_i\partial
x_j}\suml_{m=1}^{n}\sigma_1^{im}\sigma_1^{jm}\right)\,dx_kd\tau
+\nonumber
\end{eqnarray}
\begin{eqnarray}
&+&\nu\intl_{T-s}^t\intl_{X_{T-s}(\tau;\Gamma)}\suml_{k=1}^{n}
\left(\suml_{j,l}\frac{\partial F^j}{\partial
x_l}\suml_m\sigma_1^{lm}\frac{\partial \sigma_1^{jm}}{\partial
x_k}\right)\,dx_kd\tau\nonumber\\
&+&\sqrt{2\nu}\intl_{T-s}^t\intl_{X_{T-s}(\tau;\Gamma)}\suml_{k,j=1}^{n}F^j(T-\tau,x)\frac{\partial
\sigma_1^{jl}}{\partial x_k}\,dx_kdw_{\tau}^l\nonumber\\
&+&\sqrt{2\nu}
\intl_{T-s}^t\intl_{X_{T-s}(\tau;\Gamma)}\suml_{k=1}^n\left(\suml_{i,l=1}\frac{\partial
F^k}{\partial
x_i}\sigma_1^{il}\right)\,dx_k\,dW_{\tau}^{l}.\nonumber
\end{eqnarray}
Hence, because $\sigma_1$ is orthogonal matrix and $F$ satisfies
\eqref{eqn:eq-3} we have that
\begin{eqnarray*}
\frac{\partial F^k}{\partial
t}&+&\suml_{j=1}^{n}v^{j}(\frac{\partial F^k}{\partial
x_j}-\frac{\partial F^j}{\partial
x_k})+\nu\suml_{i,j=1}^n\frac{\partial^2 F^k}{\partial
x_i\partial x_j}\suml_{m=1}^{n}\sigma_1^{im}\sigma_1^{jm}\\
&=&\frac{\partial F^k}{\partial
t}+\suml_{j=1}^{n}v^{j}(\frac{\partial F^k}{\partial
x_j}-\frac{\partial F^j}{\partial x_k})+\nu\triangle
F^k=\frac{\partial p}{\partial x_k}.
\end{eqnarray*}
Therefore, it is enough to show that
\begin{eqnarray*}
\intl_{T-s}^t\intl_{X_{T-s}(\tau;\Gamma)}\suml_{k=1}^{n}
\left(\suml_{j,l}\frac{\partial F^j}{\partial
x_l}\suml_m\sigma_1^{lm}\frac{\partial \sigma_1^{jm}}{\partial
x_k}\right)\,dx_kd\tau =0.
\end{eqnarray*}
We have that $\suml_m\sigma_1^{lm}\frac{\partial
\sigma_1^{jm}}{\partial x_k}$ is antisymmetric w.r.t. indexes
$l,j$ because $\sigma_1$ is orthogonal. Hence $n=2$, it means that
it is enough to calculate
\begin{equation*}
\suml_m\sigma_1^{1m}\frac{\partial \sigma_1^{2m}}{\partial x_k}=
\cos\phi\frac{\partial}{\partial
x_k}(\sin\phi)-\sin\phi\frac{\partial}{\partial
x_k}(\cos\phi)=\frac{\partial\phi}{\partial x_k}
\end{equation*}
and, therefore,
\begin{eqnarray*}
\intl_{T-s}^t&&\intl_{X_{T-s}(\tau;\Gamma)}\suml_{k=1}^{n}
\left(\suml_{j,l}\frac{\partial F^j}{\partial
x_l}\suml_m\sigma_1^{lm}\frac{\partial \sigma_1^{jm}}{\partial
x_k}\right)\,dx_kd\tau\\
&=&\intl_{T-s}^t\intl_{X_{T-s}(\tau;\Gamma)}(\frac{\partial
F^1}{\partial x_2}-\frac{\partial F^2}{\partial x_1})\,d\phi d\tau
=\intl_{T-s}^t\intl_{X_{T-s}(\tau;\Gamma)}\psi^{-1}(\phi)\,d\phi
d\tau =0.
\end{eqnarray*}
\end{proof}
\begin{remark}
The construction of the example from  Proposition
\ref{prop:2DFlowTheSameLaw} can easily be generalized to the case
$n=3$ in the following way. Let  $\psi: \Rnu \to \Rnu$ be a
$C^1$-class diffeomorphism. Define $\phi=\psi\circ (\curl v)^1$
and
\begin{displaymath}
\sigma_1 (x) = \left(\begin{array}{ccc} \cos\phi(x) & -\sin\phi(x) & 0\\
\sin\phi(x) & \cos\phi(x) & 0\\
0 & 0 &1
\end{array}\right), x\in\Rnu^3.
\end{displaymath}
Let $(\changed{X_s(t;x)}{X_t^s})_{0\leq s\leq t\leq T}$ be a
stochastic flow corresponding to the following SDE
\begin{eqnarray}
d\changed{X_s(t;x)}{X_t^s} &=&
v(t,\changed{X_s(t;x)}{X_t^s})\,dt+\sqrt{2\nu}\sigma_1(\changed{X_s(t;x)}{X_t^s})\,dW(t),
0\leq s\leq t\leq T\label{eqn:flow-3}\\
X_s(s;x) &=& x\nonumber
\end{eqnarray}
Then the assertion of Proposition \eqref{prop:2DFlowTheSameLaw}
holds true.

Note that similar construction can be made for other components of
the $\curl v$)  but the truly three dimensional rotations
$\sigma_1$ will be considered in next paragraph.
\end{remark}

\begin{remark}
Let us note that the laws of the solutions to SDEs
\eqref{eqn:flow-2} and \del{the law of canonical flow}
\eqref{eqn:flow-1} are the same. Indeed, it is easy to see that
quadratic variations of both processes are the same.
 In the next example we
will show that it is possible to find a flow such that its
one-point motion has a law of Brownian motion.
\end{remark}
\begin{theorem}\label{thm:2DFlowDifferentLaw}
Suppose that $\nu>0$, $\delta>0$ and a divergence free vector
field $v:\mathbb{R}^2\to \mathbb{R}^2$ is of $C^{1+\delta}$ class.
Let $\phi:\mathbb{R}^2\to \mathbb{R}$ be such that\footnote{Such
$\phi$ exists because $\diver v=0$. } $v=\nabla^{\bot}\phi$.
Define
\begin{displaymath}
\sigma_1 (x) = \left(\begin{array}{cc} \cos\frac{\phi(x)}{\nu} & -\sin\frac{\phi(x)}{\nu}\\
\sin\frac{\phi(x)}{\nu} & \cos\frac{\phi(x)}{\nu}
\end{array}\right), x\in\Rnu^2,
\end{displaymath}
Let us denote by $\changed{X_s(t;x)}{\changed{X_s(t;x)}{X_t^s}}$,
$0\leq s\leq t\leq T, x\in\Rnu^2 $
 the stochastic flow of diffeomorphisms of $\Rnu^2$ of class $C^2$ corresponding to the following SDE
\begin{equation}
\label{eqn:flow-4}
\left\{
\begin{array}{rcl}
d\changed{X_s(t;x)}{\changed{X_s(t;x)}{X_t^s}} &=&
\sqrt{2\nu}\sigma_1(\changed{X_s(t;x)}{X_t^s})\,dW(t), \;0\leq
s\leq t\leq
T,\\
X_s(s;x) &=& x.
\end{array}
\right.
\end{equation}

  Assume also that $F_0\in
C^{2}(\Rnu^2)\cap L^2(\Rnu^2)$ and that
$F:[0,T]\times\Rnu^2\to\Rnu^2$ is a solution to problem
(\ref{eqn:eq-3}-\ref{eqn:eq-3b}) such that for some $\beta>0$ and
any  $\Gamma\in C^1(\mathbb{S}^1,\mathbb{R}^2)$ the condition
\eqref{eqn:UniIntegrCond} is satisfied. Denote
$Q_s(x)=\mathbb{E}(F_0(\changed{X_{T-s}(T;x)}{X_T^{T-s}(x)})\nabla
\changed{X_{T-s}(T;x)}{X_T^{T-s}(x)})$. Then $Q_s\in L^2(\Rn)\cap
C^{1+\eps}(\Rn)$,  $0<\eps<\delta$ and
\begin{equation}\label{eqn:FeynKacFormula-2}
F(s,x)=\mathrm{P}(Q_s)(x),\; s\in [0,T],\; x\in\mathbb{R}^n.
\end{equation}
\end{theorem}
%\Green{
\begin{remark} As we have already noticed above the formula
\eqref{eqn:FeynKacFormula-2} can be viewed as generalization of
Kelvin Theorem as in the Corollary
\ref{cor:LocMartingaleProperty}. Indeed, it is enough to integrate
both sides of \eqref{eqn:FeynKacFormula-2} w.r.t. arbitrary smooth
closed contour $\Gamma$.
\end{remark}%}
\begin{proof}[Proof of Theorem \ref{thm:2DFlowDifferentLaw}]
From Theorem 4.6.5 , p. 173 in \cite{Kunita} we infer that there exists a  flow $X_s(t;x),0\leq s\leq t\leq T$ for problem \eqref{eqn:flow-4} consisting of diffeomorphisms of class $C^{2+\eps}$ \del{\textbf{Furthermore, it is the process with values in $C^{2+\eps}(\Rn,\Rn)$} for any $0<\eps<\delta$.} %\Blue{I have changed the order of words here.}

Moreover, it  follows from Theorems 3.3.3, p. 94 and  4.6.5, p.
173 therein \del{of \cite{Kunita}} that for all $s\in [0,T]$,
$Q_s\in C^{1+\eps}(\Rn,\Rn)$, $0<\eps<\delta$. Let us fix $s\in
[0,T]$. We will show now that $Q_s\in L^2(\Rn,\Rn)$. Since by
Corollary 4.6.7 p. 175 of \cite{Kunita} that there exists a
positive constant $C$ such that
$$
\sup_{x\in\Rn}\mathbb{E}|\nabla
\changed{X_{T-s}(T;x)}{X_T^{T-s}(x)}|^2\leq C,
$$
 by the H\"older inequality we infer that
\begin{eqnarray}
\intl_{\Rn}|Q_s(x)|^2\,dx\leq
\intl_{\Rn}\mathbb{E}|F_0(\changed{X_{T-s}(T;x)}{X_T^{T-s}(x)})|^2\mathbb{E}|\nabla
\changed{X_{T-s}(T;x)}{X_T^{T-s}(x)}|^2\,dx\nonumber\\
\leq
C\intl_{\Rn}\mathbb{E}|F_0(\changed{X_{T-s}(T;x)}{X_T^{T-s}(x)})|^2\,dx.\label{eqn:auxIneq-1}
\end{eqnarray}

Furthermore, let us observe that the law of the one-point motion
of the flow $\changed{X_{T-s}(T;x)}{X_T^{T-s}(x)}$ is equal to the
law of the Brownian Motion (see example 6.1, p. 75 of
\cite{[Ikeda-Watanabe-1981]} for more details). Therefore, we
infer that
\begin{equation}
\intl_{\Rn}\mathbb{E}|F_0(\changed{X_{T-s}(T;x)}{X_T^{T-s}(x)})|^2\,dx=\intl_{\Rn}|S_s^{\nu}F_0(x)|^2\,dx\leq
\intl_{\Rn}|F_0(x)|^2\,dx,\label{eqn:auxIneq-12}
\end{equation}
where $\{S_t^{\nu}=e^{\nu t\triangle}\}_{t\geq 0}$ is a heat
semigroup. Combining inequalities \eqref{eqn:auxIneq-1} and
\eqref{eqn:auxIneq-12} we get
\begin{equation}
\intl_{\Rn}|Q_s(x)|^2\,dx\leq C\intl_{\Rn}|F_0(x)|^2\,dx.
\end{equation}
Similarly to Proposition \ref{prop:2DFlowTheSameLaw} we get that
$\intl_{X_{T-s}(t;\Gamma)}\suml_{k=1}^3F^k(T-t)\,dx_k,t\in[T-s,T]$
is a local martingale. Indeed, correction term in
\eqref{eqn:ContourItoFormula} due to rotation of Brownian Motion
is equal to
$\intl_{T-s}^t\intl_{X_{T-s}(\tau;\Gamma)}(\frac{\partial
F^1}{\partial x_2}-\frac{\partial F^2}{\partial x_1})\,d\phi\,
ds$, see the previous Proposition, and if $v=\nabla^{\bot}\phi$
this is exactly first order term of two dimensional equation
\eqref{eqn:eq-3}.
\end{proof}

%\del{09 August, 2007. We still need to prove that \eqref{eqn:FeynKacFormula-2} holds for all $x$. For this we could %need to know that the operator $\mathrm{P}$ is continuous in some class of Hoelder continuous functions. I will look into %this issue.} \Blue{The last issue has to be resolved but we can submit now. }
%\begin{remark}
%Under assumptions of Theorem \ref{thm:2DFlowDifferentLaw}
%formula \eqref{eqn:FeynKacFormula} can be simplified as follows
%
%where $X_t^{s}(x))$ is defined by \eqref{eqn:flow-4}.
%\end{remark}

\begin{corollary}\label{cor:2DGradientFormula-1}
Let  $\changed{(X_s(t;x))}{X_t^s(x)}$ $0\leq s\leq t\leq T,
x\in\Rnu^2 $ be  the stochastic flow corresponding to  SDE
\eqref{eqn:flow-4}. Then \del{For the flow $X_t^s(\cdot)$ defined
by \eqref{eqn:flow-4} we have that}\del{The notation below is a
bit confusing. What about using instead of yours the following
one: $(X_s^1(t,x)$? I have made these changes everywhere. If find
some places, where I missed, please do them.}
\begin{eqnarray}
d\left(
\begin{array}{cc}
\frac{\partial \changed{X_s^1(t;x)}{X_t^{1s}}}{\partial x_1} & \frac{\partial \changed{X_s^1(t;x)}{X_t^{1s}}}{\partial x_2}\\
\frac{\partial \changed{X_s^2(t;x)}{X_t^{2s}}}{\partial x_1} &
\frac{\partial \changed{X_s^2(t;x)}{X_t^{2s}}}{\partial x_2}
\end{array}
\right)&=&\frac{1}{\nu} \left(
\begin{array}{cc}
-v_2(t,\changed{X_s(t;x)}{X_t^s})\,d\changed{X_s^2(t;x)}{X_t^{2s}} & v_1(t,\changed{X_s(t;x)}{X_t^s})\,d\changed{X_s^2(t;x)}{X_t^{2s}}\\
v_2(t,\changed{X_s(t;x)}{X_t^s})\,d\changed{X_s^1(t;x)}{X_t^{1s}}
&
-v_1(t,\changed{X_s(t;x)}{X_t^s})\,d\changed{X_s^1(t;x)}{X_t^{1s}}
\end{array}
\right)\nonumber\\&& \left(
\begin{array}{cc}
\frac{\partial \changed{X_s^1(t;x)}{X_t^{1s}}}{\partial x_1} & \frac{\partial \changed{X_s^1(t;x)}{X_t^{1s}}}{\partial x_2}\\
\frac{\partial \changed{X_s^2(t;x)}{X_t^{2s}}}{\partial x_1} &
\frac{\partial \changed{X_s^2(t;x)}{X_t^{2s}}}{\partial x_2}
\end{array}
\right),\label{eqn:2DGradientFormula}
\end{eqnarray}
and
\begin{displaymath}
\left(
\begin{array}{cc}
\frac{\partial X_s^1(s;x)}{\partial x_1} & \frac{\partial X_s^1(s;x)}{\partial x_2}\\
\frac{\partial X_s^2(s;x)}{\partial x_1} & \frac{\partial
X_s^2(s;x)}{\partial x_2}
\end{array}
\right)= \left(
\begin{array}{cc}
1 & 0\\
0 & 1
\end{array}
\right).
\end{displaymath}
%Moreover,
\del{What do you mean by $\exp$ in the formula below? Have you
taken into account that the coefficients in the linear equation
above are time dependent?Decided to skip this formula for a
while.}
%\begin{displaymath}
%\left(
%\begin{array}{cc}
%\frac{\partial \changed{X_s^1(t;x)}{X_t^{1s}}}{\partial x_1} & \frac{\partial \changed{X_s^1(t;x)}{X_t^{1s}}}{\partial x_2}\\
%\frac{\partial \changed{X_s^1(t;x)}{X_t^{2s}}}{\partial x_1} & \frac{\partial
%\changed{X_s^1(t;x)}{X_t^{2s}}}{\partial x_2}
%\end{array}
%\right)= \exp\left[ \frac{1}{\nu}\left(\intl_s^t\left(
%\begin{array}{cc}
%-v_2(X_{\tau}^s)\,dX_{\tau}^{2s} & v_1(X_{\tau}^s)\,dX_{\tau}^{2s}\\
%v_2(X_{\tau}^s)\,dX_{\tau}^{1s} & -v_1(X_{\tau}^s)\,dX_{\tau}^{1s}
%\end{array}
%\right)- \intl_s^t\left(
%\begin{array}{cc}
%|v_2|^2 & |v_1|^2\\
%|v_2|^2 & |v_1|^2
%\end{array}
%\right)\,d\tau \right) \right] \del{\times}
%\end{displaymath}
%\del{\begin{displaymath}
%\left(
%\begin{array}{cc}
%1 & 0\\
%0 & 1
%\end{array}
%\right).
%\end{displaymath}}
\end{corollary}
\begin{proof}[Proof of Corollary \ref{cor:2DGradientFormula-1}]
We have by definition of the flow $(\changed{X_s(t;x)}{X_t^s})$,
$0\leq s\leq t\leq T$ that
$$
d\changed{X_s^1(t;x)}{X_t^{1s}}=\sqrt{2\nu}(\cos{\frac{\phi}{\nu}}(\changed{X_s(t;x)}{X_t^s})\,dw_t^1-\sin{\frac{\phi}{\nu}}(\changed{X_s(t;x)}{X_t^s})\,dw_t^2),
$$
$$
d\changed{X_s^1(t;x)}{X_t^{2s}}=\sqrt{2\nu}(\sin{\frac{\phi}{\nu}}(\changed{X_s(t;x)}{X_t^s})\,dw_t^1+\cos{\frac{\phi}{\nu}}(\changed{X_s(t;x)}{X_t^s})\,dw_t^2),
$$
$$
X_s(s;x)=x,x\in\Rnu^2.
$$
Taking derivative of the flow $(\changed{X_s(t;x)}{X_t^s})$,
$0\leq s\leq t\leq T$ with respect to initial condition $x$ we get
for the first component of the flow
\begin{displaymath}
d\left(
\begin{array}{c}
\frac{\partial \changed{X_s^1(t;x)}{X_t^{1s}}}{\partial x_1} \\
\frac{\partial \changed{X_s^1(t;x)}{X_t^{1s}}}{\partial x_2}
\end{array}
\right)=\sqrt{2\nu} \left(
\begin{array}{c}
(-\frac{1}{\nu}\sin{\frac{\phi}{\nu}}(\changed{X_s(t;x)}{X_t^s})\,dw^1_t-\frac{1}{\nu}\cos{\frac{\phi}{\nu}}(\changed{X_s(t;x)}{X_t^s})\,dw_t^2)(\frac{\partial\phi}{\partial
x_1}\frac{\partial \changed{X_s^1(t;x)}{X_t^{1s}}}{\partial x_1}
+\frac{\partial\phi}{\partial x_2}\frac{\partial \changed{X_s^2(t;x)}{X_t^{2s}}}{\partial x_1}) \\
(-\frac{1}{\nu}\sin{\frac{\phi}{\nu}}(\changed{X_s(t;x)}{X_t^s})\,dw^1_t-\frac{1}{\nu}\cos{\frac{\phi}{\nu}}(\changed{X_s(t;x)}{X_t^s})\,dw_t^2)(\frac{\partial\phi}{\partial
x_1}\frac{\partial \changed{X_s^1(t;x)}{X_t^{1s}}}{\partial x_2}
+\frac{\partial\phi}{\partial x_2}\frac{\partial
\changed{X_s^2(t;x)}{X_t^{2s}}}{\partial x_2})
\end{array}
\right)=
\end{displaymath}
\begin{displaymath}
=\left(
\begin{array}{c}
-\frac{1}{\nu}\,d\changed{X_s^2(t;x)}{X_t^{2s}}(v_2\frac{\partial \changed{X_s^1(t;x)}{X_t^{1s}}}{\partial x_1}-v_1\frac{\partial \changed{X_s^2(t;x)}{X_t^{2s}}}{\partial x_1}) \\
-\frac{1}{\nu}\,d\changed{X_s^2(t;x)}{X_t^{2s}}(v_2\frac{\partial
\changed{X_s^1(t;x)}{X_t^{1s}}}{\partial x_2}-v_1\frac{\partial
\changed{X_s^2(t;x)}{X_t^{2s}}}{\partial x_2})
\end{array}
\right),
\end{displaymath}
where in the last inequality we have used that
$v=\nabla^{\bot}\phi$ and definition of the flow. Similarly we can
get an equation for the gradient of the second component of the
flow. The result follows.
\end{proof}

\begin{proposition}\label{prop:FeynmanRepr}
Suppose that the vector field $v:\mathbb{R}^2\to \mathbb{R}^2$ is
of $C_0^{\infty}$ class and divergence free, i.e. $\diver v=0$.
Let $\changed{X_s(t;x)}{X_t^s(x)}$, $0\leq s\leq t\leq T$ be  the
flow corresponding to
 equation \eqref{eqn:flow-4}. Identifying $\mathbb{C}$ with $\mathbb{R}^2$ in the usual way, i.e.  $z=x_1+\imath x_2$, $x=(x_1,x_2)$, we can define a flow  $Z_s(t;z)$, $0\leq s\leq t\leq T$, $z\in \mathbb{C}$  by
$Z_s(t;z)=\changed{X_s^1(t;x)}{X_t^{1s}}+\imath \changed{X_s^1(t;x)}{X_t^{2s}}$.\\
If $F_0\in C_0^{\infty}(\Rnu^2)$ and
$F:[0,T]\times\Rnu^2\to\Rnu^2$ is a solution of equation
\eqref{eqn:eq-3}  such that for some $\beta>0$ and any smooth
closed loop $\Gamma$ condition \eqref{eqn:UniIntegrCond} is
satisfied, then
\begin{equation}
\mathbf{F}(t,z)=\mathrm{P}[\mathbb{E}(\overline{\mathbf{F}_0}(Z_{T-t}(T;z)))\frac{\partial
Z_{T-t}(T;z)}{\partial
\overline{z}}+\mathbf{F}_0(Z_{T-t}(T;z)))\frac{\partial
\overline{Z_{T-t}(T;z)}}{\partial
\overline{z}}],\label{eqn:FeynmanRepr}
\end{equation}
where $\mathbf{F}(t,z)=F^1(t,x)+\imath F^2(t,x)$ and
$\mathbf{v}(t,z)=v^1(t,x)+\imath v^2(t,x)$.\\
Moreover,  $\frac{\partial Z_s(t;z)}{\partial \overline{z}}$,
$\frac{\partial \overline{Z_s(t;z)}}{\partial \overline{z}}$
satisfy the following system of equations:
\begin{eqnarray}
d(\frac{\partial Z_s(t;z)}{\partial
\overline{z}})&=&\frac{1}{2\nu}(\mathbf{v}(t,Z_s(t;z))\frac{\partial
\overline{Z_s(t;z)}}{\partial
\overline{z}}-\overline{\mathbf{v}}(t,Z_s(t;z))\frac{\partial
Z_s(t;z)}{\partial \overline{z}})\,dZ_s(t;z)\nonumber\\
d(\frac{\partial \overline{Z_s(t;z)}}{\partial
\overline{z}})&=&\frac{1}{2\nu}(\overline{\mathbf{v}}(t,Z_s(t;z))\frac{\partial
Z_s(t;z)}{\partial
\overline{z}}-\mathbf{v}(t,Z_s(t;z))\frac{\partial
\overline{Z_s(t;z)}}{\partial
\overline{z}})\,d\overline{Z_s(t;z)}\nonumber\\
\label{eqn:FeynmanReprGradFlow} \frac{\partial Z_s(s;z)}{\partial
\overline{z}}&=&0,\;\frac{\partial \overline{Z_s(s;z)}}{\partial
\overline{z}}=1,
\end{eqnarray}
where $\overline{\cdot}$ is a complex conjugation.
\end{proposition}
\begin{proof}[Proof of Proposition \ref{prop:FeynmanRepr}]
Definition of the flow \eqref{eqn:flow-4} can be reformulated as
follows
\begin{equation}
\label{eqn:CompFlowDef} \left\{
\begin{array}{rcl}
dZ_s(t;z)(z,\overline{z})&=&\sqrt{2\nu}e^{\frac{\imath\phi(Z_s(t;z),\overline{Z_s(t;z)})}{\nu}}\,dW(t)^{\mathbb{C}},\\
Z_s(s;z)&=&z,
\end{array}
\right.
\end{equation}
where $W(t)^{\mathbb{C}}=W(t)^1+\imath W(t)^2$- wiener process in
$\mathbb{C}$. Then equation \eqref{eqn:FeynmanReprGradFlow}
immediately follow from definition \eqref{eqn:CompFlowDef}.
Formula \eqref{eqn:FeynmanRepr} is simply rewriting of formula
\eqref{eqn:FeynKacFormula}.
\end{proof}
%\begin{corollary}
%\mathbb{E}|\frac{\partial Z_s(t;z)}{\partial \overline{z}}|^2\leq
%e^{\frac{C}{\nu}\intl_s^t\Vert v(\tau)\Vert _{L^{\infty}}^2\,d\tau}
%\end{corollary}

\begin{remark}
Theorem \ref{thm:2DFlowDifferentLaw} %and
%\Blue{Instead \eqref{eqn:ContourItoFormula} you should use something else, as the latter is not at all a FK type %formula. Changed.} %formula
%\eqref{eqn:ContourItoFormula}
indicates the difference between the passive scalar  advection
equation and the vector advection equation. In the former case the
Feynman-Kac type formula does not contain a  gradient of the flow
 and hence the  solution  is completely
determined by the law of flow itself. Since the rotation of the
Brownian Motion does not change the law of the flow, we cannot
employ  the same trick for the scalar advection equation as we did
for the vector advection equation.
\end{remark}
\begin{question}
In connection with Theorem \ref{thm:2DFlowDifferentLaw} we can ask
if it possible to give a direct proof  (not through formula
\eqref{eqn:ContourItoFormula}) of the fact that the limit as
$\nu\to 0$ \del{ \eqref{eqn:FeynKacFormula-2}} exists and the
limit is a solution to the  2D Euler equations?
\end{question}

\subsection{Examples of nontrivial flows in $\mathbb{R}^3$} \label{example-3D}

In this subsection we provide nontrivial examples of the flows
which can be used in the Feynman-Kac type formula
\eqref{eqn:FeynKacFormula} in three dimensional case.

We will need the following definitions. Let $\hat{\cdot}$ be the
so called  hat-map   linear isomorphism defined by
$$
\hat{\cdot}:\Rnu^3\ni\left(
  \begin{array}{c}
    x_1 \\
    x_2 \\
    x_3 \\
  \end{array}
\right)\mapsto
 \left(
  \begin{array}{ccc}
    0 & -x_3 & x_2 \\
    x_3 & 0 & -x_1 \\
    -x_2 & x_1 & 0 \\
  \end{array}
\right)\in\mathfrak{so}(3),
$$
where $\mathfrak{so}(3)$ is the Lie algebra of antisymmetric
matrices.
 Let also $SO(3)$ be the Lie group of orthogonal matrices with determinant equal to one and let
$\exp:\mathfrak{so}(3)\ni A\mapsto e^A\in  SO(3)$ be the standard
exponential map. Let us  notice that this map is a surjection.%, see
%e.g. \Blue{Reference needed}

%\Red{The definition of the map $BCH$ relies on the fact that
%$\exp$ above is a bijection? Is this true? Now I understood that it is not true (rotation on the angle 2$\pi$), and, %therefore, I have changed definition a little; But I still need that the exponential map is surjection. I have found %this fact in internet. Though I was unable to pinpoint exact reference in some book. If you know classical book, %send me a reference.}
Denote $S=\ker (\exp)$. Define a map
$BCH:\mathfrak{so}(3)\times\mathfrak{so}(3)\to\mathfrak{so}(3)/S$
by
$$
\exp(BCH(\hat{u},\hat{v}))=\exp(\hat{u})\exp(\hat{v}),\hat{u},\hat{v}\in
\mathfrak{so}(3).
$$
%We will now find the  exact type of the "correction"  appearing in
%formula \eqref{eqn:ContourItoFormula}.
%\think{It's not clear to me what do you mean by a correction.}
%\Green{
Now we will find different form of the term
\eqref{eqn:correction_term} appearing in formula
\eqref{eqn:ContourItoFormula} due to diffusion coefficient
$\sigma$ of the flow $X_{\cdot}$.%}

\begin{proposition}\label{prop:CorrectionTerm-1} Let
$a\in C^1([0,T]\times\Rnu^3,\Rnu^3)$ and  a map $\sigma$ is
defined by  $\sigma:[0,T]\times\Rnu^3\ni (t,x) \mapsto
\exp(\widehat{a(t,x)})\in SO(3)$.   If $|a|(t,x)\not=0$, then
\begin{eqnarray}
\suml_m\sigma^{\cdot m}\frac{\partial \sigma^{\cdot m}}{\partial
x_k}=(1-\cos{|a|})\widehat{b\times\frac{\partial b}{\partial
x_k}}+\sin{|a|}\widehat{\frac{\partial b}{\partial
x_k}}+\widehat{b}\frac{\partial |a|}{\partial x_k},
\label{eqn:CorrTerm-1}
\end{eqnarray}
where $\vec{b}=\frac{\vec{a}}{|a|}$. If $|a|(t,x)=0$ then
\begin{equation}
\suml_m\sigma^{\cdot m}\frac{\partial \sigma^{\cdot m}}{\partial
x_k}=\frac{\partial \vec{a}}{\partial x_k}.\label{eqn:CorrTerm-11}
\end{equation}
\end{proposition}
\begin{remark}
We can notice that the right side if equality
\eqref{eqn:CorrTerm-1} can be rewritten as follows
$$
(1-\cos{|a|})\widehat{b\times\frac{\partial b}{\partial
x_k}}+(\sin{|a|}-|a|)\widehat{\frac{\partial b}{\partial
x_k}}+\frac{\partial \vec{a}}{\partial x_k}.
$$
Therefore it converges  to $\frac{\partial \vec{a}}{\partial x_k}$
when $|a|\to 0,|a|\not=0$. Hence, in the following considerations
we will not to single out  the case of $|a|(t,x)=0$.
\end{remark}
\begin{proof}[Proof of Proposition \ref{prop:CorrectionTerm-1}]
If $a(t,x)=0$ then formula \eqref{eqn:CorrTerm-11} immediately
follows from definition of $\sigma$. Assume that $a(t,x)\not=0$.
We will use the following Baker-Campbell-Hausdorff formula in
$\mathfrak{so}(3)$, see e.g. \cite[p. 630]{Kenth}.
\begin{proposition}If $u,v\in\Rnu^3$ then
$$
BCH(\hat{u},\hat{v})=\al\hat{u}+\beta\hat{v}+\gamma[\hat{u},\hat{v}],
$$
where $[\hat{u},\hat{v}]$ denotes the commutator of $\hat{u}$ and
$\hat{v}$, and $\al$, $\beta$, and $\gamma$ are real constants
defined by
\begin{displaymath}
\al=\frac{\sin^{-1}(d)}{d}\frac{a_1}{\theta},\beta=\frac{\sin^{-1}(d)}{d}\frac{b_1}{\phi},\gamma=\frac{\sin^{-1}(d)}{d}\frac{c_1}{\theta\phi},\label{eqn:CoeffBCHFormulae}
\end{displaymath}
where $a_1$, $b_1$, $c_1$ and $d$ are defined as
\begin{eqnarray*}
a_1 &=& \sin{\theta}\cos^2(\phi/2)-\sin{\phi}\sin^2(\theta/2)\cos{\angle(u,v)},  \\
b_1 &=& \sin{\phi}\cos^2(\theta/2)-\sin{\theta}\sin^2(\phi/2)\cos{\angle(u,v)}, \\
c_1 &=& \frac{1}{2}\sin(\theta)\sin(\phi)-2\sin^2(\theta/2)\sin^2(\phi/2)\cos{\angle(u,v)}, \\
d &=& \sqrt{a_1^2+b_1^2+2a_1b_1\cos{\angle(u,v)}+c_1^2\sin^2{\angle(u,v)}}. \\
\end{eqnarray*}
In the above formulae $\theta=|u|$, $\phi=|v|$, and $\angle(u,v)$
is the angle between the two vectors $u$ and $v$.
\end{proposition}
We have
\begin{eqnarray*}
\suml_m\sigma^{\cdot m}\frac{\partial \sigma^{\cdot m}}{\partial
x_k}&=&\exp(-\hat{a})\frac{\partial}{\partial x_k} \exp(\hat{a})=
\exp(-\hat{a})\times \liml_{\delta\to
0}\frac{1}{\delta}(\exp(\hat{a}(x+\delta
e_k))-\exp(\hat{a}(x)))\\
&=& \liml_{\delta\to
0}\frac{1}{\delta}(\exp(-\hat{a})\exp(\hat{a}(x+\delta
e_k))-\id)\\ &=& \liml_{\delta\to
0}\frac{1}{\delta}(\exp(BCH(-\hat{a},\hat{a}(x+\delta e_k)))-\id)=
\liml_{\delta\to 0}\frac{BCH(-\hat{a},\hat{a}(x+\delta
e_k))}{\delta}\\ &=& \liml_{\delta\to
0}\frac{\al(\delta)(-\hat{a}(x))+\beta(\delta)\hat{a}(x+\delta
e_k)+\gamma(\delta)[-\hat{a}(x),\hat{a}(x+\delta
e_k)]}{\delta}=(*),
\end{eqnarray*}
where in the last equality we have used Proposition
\ref{prop:CorrectionTerm-1} with $u=-\hat{a}(x)$,
$v=\hat{a}(x+\delta e_k)$. Therefore,
\begin{eqnarray*}
(*)&=&\liml_{\delta\to 0}\beta(\delta)\frac{\hat{a}(x+\delta
e_k)-\hat{a}}{\delta}+\hat{a}(x)\liml_{\delta\to
0}\frac{\beta(\delta)-\al(\delta)}{\delta}
\\&&\hspace{2truecm}
-\,\liml_{\delta\to
0}\gamma(\delta)[\hat{a}(x),\frac{\hat{a}(x+\delta
e_k)-\hat{a}(x)}{\delta}]\\
&=&\frac{\partial \hat{a}}{\partial x_k}\liml_{\delta\to
0}\beta(\delta) + \hat{a}\liml_{\delta\to
0}\frac{\beta(\delta)-\al(\delta)}{\delta}-
\widehat{(a\times\frac{\partial a}{\partial x_k})}\liml_{\delta\to
0}\gamma(\delta)
\end{eqnarray*}
So, we need to calculate the following three limits.
$$
(i)=\liml_{\delta\to 0}\beta(\delta), (ii)=\liml_{\delta\to
0}\frac{\beta(\delta)-\al(\delta)}{\delta}, (iii)=\liml_{\delta\to
0}\gamma(\delta).
$$
From \eqref{eqn:CoeffBCHFormulae} follows that we need to
calculate asymptotics of $a_1(\delta)$, $b_1(\delta)$,
$c_1(\delta)$, $d(\delta)$, $\delta\to 0$. We have
$$
\theta=|a|(x), \phi=|a|(x+\delta
e_k)=|a|(x)+\delta\frac{\partial}{\partial
x_k}|a|+\underline{o}(\delta),
$$
$$
\cos(\angle(u,v))=\frac{(-a(x),a(x+\delta
e_k))}{|a|(x)|a|(x+\delta e_k)}=-1+\bar{o}(\delta^2)
$$
\begin{eqnarray}
a_1=\sin |a|(\frac{1+\cos|a|(x+\delta e_k)}{2})-\sin |a|(x+\delta
e_k)\times\nonumber\\
(\frac{1-\cos|a|(x)}{2})(-1+\bar{o}(\delta^2))=\nonumber\\
\sin |a|(\frac{1+\cos(|a|+\delta\frac{\partial}{\partial
x_k}|a|)}{2})+(\frac{1-\cos|a|}{2})\times\nonumber\\
\sin(|a|+\delta\frac{\partial}{\partial
x_k}|a|)+\bar{o}(\delta^2)=\nonumber\\
\frac{\sin |a|}{2}(1+\cos|a|-\sin |a|\frac{\partial}{\partial
x_k}|a|\delta)+\nonumber\\
(\frac{1-\cos|a|}{2})(\sin|a|+\cos|a|\frac{\partial}{\partial
x_k}|a|\delta)+\bar{o}(\delta^2)=\nonumber\\
=\sin|a|(x)-\frac{1}{2}(1-\cos |a|)\frac{\partial}{\partial
x_k}|a|\delta+\bar{o}(\delta^2)\label{eqn:a1asymp}
\end{eqnarray}
Similarly,
\begin{eqnarray}
b_1=\sin |a|(x+\delta e_k)(\frac{1+\cos|a|}{2})-\sin
|a|\times\nonumber\\
(\frac{1-\cos|a|(x+\delta e_k)}{2})(-1+\bar{o}(\delta^2))=\nonumber\\
\sin(|a|+\delta\frac{\partial}{\partial
x_k}|a|)(\frac{1+\cos|a|}{2})+\nonumber\\\sin
|a|(\frac{1-\cos(|a|+\delta\frac{\partial}{\partial
x_k}|a|)}{2})+\bar{o}(\delta^2)=\nonumber\\
(\sin|a|+\cos|a|\frac{\partial}{\partial
x_k}|a|\delta)(\frac{1+\cos|a|}{2})+ \nonumber\\
\frac{1}{2}\sin |a|(\cos
|a|-\delta\sin |a|\frac{\partial}{\partial x_k}|a|)+\bar{o}(\delta^2)=\nonumber\\
\sin |a|+\frac{1}{2}(1+\cos |a|)\frac{\partial}{\partial
x_k}|a|\delta+\bar{o}(\delta^2)\label{eqn:b1asymp}
\end{eqnarray}
\begin{eqnarray}
c_1=1-\cos|a|+\bar{o}(\delta)\label{eqn:c1asymp}
\end{eqnarray}
\begin{eqnarray}
d=\bar{o}(\delta)\label{eqn:dasymp}
\end{eqnarray}
From \eqref{eqn:a1asymp},\eqref{eqn:b1asymp},\eqref{eqn:c1asymp}
and \eqref{eqn:dasymp} we get
\begin{eqnarray*}
(iii)=\liml_{\delta\to
0}\frac{\sin^{-1}(d)}{d}\frac{c_1}{|a|(x)|a|(x+\delta
e_k)}=\frac{1-\cos |a|}{|a|^2}
\end{eqnarray*}
\begin{eqnarray*}
(i)=\liml_{\delta\to
0}\frac{\sin^{-1}(d)}{d}\frac{b_1}{|a|(x+\delta e_k)}=\frac{\sin
|a|}{|a|}
\end{eqnarray*}
\begin{eqnarray*}
(ii)=\liml_{\delta\to 0}\frac{\sin^{-1}(d)}{d}\frac{1}{\delta}
(\frac{\sin |a|+\frac{1}{2}(1+\cos |a|)\frac{\partial}{\partial
x_k}|a|\delta+\bar{o}(\delta^2)}{|a|(x+\delta
e_k)}-\\
\frac{\sin|a|(x)-\frac{1}{2}(1-\cos |a|)\frac{\partial}{\partial
x_k}|a|\delta+\bar{o}(\delta^2)}{|a|})=\frac{|a|-\sin
|a|}{|a|^2}\frac{\partial}{\partial x_k}|a|
\end{eqnarray*}
Thus, we get
\begin{eqnarray}
\exp(-\hat{a})\frac{\partial}{\partial x_k}
\exp(\hat{a})=\frac{\sin |a|}{|a|}\frac{\partial\hat{a}}{\partial
x_k}+\frac{|a|-\sin
|a|}{|a|^2}\frac{\partial}{\partial x_k}|a|\hat{a}+\\
\frac{\cos |a|-1}{|a|^2}\widehat{a\times\frac{\partial a}{\partial
x_k}}\label{eqn:CorrTerm-2}
\end{eqnarray}
If we put $b=\frac{a}{|a|}$ and insert it in
\eqref{eqn:CorrTerm-2} we get \eqref{eqn:CorrTerm-1}.
\end{proof}

Everywhere below we assume that $v\in L^1(0,T;C_{b}^{2,\al}(\Rn
,\Rn))$ for some $\al\in (0,1)$.
\begin{corollary}\label{cor:3DcirculationDetails}
Let $(\changed{X_s(t;x)}{X_t^s})$, $0\leq s\leq t\leq T$ be the
stochastic flow corresponding to

\begin{eqnarray}
d\changed{X_s(t;x)}{X_t^s} &=&
v(t,\changed{X_s(t;x)}{X_t^s})\,dt+\sqrt{2\nu}\sigma_1(t,\changed{X_s(t;x)}{X_t^s})\,dW(t),\label{eqn:flow-5}\\
X_s(s;x) &=& x,\nonumber
\end{eqnarray}
where $\sigma_1(t,x)=\exp(\hat{a})(t,x)$, $b=\frac{a}{|a|}\in
S(2)$. Then for all $s,t\in [0,T]$ such that $s\leq t$,
\begin{eqnarray}
&&
\intl_{X_{T-s}(t;\Gamma)}\suml_{k=1}^{n}F^k(T-t,x)\,dx_k=\intl_{\Gamma}\suml_{k=1}^{n}F^k(s,x)\,dx_k
\nonumber\\ &+&
\intl_{T-s}^t\intl_{X_{T-s}(\tau;\Gamma)}\suml_{k=1}^{n}\left(\frac{\partial
F^k}{\partial t}\right.+ \left.\suml_{j=1}^{n}v^{j}(\frac{\partial
F^k}{\partial x_j}-\frac{\partial F^j}{\partial
x_k})+\nu\triangle F^k\right)\,dx_kd\tau \nonumber\\
&+&\nu\intl_{T-s}^t\intl_{X_{T-s}(\tau;\Gamma)}\big(\curl
F,(1-\cos{|a|})b\times\frac{\partial b}{\partial
x_k}+\sin{|a|}\frac{\partial b}{\partial x_k}+b\frac{\partial
|a|}{\partial x_k}\big)\,dx_kd\tau\label{eqn:ContourItoFormula-2}\\
&+&\sqrt{2\nu}
\intl_{T-s}^t\intl_{X_{T-s}(\tau;\Gamma)}\suml_{k=1}^n\left(\suml_{i,l=1}(\frac{\partial
F^k}{\partial x_i}-\frac{\partial F^i}{\partial
x_k})\sigma_1^{il}\right)\,dx_k\,dW_{\tau}^{l}.\nonumber
\end{eqnarray}
\end{corollary}
\begin{proof}[Proof of Corollary \ref{cor:3DcirculationDetails}]
Immediately follows from Proposition \ref{prop:CorrectionTerm-1}
and identity
$$
\suml_{i,j}\frac{\partial F^i}{\partial x_j}(\hat{a})^{ij}=(\curl
F,a).
$$
\end{proof}
\begin{remark}
The vector $b$ can be  interpreted as  the axis of rotation of
$\sigma$ and $\phi=|a|$ as  the angle of rotation.
\end{remark}
Now, we present a three dimensional analog of the two dimensional
result from Proposition \ref{prop:2DFlowTheSameLaw}.
\begin{proposition}\label{prop:3DFlowTheSameLaw-1}
Assume that $F_0\in C_0^{\infty}(\Rnu^3)$, $v\in
L^1(0,T;C_{b}^{2,\al}(\Rn ,\Rn))$, $\al\in (0,1)$, $v$ satisfies
condition \eqref{eqn:UniIntegrCond-2}, and $F\in L^{\infty}([0,T];
C^{2+\delta}(\Rnu^3,\Rnu^3))$ is a solution of equation
\eqref{eqn:eq-3}-\eqref{eqn:eq-3b} such that for some $\beta>0$
and any smooth closed loop $\Gamma$ condition
\eqref{eqn:UniIntegrCond} is satisfied. Let
$(\changed{X_s(t;x)}{X_t^s})$, $0\leq s\leq t\leq T$ be the
stochastic flow corresponding to
\begin{eqnarray}
d\changed{X_s(t;x)}{X_t^s} &=&
v(t,\changed{X_s(t;x)}{X_t^s})\,dt+\sqrt{2\nu}\sigma_1(t,\changed{X_s(t;x)}{X_t^s})\,dW(t),\label{eqn:flow-6}\\
X_s(s;x) &=& x,\nonumber
\end{eqnarray}
where $\sigma_1(t,x)=\exp(\hat{a}(t,x))$, $a=\curl F$. Fix $s\in
[0,T]$ and define a functions $Q_s:\mathbb{R}^3\to \mathbb{R}^3$
by $Q_s=\mathbb{E}(F_0(\changed{X_{T-s}(T;x)}{X_T^{T-s}(x)})\nabla
\changed{X_{T-s}(T;x)}{X_T^{T-s}(x)})$.

Then, $Q_s\in L^2(\mathbb{R}^3,\mathbb{R}^3)\cap
C^{\eps}(\mathbb{R}^3,\mathbb{R}^3)$,  $0<\eps<\alpha$ and
\begin{equation}\label{eqn:FeynKacFormula3D}
F(s,x)=[\mathrm{P}(Q_s)](x),\; x\in \mathbb{R}^3,\; s\in [0,T].\;
\end{equation}
\end{proposition}
\begin{proof}[Proof of Proposition \ref{prop:3DFlowTheSameLaw-1}]
In view of   \cite[ Theorem 4.6.5, p.173]{Kunita} we infer that
there exists solution  $X_s(t;x),0\leq s\leq t\leq T$ of problem
\eqref{eqn:flow-6} and $X_s(t;x),0\leq s\leq t\leq T$ is a flow of
$C^1$--diffeomorphisms. Furthermore, it is
$C^{1+\eps}(\Rn,\Rn)$-valued process for any $0<\eps<\delta$.

Moreover, it  follows from Theorems 3.3.3, p.94 and  4.6.5, p. 173
therein \del{of \cite{Kunita}} that for all $s\in [0,T]$, $Q_s\in
C^{\eps}(\Rn,\Rn)$, $0<\eps<\delta$. Let us fix $s\in [0,T]$. We
will show now that $Q_s\in L^2(\Rn,\Rn)$. Since by Corollary 4.6.7
p. 175 of \cite{Kunita} that there exists a positive constant $C$
such that
$$
\sup_{x\in\Rn}\mathbb{E}|\nabla
\changed{X_{T-s}(T;x)}{X_T^{T-s}(x)}|^2\leq C,
$$
 by the H\"older inequality we infer that
\begin{eqnarray}
\intl_{\Rn}|Q_s(x)|^2\,dx\leq
\intl_{\Rn}\mathbb{E}|F_0(\changed{X_{T-s}(T;x)}{X_T^{T-s}(x)})|^2\mathbb{E}|\nabla
\changed{X_{T-s}(T;x)}{X_T^{T-s}(x)}|^2\,dx\nonumber\\
\leq
C\intl_{\Rn}\mathbb{E}|F_0(\changed{X_{T-s}(T;x)}{X_T^{T-s}(x)})|^2\,dx.\label{eqn:auxIneq-1'}
\end{eqnarray}
Now it follows from Girsanov Theorem that
\begin{equation}
\intl_{\Rn}\mathbb{E}|F_0(\changed{X_{T-s}(T;x)}{X_T^{T-s}(x)})|^2\,dx
=\intl_{\Rn}\mathbb{\tilde{E}}(|F_0(x+\sqrt{2\nu}(W_T-W_{T-s}))|^2\mathcal{E}_{T-s}^T)\,dx,\label{eqn:auxEst-1'}
\end{equation}
where
$\mathcal{E}_{T-s}^T=e^{\intl_{T-s}^Tv(r,X_{T-s}(r;x))\,dW_r-1/2\intl_{T-s}^T|v(r,X_{T-s}(r;x))|^2\,dr}$
is a stochastic exponent. We can notice that
\begin{equation}
\mathbb{\tilde{E}}|\mathcal{E}_{T-s}^T|^2\leq
e^{2\intl_0^T|v(r)|_{L^{\infty}}(r)\,dr}\label{eqn:auxEst-2'}
\end{equation}
and, therefore, combining \eqref{eqn:auxIneq-1'},
\eqref{eqn:auxEst-1'} and \eqref{eqn:auxEst-2'} we get
\begin{equation}
\intl_{\Rn}|Q_s(x)|^2\,dx\leq
e^{\intl_0^T|v|_{L^{\infty}}(r)\,dr}\left(\intl_{\Rn}
|F_0|^4dx\right)^{\frac{1}{2}}.
\end{equation}
Now let us show that
$\intl_{X_{T-s}(t;\Gamma)}\suml_{k=1}^3F^k(T-t)\,dx_k,t\in[T-s,T]$
is a local martingale. It is enough to prove that the "correction"
term (due to nontrivial $\sigma_1$) in the formula
\eqref{eqn:ContourItoFormula-2} disappears.

Since $b=\frac{\curl F}{|\curl F|}$, $|b|=1$, $|a|=|\curl F|$ we
have
$$
(\curl F, \frac{\partial b}{\partial x_k})=|\curl
F|(b,\frac{\partial b}{\partial x_k})=0.
$$
Similarly,
$$
(\curl F, b\times\frac{\partial b}{\partial x_k})=|\curl F|(b,
b\times\frac{\partial b}{\partial x_k})=0,
$$
and
$$
(\curl F,b)\frac{\partial |\curl F|}{\partial
x_k}=\frac{1}{2}\frac{\partial |\curl F|^2}{\partial x_k}.
$$
%Hence, the rest of the proof is as in the proposition
%\ref{prop:FeynKacformula}.\Red{???}
\end{proof}
\begin{question}\label{quest:3DdifferentLawCase}
It would be interesting to generalize Theorem
\ref{thm:2DFlowDifferentLaw} to the three dimensional case. In
view of Corollary \ref{cor:3DcirculationDetails} in order to find
such generalization it is enough to prove that for any   solution
$F$  of equation \eqref{eqn:eq-3} with  $v$ being  the
corresponding  $C^{\infty}$)  vector field,
 there exists a triple
$(b,\phi,\psi)\in
(L^{\infty}([0,T],C^{\infty}(\Rnu^3,S^2)),L^{\infty}([0,T],C^{\infty}(\Rnu^3,S^1)),L^{\infty}([0,T],C^{\infty}(\Rnu^3,\Rnu)))$
such that
\begin{eqnarray}
(\cos\phi-1)(\curl F,b\times\frac{\partial b}{\partial
x_k})&+&\sin\phi(\curl F,\frac{\partial b}{\partial x_k})\nonumber\\
+(\curl F,b)\frac{\partial\phi}{\partial x_k} &+&\frac{\partial
\psi}{\partial x_k}=\frac{(v\times\curl F)^k}{\nu}\,\,\,,
k=1,2,3.\label{eqn:3DNonlinearPresentation}
\end{eqnarray}
We can notice that system \eqref{eqn:3DNonlinearPresentation} is
time independent in the sense that there are no time derivatives
of the unknown functions.  Therefore it is enough to consider the
system for every fixed time $t\in[0,T]$. If the solenoidal vector
field $v$ is two dimensional, i.e. $\diver v=0$, $v_3=0$ and  the
components $v_1$,  $v_2$ do not depend upon variable $x_3$, then
$b=(0,0,1)$, $\phi=\phi_1/\nu$, where $\phi_1$ is a stream
function for $v$, $\psi=0$, is a solution of the system
\eqref{eqn:3DNonlinearPresentation}, see Theorem
\ref{thm:2DFlowDifferentLaw}.
%system  \eqref{eqn:3DNonlinearPresentation} is undetermined
However, in the three dimensional case the problem is completely
open. One of the possibilities to narrow the problem is to
consider the case when
$F=u$ is a solution to %corresponding to the case of
the Navier-Stokes equations.
\end{question}
\begin{question}
 Another question connected with system
\eqref{eqn:3DNonlinearPresentation} is as follows. How do
variables $b$, $\phi$, $\psi$ depend upon $\nu$? Can one take the
$\nu$ to $0$ limit in the representation
\eqref{eqn:3DNonlinearPresentation}? But let us note that in the
two dimensional case  under the additional condition of
incompressibility $\diver v=0$, the  representation
\eqref{eqn:3DNonlinearPresentation} holds also in the limit
$\nu\to 0$. Indeed, in two dimensional case the stream function
corresponding to the vector field  $v$ exists because $\diver v=0$
and is independent of $F$ and $\nu$.
%Indeed, stream function has physical sense for
%two dimensional equation \eqref{eqn:eq-3} with $\nu=0$ and we have
%representation \eqref{eqn:3DNonlinearPresentation} for it.
\end{question}
\begin{remark}\label{rem:3DNonlinearPresentation-2}
The Question \ref{quest:3DdifferentLawCase} %system \eqref{eqn:3DNonlinearPresentation}
can be reformulated
in the following way. \\
\textbf{Problem A.} Find a $C^1$-class function $\sigma:\Rnu^3\to
SO(3)$ such that for any smooth closed loop $\Gamma$
\begin{equation}
\suml_{k}\intl_{\Gamma}\suml_{j}v^j\big(\frac{\partial
F^k}{\partial x_j}-\frac{\partial F^j}{\partial x_k}\big)\,dx_k=
\nu \suml_{k=1}^{n}\intl_{\Gamma} \suml_{j,l,m}\frac{\partial
F^j}{\partial x_l}\sigma^{lm}\frac{\partial \sigma^{jm}}{\partial
x_k} \,dx_k.\label{eqn:eq-1}
\end{equation}

Let $\wedge$ be the wedge product, see e.g.
\cite[p.79]{[Spivak-65]}, and denote %we introduce the  following
%notation
\begin{eqnarray*}
\al&=&\frac{1}{\nu}\suml_{j,k}v^j(\frac{\partial F^k}{\partial
x_j}-\frac{\partial F^j}{\partial x_k})\,dx_k,\\
w&=&\curl F.
\end{eqnarray*}
Suppose  $\sigma:\Rnu^3\to SO(3)$ is a $C^1$-class function. Let
us define a matrix valued function $A$,
\begin{equation}
A=d\sigma\sigma^{-1}.\label{eqn:Conn_Def}
\end{equation}
Then the matrix $A$ is antisymmetric and has the following form
\begin{equation}
A=\left(
\begin{array}{ccc}
    0 & -a_3 & a_2 \\
    a_3 & 0 & -a_1 \\
    -a_2 & a_1 & 0 \\
\end{array}
\right)\label{eqn:Conn_CompCond_1}
\end{equation}
where $a_i(x),i=1,2,3,x\in\Rnu^3$ are 1-forms. Moreover, $A$
satisfies system
\begin{equation}
dA+A\wedge A=0,\label{eqn:Conn_CompCond_2}
\end{equation}
or, in terms of 1-forms $a_i,i=1,2,3$, equivalently
\begin{eqnarray*}
  da_1 &=& a_3\wedge a_2 \\
  da_2 &=& a_1\wedge a_3 \\
  da_3 &=& a_2\wedge a_1.
\end{eqnarray*}
Furthermore, if arbitrary antisymmetric matrix $A$ of one-forms
satisfies \eqref{eqn:Conn_CompCond_2} then there exists
$\sigma:\Rnu^3\to SO(3)$ such that \eqref{eqn:Conn_Def} is
satisfied. Notice that the right part of formula \eqref{eqn:eq-1}
can be rewritten as follows
\begin{equation}
\nu\intl_{\Gamma}\suml_{i=1}^3w_ia_i. \label{eqn:Conn_Aux}
\end{equation}
Indeed,
\begin{equation*}
\suml_{k,m}\sigma^{\cdot m}\frac{\partial \sigma^{\cdot
m}}{\partial x_k}\,dx_k=d\sigma\sigma^{\bot}=d\sigma\sigma^{-1}=A.
\end{equation*}
Now we can rewrite formula \eqref{eqn:eq-1} as follows
\begin{equation}
\intl_{\Gamma}\al=-\intl_{\Gamma}\suml_{i=1}^3w_ia_i,
\end{equation}
Hence, we can reformulate the equation
\eqref{eqn:3DNonlinearPresentation} as follows
\begin{equation}
\suml_{i=1}^3w_ia_i=-\al+d\psi.\label{eqn:Conn_Eqn_1}
\end{equation}
Thus, Problem A can be solved in two stages. First, we need to
solve system
\begin{equation}
\label{eqn:Conn_Eqn_2} \left\{
\begin{array}{rcl}
  d a_1 &=& a_3\wedge a_2 \\
  d a_2 &=& a_1\wedge a_3 \\
  d a_3 &=& a_2\wedge a_1 \\
  \suml_{i=1}^3w_ia_i &=& -\al+d\psi.
\end{array}
\right.
\end{equation}
Then we need to find $\sigma:\Rnu^3\to SO(3)$ from equation
\eqref{eqn:Conn_Def}. Existence of such $\sigma$ follows from
first three equations of system \eqref{eqn:Conn_Eqn_2}.

Applying the exterior derivative operator $d$ to the last equation
of the system \eqref{eqn:Conn_Eqn_2} we can get rid of function
$\psi$ and thus we get equivalent system
\begin{equation}
\left\{
\begin{array}{rcl}
  da_1 &=& a_3\wedge a_2 \\
  da_2 &=& a_1\wedge a_3 \\
  da_3 &=& a_2\wedge a_1 \\
-d\al &=&   \suml_{i=1}^3dw_i\wedge a_i + w_1a_3\wedge a_2+
w_2a_1\wedge  a_3+w_3a_2\wedge a_1.
\end{array}
\right.
\end{equation}
This system can be reformulated in terms of matrix-valued 1-form
$A$ as follows:
\begin{equation}
\left\{
\begin{array}{rcl}
  dA+A\wedge A &=& 0 \\
  \tr(WA\wedge A+dW\wedge A) &=& 2\,d\al,
\end{array}
\right.
\end{equation}
where
\begin{displaymath}
W=\left(
\begin{array}{ccc}
    0 & -w_3 & w_2 \\
    w_3 & 0 & -w_1 \\
    -w_2 & w_1 & 0 \\
\end{array}
\right).
\end{displaymath}
Thus we have quadratic equation on the space of flat connections.
\end{remark}
Another application of Proposition \ref{prop:ContourConservation}
is a Feynman-Kac type formula for solutions of the following
equation
\begin{eqnarray}
\label{eqn:eq-4}
\frac{\partial F}{\partial t}&=&-\nu A_0 F+(v(T-\cdot)\cdot\nabla) F-(F\cdot\nabla)v(T-\cdot),\;t>0,\,x\in\Rn,\\
F(0)&=&F_0,
\end{eqnarray}
where $A_0$ is a Stokes operator, $F_0\in H$ and  $v$ satisfies
condition \eqref{eqn:IntegrCond-1}. For the simplicity sake we
formulate the result for $n=3$.
\begin{proposition}\label{prop:FeynmanKacSurface}
Let $v\in L^1(0,T;C_{b}^{2,\al}(\Rn ,\Rn))$ for some $\al\in
(0,1)$, $v$ satisfies condition \eqref{eqn:UniIntegrCond-2},
$(\changed{X_s(t;x)}{X_t^s})$, $0\leq s\leq t<\infty$ is the flow
corresponding to problem \eqref{eqn:flow-1}, $F_0\in
C_0^{\infty}(\Rn)$ and $F$ is a solution of equation
\eqref{eqn:eq-4} such that there exists $\beta>0$:
\begin{equation}
\mathbb{E}|\intl_{X_{T-s}(t;S)}F^1(T-t,x)\,dx_2\,dx_3+F^2(T-t,x)\,dx_3dx_1+F^3(T-t,x)\,dx_1dx_3|^{1+\beta}<\infty
\end{equation}
for any smooth surface $S\subset\Rnu^3$ with smooth boundary
$\Gamma$ and all $0\leq T-s\leq t\leq T$. Then it satisfies
\begin{eqnarray}
F^1(s,x)=\mathbb{E}[F_0^1(\changed{X_{T-s}(T;x)}{X_T^{T-s}(x)})(\frac{\partial
X_{T-s}^2(T;x)}{\partial x_2}\frac{\partial X_{T-s}^3(T;x)
}{\partial x_3}-\frac{\partial X_{T-s}^2(T;x)}{\partial
x_3}\frac{\partial
X_{T-s}^3(T;x)}{\partial x_2})\nonumber\\
+F_0^2(\changed{X_{T-s}(T;x)}{X_T^{T-s}(x)})(\frac{\partial
X_{T-s}^3(T;x) }{\partial x_2}\frac{\partial X_{T-s}^1(T;x)
}{\partial x_3}-\frac{\partial X_{T-s}^3(T;x)}{\partial
x_3}\frac{\partial X_{T-s}^1(T;x)}{\partial
x_2})\label{eqn:FeynmanKacsurface-1}\\
+F_0^3(\changed{X_{T-s}(T;x)}{X_T^{T-s}(x)})(\frac{\partial
X_{T-s}^1(T;x) }{\partial x_2}\frac{\partial X_{T-s}^2(T;x)
}{\partial x_3}-\frac{\partial X_{T-s}^1(T;x)}{\partial
x_3}\frac{\partial X_{T-s}^2(T;x)}{\partial x_2})\nonumber
\end{eqnarray}
\begin{eqnarray}
F^2(s,x)=\mathbb{E}[F_0^1(\changed{X_{T-s}(T;x)}{X_T^{T-s}(x)})(\frac{\partial
X_{T-s}^2(T;x) }{\partial x_3}\frac{\partial
X_{T-s}^3(T;x)}{\partial x_1}-\frac{\partial
X_{T-s}^2(T;x)}{\partial x_1}\frac{\partial
X_{T-s}^3(T;x)}{\partial x_3})\nonumber\\
+F_0^2(\changed{X_{T-s}(T;x)}{X_T^{T-s}(x)})(\frac{\partial
X_{T-s}^3(T;x) }{\partial x_3}\frac{\partial X_{T-s}^1(T;x)
}{\partial x_1}-\frac{\partial X_{T-s}^3(T;x)}{\partial
x_1}\frac{\partial X_{T-s}^1(T;x)}{\partial
x_3})\label{eqn:FeynmanKacsurface-2}\\
+F_0^3(\changed{X_{T-s}(T;x)}{X_T^{T-s}(x)})(\frac{\partial
X_{T-s}^1(T;x) }{\partial x_3}\frac{\partial X_{T-s}^2(T;x)
}{\partial x_1}-\frac{\partial X_{T-s}^1(T;x)}{\partial
x_1}\frac{\partial X_{T-s}^2(T;x)}{\partial x_3})\nonumber
\end{eqnarray}
\begin{eqnarray}
F^3(s,x)=\mathbb{E}[F_0^1(\changed{X_{T-s}(T;x)}{X_T^{T-s}(x)})(\frac{\partial
X_{T-s}^2(T;x) }{\partial x_1}\frac{\partial
X_{T-s}^3(T;x)}{\partial x_2}-\frac{\partial
X_{T-s}^2(T;x)}{\partial x_2}\frac{\partial
X_{T-s}^3(T;x)}{\partial x_1})\nonumber\\
+F_0^2(\changed{X_{T-s}(T;x)}{X_T^{T-s}(x)})(\frac{\partial
X_{T-s}^3(T;x) }{\partial x_1}\frac{\partial X_{T-s}^1(T;x)
}{\partial x_2}-\frac{\partial X_{T-s}^3(T;x)}{\partial
x_2}\frac{\partial X_{T-s}^1(T;x)}{\partial
x_1})\label{eqn:FeynmanKacsurface-3}\\
+F_0^3(\changed{X_{T-s}(T;x)}{X_T^{T-s}(x)})(\frac{\partial
X_{T-s}^1(T;x) }{\partial x_1}\frac{\partial X_{T-s}^2(T;x)
}{\partial x_2}-\frac{\partial X_{T-s}^1(T;x)}{\partial
x_2}\frac{\partial X_{T-s}^2(T;x)}{\partial x_1})\nonumber
\end{eqnarray}
\end{proposition}
\begin{proof}[Proof of Proposition \ref{prop:FeynmanKacSurface}]
The result follows from Proposition \ref{prop:FeynKacformula}.
Indeed, let $G\in L^{\infty}(0,T;L^2(\Rn,\Rn)\cap
C^{1+\eps}(\Rn,\Rn))$, $0<\eps<\al$ be a solution of equation
\eqref{eqn:eq-3}-\eqref{eqn:eq-3b}. Its existence follows from
Proposition \ref{prop:FeynKacformula}. Then $F=\curl G$ is a
solution of equation \eqref{eqn:eq-4}. For solution $G$ of
\eqref{eqn:eq-3} we have got representation by formula
\eqref{eqn:FeynKacFormula} of Feynman-Kac type. Integrating it
w.r.t. closed contour $\Gamma$ we get
\begin{equation}
\intl_{\Gamma}\suml_kG^k(s,x)\,dx_k=\mathbb{E}(\intl_{X_{T-s}(T;\Gamma)}\suml_kG_0^k(x)\,dx_k).
\end{equation}
Now, result immediately follows from Stokes Theorem.
\end{proof}
\begin{remark}
On an informal level, the Feynman-Kac type formula
\eqref{eqn:FeynmanKacsurface-1}-\eqref{eqn:FeynmanKacsurface-3} in
the case of $\nu=0$ can be seen as a solution   of the following
informal infinite dimensional first order PDE obtained by the
characteristics method. Indeed,  let us denote by $Y$ the set of
all smooth surfaces $S\subset \Rn$ with smooth boundary $\Gamma$.
Let  $TY$ be the set of all smooth vector fields on $Y$.
 If $F$ is a solution
of equation \eqref{eqn:eq-4} with parameters $\nu=0$ and $v\in
C_0^{\infty}([0,T]\times\Rn)$, then $\tilde{F}$ defined by
$$
\tilde{F}:[0,\infty)\times Y\ni
(t,S)\mapsto\intl_{S}(F(t,\cdot),\vec{n})\,d\sigma\in\Rnu,
$$
 is a solution to the following  equation
\begin{equation}
\frac{\partial\tilde{F}}{\partial
t}=D_{\tilde{v}}\tilde{F},\label{eqn:eq-5}
\end{equation}
where $D_{\tilde{v}}$ is directional derivative along the vector
field $\tilde{v}\in TY$ defined by
$$
Y\ni S\mapsto \bigcup_{x\in S}v(x)\in TY.
$$
Then, on a purely speculative level,  the solution to equation
\eqref{eqn:eq-5} obtained via the characteristics method
 %, (\Red{unclear}i.e. solution
%of \eqref{eqn:eq-5} is constant along with characteristics)
is exactly our Feynman-Kac type formula.
%\Red{Maybe we should add
%that this is purely speculative? Or is this rigorous?added in the
%beginning(in green)}
\end{remark}

\begin{remark}
 In a forthcoming publication  the authors will   consider the case of equations with less
regular velocity vector fields than those considered in the
current paper.  Transport equations with irregular velocity field
have been a subject of a great variety of works, see e.g. recent
works by Lions and  Di Perna \cite{[LionsDiPerna]},  Maniglia
\cite{[Maniglia]}, Bouchot, James and  Mancini
\cite{[BouchutJamesMancini2005]}, and references therein. Our plan
is to
 combine the results of Maniglia \cite{[Maniglia]} with
our work i.e. to find probabilistic representation of solution of
vector advection equation with irregular velocity and then study
the limit as  the viscosity $\nu$ converges to $0$.
\end{remark}

%Authors plan to investigate in the forthcoming paper the validity
%of the Feynman-Kac formulas provided here for non smooth
%velocities such as, for example, in the works of \cite{[Maniglia]}
%and \cite{[BouchutJamesMancini2005]}. Maniglia ()

\section{Proofs of results from section \ref{sec:Exist&Uniqueness}}
\begin{proof}[Proof of Proposition \ref{prop:ExistUniqSolution}]
\begin{trivlist}
\item[(i)] The  proof will be divided into three parts a), b), c).
\item[a)] Let us consider a special case when $v\in
L^{\infty}(0,T;\mathbb{L}^{3+\delta_0}(D))$. We will use Theorem
\ref{thm:LionsPDEexun} with Gelfand triple $V\subset H\cong
H'\subset V^\prime$. Denote $A(t)=\nu A+B(v(t),\cdot)$. We need to
check whether the conditions \eqref{eqn:BoundednCond} and
\eqref{eqn:CoercivCond} are satisfied. We have,
\begin{equation}
\label{eqn:QuadForm-1} \langle
A(t)f,f\rangle_{V^\prime,V}=\nu\tilde{a}(f,f)+\langle
B(v(t),f),f\rangle_{V^\prime,V},\; f\in V.
\end{equation}
The second term on the RHS of the equality \eqref{eqn:QuadForm-1}
from \eqref{eqn:NlinEstimate2} can be estimated as follows
\begin{eqnarray}
|\langle B(v(t),f),f\rangle_{V^\prime,V}|\leq
\frac{1}{2}\Vert f\Vert _{V}^2+\frac{1}{2}(\eps^{1+\delta_0/3}\Vert f\Vert _{V}^2\nonumber\\
+\frac{C_{\delta_0}}{\eps^{1+3/\delta_0}}|v(t)|_{\mathbb{L}^{3+\delta_0}(D)}^{2+\frac{6}{\delta_0}}|f|_{H}^2),\eps>0.\label{eqn:NlinEst-1}
\end{eqnarray}
Thus from the inequality  \eqref{eqn:NlinEst-1} and the continuity
of form $\tilde{a}$ we infer that,
\begin{equation}
\Vert A(t)\Vert _{\mathcal{L}(V,V^\prime)}\leq C
\nu+C_2|v(t)|_{\mathbb{L}^{3+\delta_0}(D)}.\label{eqn:OperatorANormEst}
\end{equation}
The coercivity assumption  \eqref{eqn:CoercivCond} also follows
from the inequality \eqref{eqn:NlinEstimate2}. Indeed, for $f\in
V$, $t\in [0,T]$ we have
\begin{eqnarray*}
|\langle A(t)f,f\rangle_{V^\prime,V}|=|\nu\tilde{a}(f,f)+\langle B(v(t),f),f\rangle_{V^\prime,V}|\geq\\
\frac{\nu}{2}\Vert f\Vert _V^2-\frac{C}{\nu}(\eps^{1+\delta_0/3}
\Vert f\Vert
_{V}^2+\frac{C_{\delta_0}}{\eps^{1+3/\delta}}|v(t)|_{\mathbb{L}^{3+\delta_0}(D)}^{2+\frac{6}{\delta_0}}|f|_{H}^2).
\end{eqnarray*}
By choosing $\eps>0$ such that
$\frac{\nu}{2}-\frac{C}{\nu}\eps^{1+\delta_0/3}>0$ we conclude the
proof of the coercivity condition \eqref{eqn:CoercivCond}. Thus,
by the Theorem \ref{thm:LionsPDEexun}, first statement of the
Proposition follows. \item[b)]To prove Proposition in the general
case we will show an energy inequality for solutions of equation
(\ref{eqn:DualEquation}-\ref{in_v}) when $v\in
L^{\infty}(0,T;\mathbb{L}^{3+\delta_0}(D))$. From step (a) we know
that a solution $F\in L^2(0,T;V)$ such that $F^\prime\in
L^2(0,T;V^\prime)$ exists and unique. Then, from Lemma
\ref{lem:ContinuityLem} and equality \eqref{eqn:bFormFormula-1} we
infer that
%Multiplying (\ref{eqn:DualEquation}-\ref{in_v}) on $F$ and integrating we get
\begin{eqnarray*}
\frac{1}{2}\frac{d}{dt}|F|_H^2&=&-\nu\Vert F\Vert _{V}^2+\langle f,F\rangle_{V^\prime,V}-\langle B(v,F),F\rangle_{V^\prime,V}\\
&=&-\nu\Vert F\Vert _{V}^2+\langle f,F\rangle_{V^\prime,V}+(\curl
F,v\times F)_H.
\end{eqnarray*}
Therefore, by applying the Young inequality, we infer that
\begin{eqnarray*}
|F(t)|_H^2&+&2\nu\intl_0^t|F(s)|_{V}^2\,ds-\intl_0^t(\curl
F(s),v(s)\times F(s))_H\,ds\\
&=&|F(0)|_H^2+\intl_0^t\langle f(s),F(s)\rangle_{V^\prime,V}\,ds\\
&\leq& |F(0)|_H^2+
\frac{\nu}{2}\intl_0^t|F(s)|_{V}^2\,ds+\frac{C}{\nu}\intl_0^t|f(s)|_{V^\prime}^2\,ds.
\end{eqnarray*}
The term $\intl_0^t(\curl F(s),v(s)\times F(s))_H\,ds$ can be
estimated as follows:
\begin{eqnarray}
|&\intl_0^t&(\curl F(s),v(s)\times F(s))_H\,ds|\nonumber\\
&\leq&\frac{\nu}{4}\intl_0^t|\curl F|_H^2\,ds+
\frac{C}{\nu}\intl_0^t|v(s)\times
F(s)|_H^2\,ds\nonumber\\
&\leq&\frac{\nu}{4}\intl_0^t|\curl F|_H^2\,ds+
\frac{C}{\nu}\intl_0^t(\eps^{1+\delta_0/3}|F(s)|_{V}^2+
\frac{C_{\delta_0}}{\eps^{1+\delta_0/3}}|v(s)|_{\mathbb{L}^{3+\delta_0}}^{2+6/\delta_0}|F(s)|_H^2)\,ds\nonumber\\
&\leq&(\frac{\nu}{4}+\frac{C}{\nu}\eps^{1+\delta_0/3})\intl_0^t|F(s)|_{V}^2\,ds+
\frac{C_{\delta_0}}{\nu\eps^{1+\delta_0/3}}\intl_0^t|v(s)|_{\mathbb{L}^{3+\delta_0}}^{2+6/\delta_0}|F(s)|_H^2\,ds.\label{eqn:NonlintermEst-2}
\end{eqnarray}
Let us choose $\eps>0$ such that
$\frac{\nu}{4}+\frac{C}{\nu}\eps^{1+\delta_0/3}=\frac{\nu}{2}$.
Then
\begin{eqnarray*}
|F(t)|_H^2+\nu\intl_0^t\Vert F(s)\Vert _{V}^2\,ds\leq
|F(0)|_H^2+\frac{C}{\nu}\intl_0^t|f(s)|_{V^\prime}^2\,ds\\
+\frac{C_{\delta_0}}{\nu\eps^{1+\delta_0/3}}\intl_0^t|v(s)|_{\mathbb{L}^{3+\delta_0}}^{2+6/\delta_0}|F(s)|_H^2\,ds,
t\geq 0.
\end{eqnarray*}
Hence, in view of the Gronwall Lemma, we get
\begin{eqnarray*}
|F(t)|_H^2\leq
\left(|F(0)|_H^2+\frac{C}{\nu}\intl_0^t|f(s)|_{V^\prime}^2\,ds\right)
e^{C(\delta_0,\nu)\intl_0^t|v(s)|_{\mathbb{L}^{3+\delta_0}}^{2+6/\delta_0}\,ds},
t\geq 0.
\end{eqnarray*}
Thus
\begin{eqnarray}
|F(t)|_H^2+\nu\intl_0^t\Vert F(s)\Vert _{V}^2\,ds\leq
K_1\left(|F(0)|_H^2+\frac{C}{\nu}\intl_0^t|f(s)|_{V^\prime}^2\,ds\right)\nonumber\\
\left(1+\intl_0^t|v(s)|_{\mathbb{L}^{3+\delta_0}}^{2+6/\delta_0}\,ds\right)
e^{C(\delta_0,\nu)\intl_0^t|v(s)|_{\mathbb{L}^{3+\delta_0}}^{2+6/\delta_0}\,ds},
t\geq 0.\label{eqn:EnergyInequality}
\end{eqnarray}
\item[(c)] The general case. Let $v_n\in
L^{\infty}(0,T;\mathbb{L}^{3+\delta_0}(D))$ be a sequence of
functions such that $v_n\to v$ in
$L^{2+\frac{6}{\delta_0}}(0,T;\mathbb{L}^{3+\delta_0}(D))$. Let
$F_n$ be a corresponding sequence of solutions of equation
(\ref{eqn:DualEquation}-\ref{in_v}) with $v$ being replaced by
$v_n$. Then from inequality \eqref{eqn:EnergyInequality} it
follows that the sequence $\{F_n\}_{n=1}^{\infty}$ lies in a
bounded set of $L^{\infty}(0,T;H)\cap L^2(0,T;V)$. Therefore, by
the Banach-Alaoglu Theorem there exists subsequence $\{F_{n'}\}$
and $F^*\in L^{\infty}(0,T;H)$ such that for any $q\in L^1(0,T;H)$
\begin{equation}
\intl_0^T(F_{n'}-F^*,q(s))_H\,ds\to 0\label{eqn:WeakConv-1}
\end{equation}
Similarly, from the Banach-Alaoglu Theorem it follows that we can
find a subsequence $\{F_{n''}\}$ of $\{F_{n'}\}$ convergent to
$F^{**}\in L^2(0,T;V)$ weakly i.e. for any $q\in
L^2(0,T;V^\prime)$
\begin{equation}
\intl_0^T\langle F_{n''}-F^{**},q(s)\rangle_{V^\prime,V}\,ds\to
0,\label{eqn:WeakConv-2}
\end{equation}
In particular, \eqref{eqn:WeakConv-1} and \eqref{eqn:WeakConv-2}
are satisfied for $q\in L^2(0,T;H)$. Therefore $F^*=F^{**}\in
L^{\infty}(0,T;H)\cap L^2(0,T;V)$. Put $F=F^*$. Let us now show
that $F$ satisfies equation (\ref{eqn:DualEquation}-\ref{in_v}) in
the weak sense. Let $\psi\in C^{\infty}([0,T],\Rnu)$, $\psi(1)=0$,
$h\in V$. Then by part $(a)$ of the proof we have
\begin{eqnarray}
-\intl_0^T(F_n(s),h)_H\psi'(s)\,ds&+&\intl_0^T\langle
B(v_n,F_n),h\rangle_{V^\prime,V}\psi(s)\,ds+
\nu\intl_0^T\tilde{a}(F_n(s),h)\psi(s)\,ds\nonumber\\
&=&(F_0,h)_H\psi(0)+\intl_0^T\langle
f(s),h\rangle_{V^\prime,V}\psi(s)\,ds.\label{eqn:WeakSolutApprox}
\end{eqnarray}
Convergence of the first term, respectively third term,
 follows immediately from \eqref{eqn:WeakConv-1}, respectively
\eqref{eqn:WeakConv-2}. For the second term we have
\begin{eqnarray*}
|\intl_0^T\langle
B(v_n,F_n)-B(v,F),h\rangle_{V^\prime,V}\psi(s)\,ds|\leq
|\intl_0^T\langle B(v_n-v,F_n),h\rangle_{V^\prime,V}\psi(s)\,ds|\\
+|\intl_0^T\langle
B(v,F_n-F),h\rangle_{V^\prime,V}\psi(s)\,ds|=I_n+II_n.
\end{eqnarray*}
Let $\eps>0$ be fixed. For any $\eps_2,\eps_3>0$ we have, by
inequality \eqref{eqn:NlinEstimate2}, the following inequalities
\begin{eqnarray*}
I_n &\leq & \eps_3\intl_0^T|\curl F_n|_H^2\,ds+
\frac{C}{\eps_3}\intl_0^T(\eps_2|h|_{V}^2+\frac{C}{\eps_2}|v_n-v|_{\mathbb{L}^{3+\delta_0}(D)}^{2+\frac{6}{\delta_0}}|h|_H^2)|\psi|^2\,ds\\
&=&\eps_3\Vert F_n\Vert
_{L^2(0,T;V)}^2+\frac{C\eps_2}{\eps_3}|h|_{V}^2\intl_0^T|\psi|^2\,ds
+\frac{C|h|_H^2}{\eps_3\eps_2}\intl_0^T|v_n-v|_{\mathbb{L}^{3+\delta_0}(D)}^{2+\frac{6}{\delta_0}}|\psi|^2\,ds.
\end{eqnarray*}
Taking into account boundedness of the sequence
$\{F_n\}_{n=1}^{\infty}$ in $L^2(0,T;V)$ and the convergence of
$\{v_n\}_{n=1}^{\infty}$ to $v$ in
$L^{2+\frac{6}{\delta_0}}(0,T;\mathbb{L}^{3+\delta_0}(D))$, we can
choose $\eps_2,\eps_3$ and $N=N(\eps)$ in such way that
$I_n\leq\frac{\eps}{2}$, for $n\geq N$.

For $II_n$ we have $II_n=|\intl_0^T\langle F_n-F,\curl (v\times
h)\rangle_{V^\prime,V}\psi(s)\,ds|$. From inequality
\eqref{eqn:EstimGag-NirType} it follows that $v\times h\in
L^2(0,T;H)$. Therefore, $\curl (v\times h)\in L^2(0,T;V^\prime)$
and the convergence of $II_n$ to $0$ follows from inequality
\eqref{eqn:WeakConv-2}. The uniqueness of $F$ follows from the
energy inequality \eqref{eqn:EnergyInequality}. It remains to show
that $F\in C([0,T],H_w)$. Let us show that $F\in
C([0,T],V^\prime)$. Then, since $F\in L^{\infty}(0,T;H)$, it
immediately follows from  \cite[Lemma 1.4, p.178]{Temam_2001} that
$F\in C([0,T],H_w)$. To prove that $F\in C([0,T],V^\prime)$ it is
enough to show that $F'\in L^1(0,T;V^\prime)$. Indeed, we have
that $F\in L^{\infty}(0,T;H)\subset L^1(0,T;V^\prime)$ and by
\cite[Lemma 1.1, p.169]{Temam_2001} the result
 follows. We have
\begin{eqnarray}
|F'|_{L^{1+\frac{3}{2\delta_0+3}}(0,T;V^\prime)}^{1+\frac{3}{2\delta_0+3}}&=&|AF|_{L^{1+\frac{3}{2\delta_0+3}}(0,T;V^\prime)}^{1+\frac{3}{2\delta_0+3}}=\intl_0^T|A(s)F(s)|_{V^\prime}^{1+\frac{3}{2\delta_0+3}}\,ds\nonumber\\
&\leq&\intl_0^T|A(s)|_{\mathcal{L}(V,V^\prime)}^{1+\frac{3}{2\delta_0+3}}|F(s)|_V^{1+\frac{3}{2\delta_0+3}}\,ds\nonumber\\
&\leq& (\intl_0^T |F(s)|_V^2\,ds)^{\frac{\delta_0+3}{2\delta_0+3}}
(\intl_0^T|A(s)|_{\mathcal{L}(V,V^\prime)}^{2+\frac{6}{\delta_0}}\,ds)^{\frac{\delta_0}{2\delta_0+3}}\nonumber\\
&\leq& (\intl_0^T |F(s)|_V^2\,ds)^{\frac{\delta_0+3}{2\delta_0+3}}
(\intl_0^T(C_1\nu+C_2|v(s)|_{\mathbb{L}^{3+\delta_0}(D)})^{2+\frac{6}{\delta_0}}\,ds)^{\frac{\delta_0}{2\delta_0+3}}\nonumber\\
&\leq&
C|F|_{L^2(0,T;V)}^{\frac{2\delta_0+6}{2\delta_0+3}}(C_1(\nu,T,\delta_0)+|v|_{L^{2+\frac{6}{\delta_0}}(0,T;\mathbb{L}^{3+\delta_0}(D))}^{\frac{2\delta_0+6}{2\delta_0+3}})<\infty,\label{eqn:TimeDerivativeEst}
\end{eqnarray}
where the second inequality follows from the H\"older inequality
and the third one follows from the inequality
\eqref{eqn:OperatorANormEst}. Thus, first statement of the
Proposition \ref{prop:ExistUniqSolution} is proved. \item[(ii)] To
prove [ii] we follow an idea from \cite{Brz_1991} and
\cite{Brz_Li_2006}.
\begin{lemma}\label{lem:LebegIntegCont}
Let $g:[0,T]\to\Rnu$ be measurable function such that
$\intl_0^T|g(s)|\,ds<\infty$. Then for any $\delta>0$ there exists
a partition $\{T_i\}_{i=1}^n$ of interval $[0,T]$ such that
$\intl_{T_i}^{T_{i+1}}|g(s)|\,ds<\delta$, $i=1,\ldots,n$.
\end{lemma}
\begin{proof} Follows easily from  \cite[Theorem 8.17]{Rudin_1987_RCA}.
\end{proof}
\begin{trivlist}
\item{Existence of a local solution.} Let $X_T=\{F \in L^2(0,T;
D(A)): F^\prime \in L^2(0,T;H)\}$ be a Banach space endowed with a
norm
$$|F|_{X_T}^2=\nu^2|F|_{L^2(0,T;D(A))}^2+|F^\prime|_{L^2(0,T;H)}^2.$$
We will prove the following result.
\begin{lemma}\label{lem_aux-1}
If $v$ satisfies assumption \eqref{eqn:IntegrCond-1}, $z\in X_T$
then $B(v(\cdot),z)\in L^2(0,T;H)$.
\end{lemma}

In view of Proposition \ref{prop:CoercivSemigroup}  and the above
Lemma, a map
 $\Phi_T:X_T\to X_T$ defined by $\Phi_T(z)=G$ iff  $G$ is  the unique solution
solution of the problem
\begin{equation}
G'+\nu AG=f-B(v(t),z),\;G(0)=F_0,\label{eqn:PhiMapDef}
\end{equation}
is well defined.
\begin{proof}[Proof of Lemma \ref{lem_aux-1}]
 From inequality
\eqref{eqn:NlinEstimate3} we have
\begin{equation*}
\Vert B(v(\cdot),z)\Vert _{L^2(0,T;H)}^2\leq
C_1(\eps,\delta_0)\Vert z\Vert _{L^2(0,T;\mathbb{H}^2(D))}^2
+C_2(\eps,\delta_0)|z|_{C([0,T];V)}^2
|v|_{L^{2+\frac{6}{\delta_0}}(0,T;\mathbb{L}^{3+\delta_0}(D))}.
\end{equation*}
Thus the result follows from Lemma \ref{lem:ContinuityLem}.
\end{proof}

We will  show that there exists  $T_1\leq T$ such that
$\Phi_{T_1}$ is a strict  contraction. By Proposition
\ref{prop:CoercivSemigroup} and inequality
\eqref{eqn:NlinEstimate3} we have, for all $t\in [0,T]$,
\begin{eqnarray*}
\Vert \Phi_t(z_1)-\Phi_t(z_2) \Vert_{X_t}^2 &\leq &C_1\Vert
B(v,z_1-z_2)\Vert _{L^2(0,t;H)}^2\leq
C_1\eps^{1+\delta_0/3}|z_1-z_2|_{L^2(0,t;D(A))}^2\\
&+&C_1\frac{C_{\delta}}{\eps^{1+3/\delta}}|z_1-z_2|_{C(0,t;V)}^2|v|_{L^{2+6/\delta_0}(0,T;\mathbb{L}^{3+\delta_0}(D))}\\
&\leq&  C_1\eps^{1+\delta_0/3}|z_1-z_2|_{X_t}^2\\
&+&
C_1\frac{C_{\delta}}{\eps^{1+3/\delta}}|z_1-z_2|_{X_t}^2|v|_{L^{2+6/\delta_0}(0,t;\mathbb{L}^{3+\delta_0}(D))}).
\end{eqnarray*}
Now let us choose $\eps>0$ that $C_1\eps^{1+\delta_0/3}=1/2$ and
denote $K=C_1\frac{C_{\delta}}{\eps^{1+3/\delta}}$. We have
\begin{equation}
\Vert \Phi_t(z_1)-\Phi_t(z_2)\Vert
_{X_t}^2\leq(1/2+K|v|_{L^{2+6/\delta_0}(0,t;\mathbb{L}^{3+\delta_0}(D))})|z_1-z_2|_{X_t}^2,
t\in [0,T].\label{eqn:ContractEst-1}
\end{equation}
Choose $t=T_1$ such that
$|v|_{L^{2+6/\delta_0}(0,T_1;\mathbb{L}^{3+\delta_0}(D))}\leq
d=\frac{1}{3K}$ then $\Phi_{T_1}$ is an affine contraction map and
by the Banach Fixed Point Theorem there exists a fixed point $F\in
X_{T_1}$ of $\Phi_{T_1}$. Obviously $F$ is a solution of problem
(\ref{eqn:DualEquation}-\ref{in_v}) on interval $[0,T_1]$.
\item{Existence of a global solution.} From Lemma
\ref{lem:LebegIntegCont} and assumption \eqref{eqn:IntegrCond-1}
it follows that we can find partition
$0=T_0<T_1<\ldots<T_{k-1}<T_k=T$ of interval $[0,T]$ such that
$|v|_{L^{2+6/\delta_0}(T_i,T_{i+1};\mathbb{L}^{3+\delta_0}(D))}<\frac1{3K}$,
$i=0,\ldots,k-1$. Therefore, we can use the inequality
\eqref{eqn:ContractEst-1} and the Banach Fixed Point Theorem
iteratively to define global solution.
\end{trivlist}
\item[(iii)] To proof the statement in part  [iii] we will use a
method suggested by Temam in \cite{Temam_1981}. We will consider
only the case $k=1$. General case follows by induction. Let us
recall that
$$A(t)=\nu A+B(v(t),\cdot).$$ By differentiating the equation
(\ref{eqn:DualEquation}-\ref{in_v}) w.r.t. $t$ (in weak sense) we
find that $F^\prime$ is a solution of
\begin{equation*}
\frac{dF^\prime}{dt}=-A(t)F^\prime+B(v^\prime(t),F)+f^\prime, t\in
[0,T].
\end{equation*}
Now from the assumptions of the statement in part [ii] it follows
that it is enough to prove that $B(v^\prime(\cdot),F)\in
L^2(0,T;H)$ and then use the already proven statement in part [i].
From inequality \eqref{eqn:NlinEstimate3} we have
\begin{eqnarray*}
\intl_0^T|B( v^\prime(t),F)|_H^2\,dt &\leq&
\eps^{1+\delta_0/3}\intl_0^T\Vert \curl F\Vert
_V^2\,dt+\frac{C_{\delta_0}}{\eps^{1+3/\delta_0}}
\intl_0^T|v^\prime(t)|_{\mathbb{L}^{3+\delta_0}(D)}^2|\curl
F|_H^2\,dt
\\
&\leq & \eps^{1+\delta_0/3}|F|_{L^2(0,T;D(A))}
+\frac{C_{\delta_0}}{\eps^{1+3/\delta_0}}\Vert F\Vert
_{C([0,T];V)}^2\intl_0^T
|v^\prime(t)|_{\mathbb{L}^{3+\delta_0}(D)}^{2+\frac{6}{\delta_0}}\,dt<\infty.
\end{eqnarray*}
Note that $F\in C([0,T];V)$ by Lemma \ref{lem:ContinuityLem}.
\end{trivlist}
\end{proof}
\begin{proof}[Proof of Proposition \ref{prop:ExistUniqSolution-2}]
The proof is very similar to the proof of the previous
Proposition.
\begin{trivlist}
\item[(i)] The  proof will be divided into three parts a), b), c).
\item[\textbf{a)}] First we consider a special case when $v\in
L^{\infty}(0,T;\mathbb{L}^{3+\delta_0}(D))$. We will use Theorem
\ref{thm:LionsPDEexun} with Gelfand triple $V\subset H\cong
H'\subset V^\prime$. Denote $B(t)=\nu A+\curl(v(t)\times\cdot)$.%,$D(B(t))=\Red{?}$.
We need to check whether the conditions \eqref{eqn:BoundednCond}
and \eqref{eqn:CoercivCond} are satisfied. We have
\begin{eqnarray}
\langle B(t)f,f\rangle_{V^\prime,V}&=&\nu\tilde{a}(f,f)+\langle
\curl(v(t)\times
f),f\rangle_{V^\prime,V}\nonumber\\
&=&\nu \tilde{a}(f,f)+\langle v(t)\times f,\curl
f\rangle_{V^\prime,V},\;t\in[0,T],\; f\in V.\label{eqn:QuadForm-2}
\end{eqnarray}
Now we can use the inequality \eqref{eqn:NlinEst-1} and the
continuity of the form $\tilde{a}$ to get
$$\Vert B(t)\Vert _{\mathcal{L}(V,V^\prime)}\leq C
\nu+C_2|v(t)|_{\mathbb{L}^{3+\delta_0}(D)}.$$  The coercivity
condition \eqref{eqn:CoercivCond} can be proved in the same way as
in the proof of Proposition \ref{prop:ExistUniqSolution}.
Therefore, by Theorem \ref{thm:LionsPDEexun} first statement of
the Proposition is proved in  our special case. \item[\textbf{b)}]
To prove Proposition in the general case we will, as before, show
an energy inequality for solutions of the problem
(\ref{eqn:DualEquation-2}-\ref{in_v-2}) when $v\in
L^{\infty}(0,T;\mathbb{L}^{3+\delta_0}(D))$. From Step
\textbf{(a)} we know that there exists a unique solution $G\in
L^2(0,T;V)$ such that $G^\prime\in L^2(0,T;V^\prime)$. Then, from
Lemma \ref{lem:ContinuityLem} it follows that $G\in C([0,T];H)$
and
%Multiplying (\ref{eqn:DualEquation-2}) on $G$ and integrating we get
\begin{eqnarray*}
\frac{1}{2}\frac{d}{dt}|G|_H^2&=&-\nu\Vert G\Vert _{V}^2+\langle
f,G\rangle_{V^\prime,V}-\langle v\times
G,\curl G\rangle_{V^\prime,V}\\
&=&-\nu\Vert G\Vert _{V}^2+\langle f,G\rangle_{V^\prime,V}+(\curl
G,v\times G)_H
\end{eqnarray*}
Therefore, by the Young inequality,
\begin{eqnarray*}
|G(t)|_H^2&+&2\nu\intl_0^t|G(s)|_{V}^2\,ds-\intl_0^t(\curl
G(s),v(s)\times G(s))_H\,ds\\
&=&|G(0)|_H^2+
\intl_0^t\langle f(s),G(s)\rangle_{V^\prime,V}\,ds\\
&\leq& |G(0)|_H^2+
\frac{\nu}{2}\intl_0^t|G(s)|_{V}^2\,ds+\frac{C}{\nu}\intl_0^t|f(s)|_{V^\prime}^2\,ds,
t\in [0,T].
\end{eqnarray*}
The term $\intl_0^t(\curl G(s),v(s)\times G(s))_H\,ds$ can be
estimated in the same way as in Proposition
\ref{prop:ExistUniqSolution}, see \eqref{eqn:NonlintermEst-2}.
Thus we infer that $G$ satisfies the following inequality, for
$t\in [0,T]$,
\begin{eqnarray}
\label{eqn:EnergyInequality-2Aux} |G(t)|_H^2&+&\nu\intl_0^t\Vert
G(s)\Vert _{V}^2\,ds\leq
K_1\left(|G(0)|_H^2+\frac{C}{\nu}\intl_0^t|f(s)|_{V^\prime}^2\,ds\right)\nonumber\\
&&\hspace{2truecm}\left(1+\intl_0^t|v(s)|_{\mathbb{L}^{3+\delta_0}}^{2+6/\delta_0}\,ds\right)\,
e^{C(\delta_0,\nu)\intl_0^t|v(s)|_{\mathbb{L}^{3+\delta_0}}^{2+6/\delta_0}\,ds}.
\end{eqnarray}
\item[\textbf{c)}] The general case. Now, let
$\{v_n\}_{n=1}^{\infty}$ be an
$L^{\infty}(0,T;\mathbb{L}^{3+\delta_0}(D))$-valued sequence of
functions such that $v_n\to v\in
L^{2+\frac{6}{\delta_0}}(0,T;\mathbb{L}^{3+\delta_0}(D)),n\to\infty$
in $L^{2+\frac{6}{\delta_0}}(0,T;\mathbb{L}^{3+\delta_0}(D))$. Let
$\{G_n\}_{n=1}^{\infty}$ be corresponding sequence of solutions of
the problem (\ref{eqn:DualEquation-2}-\ref{in_v-2}). Then from
\eqref{eqn:EnergyInequality-2Aux} it follows that sequence
$\{G_n\}_{n=1}^{\infty}$ lie in a bounded set of
$L^{\infty}(0,T;H)\cap L^2(0,T;V)$. Using the same argument as in
the proof of Proposition \ref{prop:ExistUniqSolution} we can find
subsequence $\{G_{n'}\}_{n'=1}^{\infty}$ weakly convergent to
$G\in L^{\infty}(0,T;H)\cap L^2(0,T;V)$ which solves the problem
(\ref{eqn:DualEquation-2}-\ref{in_v-2}) in a weak sense. Moreover,
it follows from inequality \eqref{eqn:EnergyInequality-2Aux}, that
the function $G$ satisfies energy inequality
\eqref{eqn:EnergyInequality-2}. Uniqueness of the solution of the
problem
(\ref{eqn:DualEquation-2}-\ref{in_v-2}) %\Red{of which problem ?}
follows from the energy inequality \eqref{eqn:EnergyInequality-2}.
The only difference with the previous Proposition is that now we
can prove that $G'\in L^2(0,T,V^\prime)$. Indeed, we have
\begin{eqnarray*}
\Vert G'\Vert _{L^2(0,T,V^\prime)}^2 &=& \Vert BG\Vert
_{L^2(0,T,V^\prime)}^2\\
&\leq&\intl_0^T|\nu AG+\curl(v(t)\times
G(t))|_{V^\prime}^2\,dt\\
&\leq&\nu^2\Vert G\Vert _{L^2(0,T,V)}^2+\intl_0^T|v(t)\times
G(t)|_H^2\,dt\\
&\leq& \nu^2\Vert G\Vert _{L^2(0,T,V)}^2
+\intl_0^T(C_1|G(t)|_V^2+C_2|v(t)|_{\mathbb{L}^{3+\delta_0}}^{2+6/\delta_0}|G(t)|_H^2)\,dt\\
&\leq& C_3\Vert G\Vert _{L^2(0,T,V)}^2+C_2\Vert G\Vert
_{L^{\infty}(0,T,H)}^2|v|_{L^{2+\frac{6}{\delta_0}}(0,T;\mathbb{L}^{3+\delta_0}(D))}<\infty.
\end{eqnarray*}
Thus, the first statement of Proposition is proved. Statements
[ii] and [iii] can be proved in the same way as in the proof of
Proposition \ref{prop:ExistUniqSolution}.
%%% For dissertation
%\Blue{A reader can be lost here. What are we proving now? what
%does the (i) below refer to? There was  a mistake here. I changed
%(i) to (ii) and (ii) to (iii) below.}

\item[(ii)]
\begin{trivlist}
\item{\textit{Existence of a local solution.}} Let $X_T=\{F \in
L^2(0,T; D(A)): F^\prime \in L^2(0,T;H)\}$ be a Banach space
endowed with a norm
$$|F|_{X_T}^2=\nu^2|F|_{L^2(0,T;D(A))}^2+|F^\prime|_{L^2(0,T;H)}^2.$$
We will prove the following result.
\begin{lemma}\label{lem_aux-2}
If  $v\in L^2(0,T;V)$, $z\in X_T$ then $\curl(v(t)\times z)\in
L^2(0,T;H)$. %\Red{What is $z$ here?Changed.}
\end{lemma}
In view of Proposition \ref{prop:CoercivSemigroup}  and the above
Lemma, a map $\Phi_T:X_T\to X_T$ defined by $\Phi_T(z)=G$ iff  $G$
is  the unique solution of the problem
\begin{equation}
G'+\nu AG=f-\curl(v(t)\times z),G(0)=F_0\in
V\label{eqn:PhiMapDef-2}
\end{equation}
is well defined.
\begin{proof}[Proof of Lemma \ref{lem_aux-2}]
We have:
\begin{eqnarray}
\Vert \curl(v(t)\times z)\Vert _{L^2(0,T;H)}^2&\leq& C(\Vert
z\nabla
v\Vert _{L^2(0,T;H)}^2+\Vert v\nabla z\Vert _{L^2(0,T;H)}^2)\nonumber\\
&\leq&
C|z|_{C([0,T];V)}|v|_{L^2(0,T;V)}\label{eqn:NonlinearTermEst}
\end{eqnarray}
and the result follows from Lemma \ref{lem:ContinuityLem}.
%\Red{You mean that we apply Lemma \ref{lem:ContinuityLem} to
%function $z$? But what are the assumptions about $z$?Changed. See
%above.}
\end{proof}
Now we will show that there exists  $T_1\in (0, T]$ such that
$\Phi_{T_1}$ is a strict contraction. Let us fix $t\in (0,T]$ and
take $z_1,z_2 \in X_T$. Then,  by Proposition
\ref{prop:CoercivSemigroup} and Lemma \ref{lem:EstimGag-NirType},
we have
\begin{eqnarray*}
\Vert \Phi_t(z_1)-\Phi_t(z_2)\Vert _{X_t}^2&\leq& C_1\Vert \curl(v(t)\times (z_1-z_2))\Vert _{L^2(0,t;H)}^2\\
&\leq& C|z_1-z_2|_{C(0,t;V)}^2|v|_{L^2(0,t;V)}^2\leq
C|z_1-z_2|_{X_t}^2|v|_{L^2(0,t;V)}^2.
\end{eqnarray*}
 Let us choose $T_1\in  (0,T]$ such that $C|v|_{L^2(0,T_1;V)}<1/2$. Then
$\Phi_{T_1}$ is a strict contraction map and hence by the Banach
Fixed Point Theorem there exists a unique  $F\in X_{T_1}$ that is
a fixed point of $\Phi_{T_1}$. By the definition of the $\Phi_T$
it follows that  $F\in X_{T_1}$ is a solution of problem
(\ref{eqn:DualEquation}-\ref{in_v}) on interval $[0,T_1]$. Notice
also that $F(T_1)\in V$. Therefore, the map $\Phi_T$ with initial
data $F(T_1)$ is well defined on interval $[T_1,T]$.

\item{\textit{ Existence of a global solution.}} From Lemma
\ref{lem:LebegIntegCont} and assumption \eqref{eqn:IntegrCond-1}
it follows that we can find a partition
$0=T_0<T_1<\ldots<T_{k-1}<T_k=T$ of the interval $[0,T]$ such that
$|v|_{L^2(T_i,T_{i+1};V)}<1/2$, $i=0,\ldots,k-1$. Therefore, we
can use  inequality \eqref{eqn:ContractEst-1} and the Banach Fixed
Point Theorem iteratively to define a global solution
(\ref{eqn:DualEquation-2}-\ref{in_v-2}). %\Red{Have you mentioned
%that $F(T_1)$ is a well defined element of the space $?$ and so
%the map $\Phi_T$ on the time interval $[T_1,T_2]$  (which differs
%from the one on $[0,T_1]$ by the choice of the initial data) is
%well defined?Changed.}
\end{trivlist}
\item[(iii)]  We will consider only the case $k=1$. General case
follows by induction. We differentiate equation
(\ref{eqn:DualEquation-2}-\ref{in_v-2}) w.r.t. $t$ (in the weak
sense) %\Red{Misha, it seems in almost every sentence  a word "a"
%or "the" is missing. Something has to be done about this!}
and get
an equation for the function $G^\prime$:
\begin{eqnarray*}
\frac{\partial G^\prime}{\partial t}(t) &=& -\nu
AG^\prime(t)-\curl{G^\prime(t)\times v(t)}-\curl(G(t)\times
v^\prime(t))+f^\prime(t)\\
G^\prime(0) &=& -\nu A G_0-\curl{v(0)\times G_0}+f(0),t\in[0,T].
\end{eqnarray*}
Now from the assumptions of the statement [ii] it follows that it
is enough to prove that $\curl(G\times v^\prime(t))\in L^2(0,T;H)$
and then use the already proven statement in part [i]. By the
inequality \eqref{eqn:NonlinearTermEst}  we have
\begin{equation*}
\Vert \curl(G\times v^\prime)\Vert _{L^2(0,T;H)}^2\leq
C|G|_{C([0,T];V)}|v^\prime|_{L^2(0,T;V)}<\infty.
\end{equation*}
Note that $G\in C([0,T];V)$ by Lemma \ref{lem:ContinuityLem}.
%\Red{What is the logical sequence in the above argument?Changed.}

%%% for dissertation
\end{trivlist}
\end{proof}
\begin{proof}[Proof of Theorem \ref{thm:Duality}]
\begin{trivlist}
\item{$1^{st}$ Step.} %\Red{Was unable to find normal version of
%different "first".}
Fix $\delta_0>0$. Let us prove the theorem in the case of smooth
initial data $F_0\in C^{\infty}(\overline{D})\cap H$ and vector
field $v\in C_b^{\infty}([0,T]\times\overline{D})\cap
L^{2+\frac{6}{\delta_0}}(0,T;\mathbb{L}^{3+\delta_0}(D))$.
%\Red{Can we be more specific about the assumptions on $v$. Also we
%should mention that we fix now $\delta_0$.Changed.}
For each $\eps>0$ we can find $F_0^{\eps}\in
C^{\infty}(\overline{D})\cap H$, $G_0^{\eps}\in
C^{\infty}(\overline{D})\cap H$, $v^{\eps}\in
C_b^{\infty}([0,T]\times\overline{D})\cap
L^{2+\frac{6}{\delta_0}}(0,T;\mathbb{L}^{3+\delta_0}(D))$ such
that as $\eps\to 0$, $F_0^{\eps}\to F_0$  in $H$, $G_0^{\eps}\to
G_0$
 in $H$ and $v^{\eps}\to v$ in
$L^{\infty}(0,T;\mathbb{L}^{3+\delta_0}(D))$. It follows from
Corollaries \ref{cor:InfReg-1} and \ref{cor:InfReg-2} that there
exists solutions $F^{\eps}\in C([0,T];H)\cap
C^{\infty}((0,T]\times\overline{D})$, $G^{\eps}\in C([0,T];H)\cap
C^{\infty}((0,T]\times\overline{D})$) that are solutions to the
following problems
\begin{eqnarray*}
\frac{\partial F^{\eps}}{\partial t}(t) &=& -\nu A F^{\eps}(t)-P(v^{\eps}(t)\times\curl F^{\eps}(t))\\
F^{\eps}(0,\cdot) &=& F_0^{\eps},t\in [0,T]\\
\frac{\partial G^{\eps}}{\partial t}(t) &=& -\nu A G^{\eps}(t)+\curl{(v^{\eps}(T-t)\times G^{\eps}(t))}\\
G^{\eps}(0,\cdot) &=& G_0^{\eps},t\in [0,T]
\end{eqnarray*}
Therefore, for all $t\in (0,T]$ we have %\Red{Please add letter "t"
%in those equations (and below as well).Changed.}
\begin{eqnarray*}
&&\frac{d}{dt}(F^{\eps}(t),G^{\eps}(T-t))_{\mathbb{L}^2(D)}\\
&=& (\frac{d}{dt}F^{\eps}(t),G^{\eps}(T-t))_{\mathbb{L}^2(D)}-
(F^{\eps}(t),\frac{d}{dt}G^{\eps}(T-t))_{\mathbb{L}^2(D)}\\
&=&\nu(\mathrm{P}\triangle
F^{\eps}(t),G^{\eps}(T-t))_{\mathbb{L}^2(D)}
-(\mathrm{P}(v(t)\times\curl F^{\eps}(t)),G^{\eps}(T-t))_{\mathbb{L}^2(D)}\\
&-& \nu(F^{\eps}(t),\mathrm{P}\triangle
G^{\eps}(T-t))_{\mathbb{L}^2(D)}-(F^{\eps}(t),\curl{(v(t)\times
G^{\eps}(T-t))})_{\mathbb{L}^2(D)}\\
&=&K_1(t)-K_2(t)-K_3(t)-K_4(t)
\end{eqnarray*}
It follows from the fact that
$\diver{F^{\eps}}=\diver{G^{\eps}}=0$, $F^{\eps}|_{\partial
D}=G^{\eps}|_{\partial D}=0$ and the integration by parts formula
that
$(F^{\eps},\nabla\psi)_{\mathbb{L}^2(D)}=(G^{\eps},\nabla\psi)_{\mathbb{L}^2(D)}=0$
for any $\psi\in C^{\infty}(\overline{D})$. Thus, we have
\begin{eqnarray}
K_1(t) &=&(\mathrm{P}\triangle
F^{\eps}(t),G^{\eps}(T-t))_{\mathbb{L}^2(D)}=(\triangle
F^{\eps}(t),G^{\eps}(T-t))_{\mathbb{L}^2(D)},\nonumber\\
K_2(t) &=&(\mathrm{P}(v(t)\times\curl
F^{\eps}(t)),G^{\eps}(T-t))_{\mathbb{L}^2(D)}\nonumber\\ &=&
(v(t)\times\curl
F^{\eps}(t),G^{\eps}(T-t))_{\mathbb{L}^2(D)},t\in(0,T]\label{eqn:ii-equality}
\end{eqnarray}
and
$$
K_3(t)=(F^{\eps}(t),\mathrm{P}\triangle
G^{\eps}(T-t))_{\mathbb{L}^2(D)}=(F^{\eps}(t),\triangle
G^{\eps}(T-t))_{\mathbb{L}^2(D)},t\in(0,T]
$$
Therefore, by the Green Formula we get
$K_1(t)-K_3(t)=0,t\in(0,T]$. From
\eqref{eqn:Wedge-1},\,\,\eqref{eqn:ii-equality} and the formula
\begin{equation*}
\intl_Du\curl vdx-\intl_Dv\curl udx=\intl_{\partial D}(u\times
v,\overrightarrow{n})\,d\sigma
\end{equation*}
we infer that
\begin{eqnarray*}
K_2(t) &=&(v(t)\times\curl
F^{\eps}(t),G^{\eps}(T-t))_{\mathbb{L}^2(D)}\\
&=&-(\curl F^{\eps}(t)\times v(t),G^{\eps}(T-t))_{\mathbb{L}^2(D)}\\
&=&-(\curl F^{\eps}(t),v(t)\times
G^{\eps}(T-t))_{\mathbb{L}^2(D)}=-K_4(t),t\in(0,T].
\end{eqnarray*}
Thus,
$\frac{d}{dt}(F^{\eps}(t),G^{\eps}(T-t))_{\mathbb{L}^2(D)}=0,t\in(0,T]$.
Also, by the regularity of $F^{\eps}$, $G^{\eps}$ it follows that
$(F^{\eps}(\cdot),G^{\eps}(T-\cdot))_{\mathbb{L}^2(D)}\in
C^{\infty}((0,T])\cap C([0,T])$. As a result we get equality
\eqref{eqn:Duality}.

\item{$2^{nd}$ step.} \del{Let us} Let us suppose that we have
showed that %\Red{What quantifier we apply to $t$ here?Changed.}
$F_{\eps}(t)\to F(t),t\in [0,T]$ in weak topology of $H$ and
$G_{\eps}\to G$ in $C([0,T],H)$. Then we have
\begin{eqnarray*}
|(F(t),G(T-t))&-&(F^{\eps}(t),G^{\eps}(T-t))|\\
&=&|(F-F^{\eps}(t),G(T-t))+(F^{\eps}(t),G-G^{\eps}(T-t))|\\
&\leq&|(F-F^{\eps}(t),G(T-t))|+|F^{\eps}(t)|_H|G-G^{\eps}(T-t)|_H\\
&\leq&|(F-F^{\eps}(t),G(T-t))|+|F_0^{\eps}|_H\supl_{s\in
[0,T]}|G-G^{\eps}(s)|_H\stackrel{\eps\to 0}{\rightarrow}0,t\in
[0,T]
\end{eqnarray*}
i.e. $(F(t),G(T-t))_{H}=\liml_{\eps\to
0}(F^{\eps}(t),G^{\eps}(T-t))_{H},t\in [0,T]$ and the result
follows from first step. In order to show weak convergence of
$F_{\eps}(t)$ to $F(t)$, $t\in [0,T]$ let us first notice that by
the Banach-Alaoglu Theorem, $F_{\eps}$ converges to $F$ weakly-*
in $L^{\infty}(0,T;H)$. The proof of  this claim can be performed
in exactly the same as the proof of the convergence of $F_n \to F$
in Proposition \ref{prop:ExistUniqSolution}). Also, we have from
the Banach-Alaoglu Theorem that $F^{\eps}(t)$ weakly-* convergent
to some $\Psi(t)\in H,t\in [0,T]$. We will show that $\Psi=F$. Fix
$\xi\in V$. Let us denote $g(t)=(\Psi(t)-F(t),\xi)_H$, $t\in
[0,T]$. Since $V$ is dense in $H$ it is enough to show that $g=0$.
Now we will show that $g\in C([0,T])$. From  the part (i) of
Proposition \ref{prop:ExistUniqSolution} we infer that $F\in
C([0,T],V^\prime)$. Thus, $(F(\cdot),\xi)_H=\lb
F(\cdot),\xi\rb_{V^\prime,V}\in C([0,T])$. Furthermore, for $t\in
[0,T]$ we have
\begin{eqnarray}
&&|(F_{\eps}(t),\xi)-(F_{\eps}(s),\xi)|\leq \intl_s^t|\lb F_{\eps}'(r),\xi\rb |dr\label{eqn:HolderContinuity-1}\\
&\leq &(\intl_0^T |\lb F_{\eps}'(r),\xi\rb
|^{1+\frac{3}{2\delta_0+3}}\,dr)^{\frac{2\delta_0+3}{2\delta_0+6}}|t-s|^{\frac{3}{2\delta_0+6}}
\nonumber\\
&\leq& |F_{\eps}'|_{L^{1+\frac{3}{2\delta_0+3}}(0,T;V^\prime)}|\xi|_V|t-s|^{\frac{3}{2\delta_0+6}}\nonumber\\
&\leq &
C|F_{\eps}|_{L^2(0,T;V)}(C_1(\nu,T,\delta_0)+|v_{\eps}|_{L^{2+\frac{6}{\delta_0}}(0,T;\mathbb{L}^{3+\delta_0}(D))}^{\frac{2\delta_0+6}{2\delta_0+3}})^{\frac{2\delta_0+3}{2\delta_0+6}}|\xi|_V|t-s|^{\frac{3}{2\delta_0+6}}\nonumber\\
&\leq & C|F_{\eps}^0|_H(C(\nu,T,\delta_0)+|v_{\eps}|_{L^{2+\frac{6}{\delta_0}}(0,T;\mathbb{L}^{3+\delta_0}(D))})|\xi|_V|t-s|^{\frac{3}{2\delta_0+6}}\nonumber\\
&\leq&
C|F^0|_H(C(\nu,T,\delta_0)+|v|_{L^{2+\frac{6}{\delta_0}}(0,T;\mathbb{L}^{3+\delta_0}(D))})|\xi|_V|t-s|^{\frac{3}{2\delta_0+6}},0<s\leq
t<T.\nonumber
\end{eqnarray}
In the above sequence of inequalities,  the first one  follows
because $(F_{\eps}(\cdot),\xi)\in C^{\infty}((0,T))$, the second
one from the H\"older inequality and the fourth one  from the
inequality \eqref{eqn:TimeDerivativeEst}.

Taking the $\eps \todown 0$ limit  in
\eqref{eqn:HolderContinuity-1} we immediately get
\begin{equation}
|(\Psi(t),\xi)-(\Psi(s),\xi)|\leq
C(F_0,v,\nu,\delta_0,T)|\xi|_V|t-s|^{\frac{3}{2\delta_0+6}},\label{eqn:HolderContinuity-2}
\end{equation}
where
$C(F_0,v,\nu,\delta_0,T)=C|F^0|_H(C(\nu,T,\delta_0)+|v|_{L^{2+\frac{6}{\delta_0}}(0,T;\mathbb{L}^{3+\delta_0}(D))})$. %\Red{Please
%finish.Changed.}
Hence, $\Psi\in C([0,T],V^\prime)$ and,
consequently, $g\in C([0,T])$. Therefore it is enough to prove
that $g(t)=0$ for a.a. $t\in [0,T]$. We have already observed that
\begin{equation}
\lim_{\eps\todown 0} \intl_0^T(F^{\eps}(s)-F(s),q(s))_H\,ds =0,
\,\mbox{for all } \,q\in L^1(0,T;H).\label{eqn:AuxWeakConv-1}
\end{equation}
Take any $f\in L^1(0,T)$ and put $q=\xi f$,
$g_{\eps}=(F^{\eps}(\cdot)-F(\cdot),\xi)_H$. Then  condition
\eqref{eqn:AuxWeakConv-1} can be rewritten as follows
\begin{equation}
\intl_0^Tg_{\eps}(s)f(s)\,ds\to 0,\,\mbox{for all } \, f\in
L^1(0,T).\label{eqn:AuxWeakConv-2}
\end{equation}
On the other hand, it follows from definition of $g$ that
$g_{\eps}$ is convergent to $g$ pointwise. Let us show that
\eqref{eqn:AuxWeakConv-2} and pointwise convergence of $g_{\eps}$
imply that $g=0$ a.e.. By the Egorov Theorem, see e.g.
 \cite[Theorem 2.2.1, p. 110]{[Bogachev-2007]},  for any $l>0$ there exists
a measurable set $A_l\subset [0,T]$ such that $\la(A_l)<l$  and
$g_{\eps}\to g$  uniformly  on $[0,T]\setminus A_l$. Here  $\la$
denotes the Lebesgue measure. Hence by  \eqref{eqn:AuxWeakConv-2}
we infer that $g(t)=0$, for a.e. $t\in[0,T]\setminus A_l$ \del{
Hence, $\la(\{t:g(t)\neq 0\})<l$ for any $l>0$.} and
consequently\del{, we infer that}  $g(t)=0$ for a.e. $t\in[0,T]$.

Thus, it remains to show that $G_{\eps}\to G$ in $C([0,T],H)$.
Denote $R^{\eps}=G^{\eps}-G$. Then $R^\eps$ is a solution to the
following problem.
\begin{eqnarray*}
\frac{\partial R^{\eps}}{\partial t}(t) &=& -\nu A
R^{\eps}(t)+\curl(v^{\eps}(T-t)\times R^{\eps}(t))
+\curl((v^{\eps}(T-t)-v(T-t))\times G(t))\\
R^{\eps}(0,\cdot) &=& G_0^{\eps}-G_0, t\in [0,T].
\end{eqnarray*}
Applying  the energy inequality \eqref{eqn:EnergyInequality-2} to
the function $R^\eps$ we infer that for any $\tau>0$
\begin{eqnarray}
|R^{\eps}|_{C([0,T];H)}^2&\leq&
C(|v^{\eps}|_{L^{2+\frac{6}{\delta_0}}(0,T;\mathbb{L}^{3+\delta_0}(D))})
(|G_0^{\eps}-G_0|_{H}^2+|\curl((v-v^{\eps})\times
G)|_{L^2(0,T;V^\prime)}^2)\nonumber\\
&\leq&
C(|v^{\eps}|_{L^{2+\frac{6}{\delta_0}}(0,T;\mathbb{L}^{3+\delta_0}(D))})
(|G_0^{\eps}-G_0|_{H}^2+|(v-v^{\eps})\times G|_{L^2(0,T;H)}^2)\nonumber\\
&\leq&
C(|v|_{L^{2+\frac{6}{\delta_0}}(0,T;\mathbb{L}^{3+\delta_0}(D))})(|G_0^{\eps}-G_0|_{H}^2+\tau^{1+\delta_0/3}\intl_0^T|G(s)|_V^2\,ds+\nonumber\\
&&\frac{C_{\delta_0}}{\tau^{1+3/\delta_0}}|G|_{C([0,T];H)}^2|v^{\eps}-v|_{L^{2+\frac{6}{\delta_0}}(0,T;\mathbb{L}^{3+\delta_0}(D))}),\label{eqn:residueEst-1}
\end{eqnarray}
where last inequality of \eqref{eqn:residueEst-1} follows from
Lemma \ref{lem:EstimGag-NirType}. Now, from the convergences
$v^{\eps} \to v$ in
$L^{2+\frac{6}{\delta_0}}(0,T;\mathbb{L}^{3+\delta_0}(D))$,
$G_0^{\eps} \to G_0$ in $H$ and inequalities
\eqref{eqn:residueEst-1} we get the result.
\end{trivlist}
\end{proof}

\end{document}